\documentclass{article}
\usepackage{TeX_Math_Preamble}
\title{Bridging persistence of discrete and continuous 1D maps: Total variation, permutation equivalence, and diagonal affine transformations}
\author{Matthew Blake\thanks{blake27@msu.edu.}, Firas A.~Khasawneh\thanks{khasawn3@msu.edu.}}

\date{Department of Computational Mathematics, Sciences and Engineering; Michigan State University, East Lansing, MI 48824}

\begin{document}
\maketitle

\begin{abstract}
    Fundamental to the field of topological signal processing (TSP) is the \textit{persistence diagram} (PD), or barcode, of a signal, which gives noise-resistant measures of topologically significant features as distinct from noise. Existing literature proves that the total persistence computed from these diagrams is equal to the total variation, for certain continuous signals on compact domains. We provide a closed formula for this shared quantity, and generalize to all continuous and discrete functions on bounded domains, and to certain functions on non-compact 1D manifolds. Extending the identification of continuous and discrete persistence, we investigate the similarities of functions with the same PDs, and prove that, for any $n$, the persistence diagram can be calculated as a function of the sequence of critical values; we compute these formulas for $n \leq 6$. Finally, we examine which function transformations cause predictable changes in the PD of the transformed function, characterizing affine and linear transformations. %
\end{abstract}

\begin{multicols}{2}
\begin{description}
    \item[\textbf{\ref*{sec:Introduction}}]  \nameref{sec:Introduction} (\pageref*{sec:Introduction})
    \item[\textbf{\ref*{sec:Contributions}}] \nameref{sec:Contributions}
    (\pageref*{sec:Contributions})
    \item[\ref*{sec:Domains, Level Sets, and Persistence}] \nameref{sec:Domains, Level Sets, and Persistence} (\pageref*{sec:Domains, Level Sets, and Persistence})
    \begin{description}
        \item[\ref*{subsec:Discrete and Continuous Maps}] \nameref{subsec:Discrete and Continuous Maps}
        \item[\ref*{subsec:Domains, Genericity, and Critical Points}] \nameref{subsec:Domains, Genericity, and Critical Points}
        \item[\ref*{subsec:Level sets}] \nameref{subsec:Level sets}
        \item[\ref*{subsec:Persistence Diagrams}] \nameref{subsec:Persistence Diagrams}
    \end{description}
    
    \item[\ref*{sec:Total Variation and Total Persistence}] \nameref{sec:Total Variation and Total Persistence} (\pageref*{sec:Total Variation and Total Persistence})
    \begin{description}
        \item[\ref*{subsec:Computation and Comparison}] \nameref{subsec:Computation and Comparison} 
        \item[\ref*{subsec:Formula Derivation}] \nameref{subsec:Formula Derivation}
        \item[\ref*{subsec:Circular Domains}] \nameref{subsec:Circular Domains}
        \item[\ref*{subsec:Non-Separating Functions}] \nameref{subsec:Non-Separating Functions}
        \item[\ref*{subsec:Discrete Functions}] \nameref{subsec:Discrete Functions}
        \item[\ref*{subsec:Non-Tame Functions}]
        \nameref{subsec:Non-Tame Functions}
    \end{description}
    \item[\ref*{sec:Permutation Equivalence}] \nameref{sec:Permutation Equivalence} (\pageref*{sec:Permutation Equivalence})
    \begin{description}
        \item[\ref*{subsec:Permutation Theorem}]
        \nameref{subsec:Permutation Theorem}
        \item[\ref*{subsec:Persistence Diagram Computation}] \nameref{subsec:Persistence Diagram Computation}
        \item[\ref*{subsec:Persistence_Diagram_Preimage}] \nameref{subsec:Persistence_Diagram_Preimage}
    \end{description}

    \item[\ref*{sec:Signal Transformations}] \nameref{sec:Signal Transformations} (\pageref*{sec:Signal Transformations})
    \begin{description}
        \item[\ref*{subsec:Diagonal Affine Transformation Theorem}] \nameref{subsec:Diagonal Affine Transformation Theorem}
        \item[\ref*{subsec:Function Addition}] \nameref{subsec:Function Addition}
        \item[\ref*{subsec:Integration and Differentiation}] \nameref{subsec:Integration and Differentiation}
    \end{description}

    \item[\ref*{sec:Conclusions and Further Directions}] \nameref{sec:Conclusions and Further Directions} (\pageref*{sec:Conclusions and Further Directions})
    \item[\ref*{sec:app_Lemmas}] \nameref{sec:app_Lemmas} (\pageref*{sec:app_Lemmas})
    \item[\ref*{sec:app_Generalizing Equivalence Domain}]\nameref{sec:app_Generalizing Equivalence Domain} (\pageref*{sec:app_Generalizing Equivalence Domain})
\begin{description}
    \item[\ref*{subsec:app_Formula Derivation Proof}] \nameref{subsec:app_Formula Derivation Proof}
     \item[\ref*{subsec:app_Non-separating Functions}] \nameref{subsec:app_Non-separating Functions}
     \item[\ref*{subsec:app_Discrete Functions}]
     \nameref{subsec:app_Discrete Functions}
     \item[\ref*{subsec:app_Functions on the Circle}] \nameref{subsec:app_Functions on the Circle}
     \item[\ref*{subsec:app_Non-Tame Functions}] \nameref{subsec:app_Non-Tame Functions}
     \item[\ref*{subsec:app_Open and Unbounded Domains}]
     \nameref{subsec:app_Open and Unbounded Domains}
 \end{description}
 \item[\ref*{sec:app_Persistence Diagram Computed Examples}] \nameref{sec:app_Persistence Diagram Computed Examples} (\pageref*{sec:app_Persistence Diagram Computed Examples})
\item[\ref*{sec:app_Diagonal Affine Transformation Theorem Cases}] \nameref{sec:app_Diagonal Affine Transformation Theorem Cases} (\pageref*{sec:app_Diagonal Affine Transformation Theorem Cases})
    \item[\ref*{sec:app_Direct Computation}] \nameref{sec:app_Direct Computation} (\pageref*{sec:app_Direct Computation})
\begin{description}
    \item[\ref*{subsec:app_PD Preimage Illustrations}] \nameref{subsec:app_PD Preimage Illustrations}
    \item[\ref*{subsec:app_Computed PD Flowcharts}]\nameref{subsec:app_Computed PD Flowcharts}
\end{description}
    
\end{description}
\end{multicols}

\section{Introduction}
\label{sec:Introduction}
Repeated measurement of a quantity over time yields a collection of data which is variously called a signal or time series. Applications of these time series include biology, physiology \cite{Biswas2023}; gene expression \cite{Dequeant2008}; microphone readings, EEGs, neural activation \cite{Belton2026}; and financial markets \cite{Jesus2025}, among many others. In many cases, the data has a ``shape'' which conveys some information on the underlying mechanisms that generated it. Sublevel persistent homology allows us to easily identify and characterize these underlying shapes, and distinguish them from the noise. Depending on the domain, this characterization can be rather complicated; for example, point clouds of arbitrary dimension have many different types of filtrations, such as the \v{C}ech complex or the Vietoris-Rips complex \cite[p.27-28]{dey2022computational}. However, the case of one-dimensional signals is much simpler. 

When working with 1D simplicial homology, we use the elder rule to pair up local minima and local maxima, and consequently most of the ``significant" variation (that is, not noise) is concentrated in relatively few pairs, while there are other (typically more numerous) pairs which have a smaller absolute difference in critical value, which correspond to finer details of a signal, which often have less significance or are more attributable to noise. These pairings are graphed in two dimensions (birth time, death time) in a figure called a persistence diagram, which is the primary object of inquiry of this paper. Each pairing has a persistence equal to the amount of time it exists for -- the absolute difference of its coordinates-- and the persistence diagram has a total persistence equal to the sum of the persistence of all of its points.

A previous result from Biswas et al.\ equates the total persistence with the total variation of the function: that is, the total amount of change the signal undergoes in the $y$-direction \cite{Biswas2023}. They demonstrate that this equality holds for continuous functions defined on compact intervals, which are \textit{tame} and \textit{separating}, meaning that all critical points have unique critical values\footnote{Biswas et al.\ extend this result to compact geometric networks, but this is of less interest in this context of 1D signals}. Biswas et al.\ call such functions \textit{generic}.\footnote{Biswas et al.\ later define a function $f$ on the circle to be generic if ``there is no non-empty open interval along which $f$ is constant, so all critical points are isolated, and the values of these critical points are distinct". The implication in this definition is not completely sound, because there exist some functions, such as the continuous extension of $f(\theta) = \theta^2 \sin \frac 1\theta, \theta \in (0,2\pi]$ to $ f(0) = 0$, which are nonconstant yet for which not all critical points are isolated (namely, $\theta = 0$).} While they reason that Cohen-Steiner's stability condition \cite{CohenSteiner2005} implies that any non-separating function can be arbitrarily approximated by a separating function, we use a framework inspired by Belton and Essl \cite{Belton2026} to directly evaluate persistence at non-separating functions, without the need for such approximation.

Moreover, we directly compute the total variation and total persistence as a weighted alternating sum of the sequence of critical values of the function, and we extend this formula to apply to functions with infinite critical points, and even to certain functions on noncompact 1D manifolds. We also extend the equivalence and apply the formula to discrete time series, arguing that the same tools can be applied to examine persistence of discrete or continuous signals.

Building upon this, we establish that the persistence diagram of a function is fully determined by its sequence of critical values, and that any two functions with the same sequence of ordinal critical values in the function will have the same ordinal critical values paired in the persistence diagram. We call this result the Permutation Theorem. We then describe the persistence diagram of a function in terms of this sequence, which relies upon a consistent and calculable but unpredictable permutation function, which has been computed up to $n=6$. Additionally, we characterize two operations under which, when applied to a permutation of a given length, the persistence diagram is invariant, and conjecture that any two permutations with the same persistence diagram can be connected by a sequence of these operations. 

Our final section describes the changes to the persistence diagram of a function when the function undergoes various transformations. First, we fully characterize the changes to a function's persistence diagram when it undergoes a \textit{diagonal affine transformation}-- that is, any combination of vertical or horizontal translations, scalings, and reflections. Looking at other transformations, we try to characterize the set of functions with a particular persistence diagram. With this in mind, we use the stability condition of Cohen-Steiner et al.\ to examine the conditions under which the addition of two functions can be seen in the persistence diagram of the first function. Finally, we show that the success with characterizing the effects of diagonal affine transformations cannot be extended to all linear transformations, as there is nothing in general that can be said comparing the persistence diagrams of a function and its derivative or antiderivative, or their discrete analogs.

\section{Contributions}
\label{sec:Contributions}
\begin{description}
    \item[\ref{sec:Total Variation and Total Persistence}]We generalize the equation of total variation (tVar) and total persistence (tPer) of a function on the interval or circle, from the set of generic continuous functions with finitely many critical points \cite{Biswas2023}, to the set functions which are continuous or discrete on a bounded domain, or else piecewise $C^1$ and absolutely convergent continuous functions on an unbounded domain with two asymptotes. In each of these cases, we provide an explicit formula for the total variation and persistence as a function of its critical values. 
    \item[\ref{subsec:Permutation Theorem}] We prove that the persistence diagram is determined completely by the relative ordering of the sequence of critical values.
    \begin{description}
    \item[\ref{subsec:Persistence Diagram Computation}] We further characterize the persistence diagram of a signal as determined by the relative values of the sequence of local extremum values. We write each of these ordinal sequences as a permutation, and compute the formulas for the persistence diagrams for small $n \leq 6$, though we do not have a general formula for larger $n$. 
    \item[\ref{subsec:Persistence_Diagram_Preimage}] We characterize two operations on the ordering of critical points which preserve persistence diagrams (reversal and \textit{persistence transposition}), and conjecture that all functions with the same persistence diagram can be connected by a sequence of such transformations.
    \end{description}
    \item[\ref{subsec:Diagonal Affine Transformation Theorem}] We completely characterize the effects on the persistence diagrams of diagonal affine transformations upon a signal or time series.
    \begin{description}
    \item[\ref{subsec:Function Addition}] We show that there is no general rule for the effect of adding two functions under their persistence diagrams, but found specific conditions under which this effect can be characterized.
    \item[\ref{subsec:Integration and Differentiation}] We prove that there is no meaningful general relationship between the persistence diagram of a [discrete, continuous] function, and the persistence diagram of its [difference, derivative] or [summation, integral]. We show this by demonstrating that any two functions can be connected by a series of [differences, derivatives], [summations, integrals], and function changes which do not affect the persistence diagram. 
    \end{description}

    \end{description}

\section{Domains, Level Sets, and Persistence}
\label{sec:Domains, Level Sets, and Persistence}
\subsection{Discrete and Continuous Maps}
\label{subsec:Discrete and Continuous Maps}

The terms \textbf{signal} and \textbf{time series} are often used interchangeably to refer to one dimensional functions, which typically track a single response variable as it changes over time. We refer to functions whose domain is disconnected as \textbf{discrete}; for simplicity, the domain is typically assumed to be the set of integers $\{1,\ldots, n\}$, which we will denote $I_n$. Conversely, we refer to \textbf{continuous functions} as those which meet the traditional real analysis definition of continuity: At each point $x_0$ , for all $\epsilon > 0$, there exists a $\delta > 0$ such that, given any $x$ in the neighborhood of $x_0$, then $|x_0-x| < \delta \implies |f(x_0)-f(x)|<\epsilon$. For a continuous function $f$, points $x_0$ for which there exist points $x_0+\epsilon, x_0-\epsilon$ in the domain of $f$, for sufficiently small $\epsilon>0$, are called \textbf{interior} points; points for which only one of $f(x-\epsilon)$ or $f(x+\epsilon)$ are defined are called \textbf{boundary} points. On the other hand, if $f$ is a discrete function, we call a point $x_i$ for which $f(x_i)$ is defined a \textit{boundary point} if either $f(x_i-k)$ or $f(x_i+k)$ fails to exist for \textit{all} $k >0$. Conversely, $x_i$ is an interior point of $f$ if there exist $k_1,k_2>0$ such that $f(x_i), f(x_i-k_1)$, and $f(x_i+k_2)$ exist. We note that if the domain of a discrete function $f$ is $I_n$, then the point $(i,f(i))$ is a boundary point if and only if $i=1$ or $n$.  

\subsection{Domains, Genericity, and Critical Points}
\label{subsec:Domains, Genericity, and Critical Points}
For the purposes of this paper, the domain of a continuous function will be either $\bb I = [0,1]$, the \textbf{unit interval}, or $\bb S^1$, the \textbf{unit circle}. While the former domain has boundary points at $0$ and $1$ (or more generally, $x_1$ and $x_n$, the smallest and largest $x$ values), a function defined on the circle has no boundary points, and so every point is an interior point. We will only consider functions which are \textit{bounded}; that is, circles and finite intervals. 

Just as in calculus, we define the set of \textbf{critical points} of a continuous function $f$ to be the set of points for which the derivative $f'$ of $f$ is either zero or undefined. Thus, all boundary points are critical points, and most signals will also have interior critical points. However, when examining either persistence or variation, we only care about those critical points at which the derivative $f'$ actually changes sign, and not a \textit{degenerate} or \textit{saddle} point, at which the derivative reaches a local maximum or minimum at zero, and thus reaches it without changing sign. Therefore, unless otherwise specified, we will use the term \textit{critical point} to refer only to \textbf{local extrema}, at which the sign of $f'$ actually changes [or alternatively, only exists on one side, as all boundary points are local extrema]; that is, we will disregard inflection points. The local extrema will be classified as \textbf{local maxima} or \textbf{local minima} according to the direction of the sign change. The \textbf{global} maximum or minimum of a function is the maximal maximum or the minimal minimum, respectively. 

For a discrete function, we define the critical points in terms of differences, not derivatives: That is, a point $x_i$ is a critical point of a discrete function $f$ if
\[\text{sign}(f(x_i) - f(x_{i-1})) \neq \text{sign}(f(x_{i+1}) - f(x_i)),\]
or if either $f(x_{i-1})$ or $f(x_{i+1})$ does not exist. 

For a given critical \textit{point} $c_i = (x_i,y_i)$, we refer to the $y$-value $y_i$ as the \textbf{critical value} associated with the critical point. Because each $x_i$ has exactly one point associated with it, we will simply refer to $x_i$ as the \textit{point}, and $y_i$ as the \textit{value}, to distinguish them. \textit{Extremum}, \textit{maximum}, and \textit{minimum} technically refer to the critical value, but can sometimes refer to the point itself when it is clear from the context, such as for \textit{generic} functions. 

Following the example of \citet{Biswas2023}, we will define a \textbf{generic} function as a function with a bounded domain, and whose critical points are all isolated and distinct in value. Thus, as we are assuming a bounded domain, we will call a function $f$ \textbf{non-generic} if it fails to be either tame or separating: that is, if $f$ has an infinite number of critical points (which could be all isolated, or include a horizontal segment), or has multiple distinct critical points with the same critical value. \\

We will use the convention that a critical point has coordinates $(x_i,y_i)$: that is, $y_i = f(x_i)$.  We will sometimes still write out $y_i$ as $f(x_i)$ for clarity, especially when multiple functions are involved. When we want the $y$ levels to be sorted, we will denote $y_{(i)}$ as the $i$th critical value, such that $y_{(i)}<y_{(j)} \iff i < j$. 

\subsection{Level sets}
\label{subsec:Level sets}
\begin{figure}[h]
    \centering
    \includegraphics[width=0.25\linewidth]{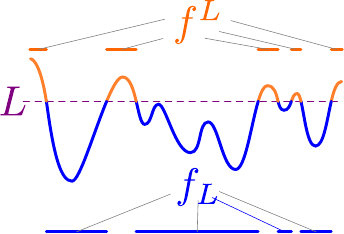}
    \caption{The sublevel set $f_L$ and superlevel set $f^L$ illustrated for some value $L \in M$.} %
    \label{fig:sub_and_superlevel_sets}
\end{figure}

Let $f : M \to \rr$ be a smooth ($C^\infty$), tame Morse function; that is, $f$ is an infinitely differentiable function with finitely many critical points in a bounded interval, which each have unique critical values. Given this function $f$, and an interval $I \subseteq \rr$, the \textbf{interval level set} of $f$ with respect to $I$ is $\mu_I = f^{-1} (I) =\{x \in M \ | \ f(x) \in I \}$. In particular, if $I =(-\infty, \alpha]$, this is a \textbf{sublevel set}; if $I = [\alpha, \infty)$, it is a \textbf{superlevel set}. If $I = [\alpha, \alpha ] = \{\alpha\}$ is a singleton, then we call $\mu_I$ the \textbf{level set} of $f$ at $\alpha$. The sublevel and superlevel sets of a function $f$ are shown in Figure \ref{fig:sub_and_superlevel_sets}. 

As our function is defined on a finite interval $M = [m^-,m^+]$, the extreme value theorem tells us that range of $f$ is bounded on a closed interval, which we will call $[y^-, y^+]$. We shall assume $M$ is a finite interval $[m^-, m^+]$ from this point forward, unless indicated otherwise, and therefore that there exist $y^- := \min _{x \in M} f(x)$ and $y^+ := \max _{x \in M} f(x)$.\footnote{If $f$ has $n$ critical values, then we can write $y^+ = y_{(n)}$ and $y^- = y_{(1)}$; we use the former notation to emphasize that the extrema are global.} With this in mind, we will notate the sublevel sets $f_\alpha := f^{-1} [y^-, \alpha]$ and superlevel sets $f^\alpha := f^{-1} [\alpha, y^+]$, for all values $\alpha \in [y^-, y^+]$.\\

To characterize the sublevel persistence of a function $f$, we must characterize the sublevel sets of every level $\alpha$. A standard result \cite[Theorem~1.3]{dey2022computational} in topological data analysis says that whenever an interval $[a,b]$ contains no critical points of a \textit{tame Morse} function $f$, then its sublevel set $f_b$ is \textbf{diffeomorphic} to $f_a$; that is, smoothly homeomorphic in both directions. %

As our function $f$ meets these criteria, we can partition the set of all non-critical points in the range of $f$ into finitely many open intervals within which the sublevel sets of all closed intervals are diffeomorphic, and within which the superlevel sets of all closed intervals are diffeomorphic. \\

Therefore, to characterize the sublevel sets of a function $f$ up to homeomorphism (or in fact diffeomorphism), one need only characterize these intervals, and of course the critical points which separate them. 

Furthermore, at a non-boundary critical point, the sublevel set of the critical value $y_0$ of $f$ is  diffeomorphic to the sublevel set of $y_0+\epsilon$ for sufficiently small $\epsilon$; the superlevel set of $y_0$ is likewise diffeomorphic to the superlevel set of $y_0-\epsilon$. This is a consequence of the fact that the sublevel and superlevel sets take $f^{-1}$ over an interval which includes $\alpha$ (that is, closed at least on the side of $\alpha$). 

Therefore, to characterize the level sets of points in between critical values $y_{(i)} <y<y_{(i+1)}$, it suffices to compute the sublevel set of $f$ at $y_{(i)}$ and the superlevel set of $f$ at $y_{(i+1)}$. Repeating this process for all critical points, we complete the task of characterizing the sublevel and superlevel sets of the range of $f$ by computing merely the sublevel and superlevel sets of its critical points. \\

We have reduced our task, the characterization of the persistence of a function, to a finite number of calculations. However, we still wish to find an effective way to visualize the differences between the directional level sets at (and therefore between) various critical points, which is not always easy to do%
. This is the purpose of the \textbf{persistence diagram}. A persistence diagram will show \textbf{topological features}, identifying longer-lasting signals which often represent physically meaningful phenomena, as opposed to shorter-lived features, which are more likely due to noise. In particular, the sublevel persistence will track zero-dimensional features, which are connected components of the sublevel sets. The superlevel persistence, conversely, will track one-dimensional features (that is, loops), when the points in the superlevel set are identified with one another, shown in Figure \ref{fig:DrK_Rel_Hom}.  \\ %

\subsection{Persistence Diagrams}
\label{subsec:Persistence Diagrams}

The first type of persistence diagram, the \textbf{ordinary persistence diagram}, displays the sublevel persistence of the function $f$ at various points $\alpha$. In particular, we track the 0-dimensional homology-- that is, the number of connected components of the curve's level sets as we increase from $f_{y^-}$ to $f_{y^+}$. As we are looking at the evolution of sublevel sets of $\alpha$ from $y^-$ to $y^+$, we will imagine a slider starting at $f^{-1}(y^-)$; that is, $f^{-1}[y^-, \alpha] = f_\alpha$ where $\alpha = y^-$. As $\alpha$ increases, the interval $[y^-,\alpha]$ grows continuously, and so too does $f_\alpha$, until $f_{y^+} = [m^-, m^+]$. This continuous growth allows us to imagine $\alpha$ as a sort of ``time'', and, when working with sublevel persistence, we describe points in $M$ contained in $f_{\alpha_1}$ as being reached ``sooner'' and points in $f_{\alpha_2}$ but not $f_{\alpha_1}$ being reached ``later'' if $\alpha_1 < \alpha_2$. Because superlevel persistence begins at $y^+$ and moves towards $y^-$, when speaking of superlevel persistence the above terminology applies when $\alpha_1 > \alpha_2$. \\

With this framework, at each non-boundary critical point of $f$, we track the \textit{birth} or \textit{death} of a topological feature: A feature is ``born'' at a critical point if it does not exist before the critical point is reached, and exists after; a feature similarly ``dies'' if it exists (distinctly) before reaching a critical point, but merges with another feature upon reaching the critical point. By convention, when two features meet, the ``older'' feature (that is, the feature which was born before the other) continues on, while the ``younger'' dies; this convention is called the \textit{elder rule}. We will use the notation $k_C$ to denote the connected sublevel component which was born at $f_C$, and $k_C^D$ if $k_C$ dies at a critical point with value $D$. A capital $K$ likewise signifies loops in the superlevel persistence. 

The birth and death times of each topological feature are recorded and plotted as an ordered pair. The birth times are plotted on the $x$ axis, and the death times on the $y$ axis. Accordingly, in the ordinary persistence diagram, which tracks sublevel persistence, all the points are plotted above the line $x=y$, as the birth time of every point precedes its death time. \\

\begin{figure}
    \centering
    \includegraphics[width=\textwidth]{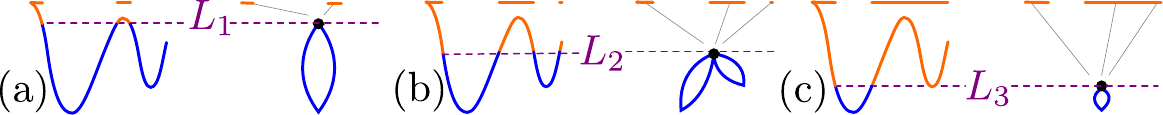}
    \caption{In the superlevel persistence, we identify (``pinch'' together) all parts (in orange) of the domain mapped above the dotted $L$ line, and we count the number of loops this creates in the portion of the domain mapped below the level $L$ (in blue). Any part of the function below the line $L$ which does not create a loop (a ``tail'') is disregarded.}
    \label{fig:DrK_Rel_Hom}
\end{figure}

The \textbf{relative persistence diagram} tracks the relative persistence of superlevel sets of $f$ and therefore considers $f^{-1} [\alpha, y^+] $ as $\alpha$ varies from $y^+$ to $y^-$. Instead of tracking connected components, the superlevel persistence tracks the number of loops created when all of the points in the superlevel set are identified with one another, or ``pinched together'', as seen in Figure \ref{fig:DrK_Rel_Hom}. This is tracked similarly to the sublevel persistence: Whenever the superlevel set reaches a new local maximum, a new loop is ``born'', with one exception we will explain momentarily. When an interior local minimum is reached, two loops join, and the younger one ``dies''; like in sublevel persistence, no deaths occur at boundary points. The one exception which distinguishes the computation of superlevel persistence from that of sublevel persistence is that, when $f$ is defined on an interval, the first maximum reached does not birth a connected loop, and therefore the global maximum is not represented in superlevel persistence. A consequence of this is that the  number of loops which superlevel persistence will count at a given level $L \in M$ will always be one less than the number of connected components present in the superlevel set whenever $M$ is an interval; conversely, when $M = \bb S^1$, the unit circle, then the number of loops counted will be equal to the number of connected components. 

Analogous to the sublevel persistence, as we are moving from $y^+$ to $y^-$, a ``later'' time of death is associated with a lesser value of $\alpha$. Thus, when we plot the birth times against the death times of superlevel features, the death time is less than the birth time, and so every point is plotted \textit{below} the dashed diagonal line, on which the birth time equals the death time. 

\begin{figure}
    \centering
    \includegraphics[width=.8\linewidth]{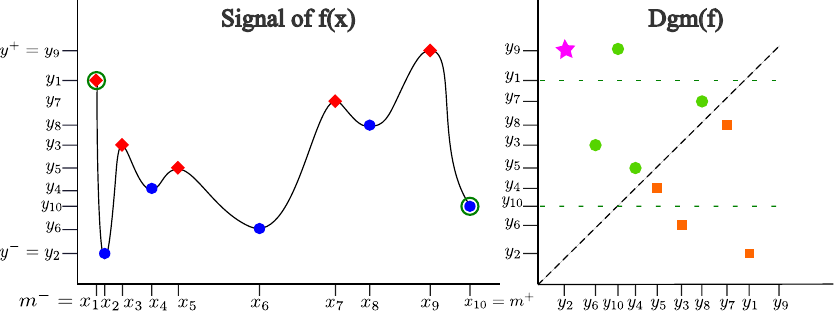}
    \caption{The left figure maps the signal $y$ over time $x$. The right side is the persistence diagram $\dgm (f)$.
    On the left, blue dots = local minima, red diamonds = local maxima, and green circles = boundary points. On the right, pink star = $\ess(f)$, green dots = $\ord(f)$, and orange squares = $\rel (f)$.} \label{fig:first_sample}
\end{figure}

Finally, the \textbf{essential persistence diagram} displays those features which are born and do not die. On intervals, this consists of a single point at $(y^-, y^+)$; on circles, as there is an inherent loop in the structure, the essential persistence diagram contains an additional point at $(y^+, y^-)$. \\

\begin{figure}[h]
    \centering
    \includegraphics[width=.75 \linewidth]{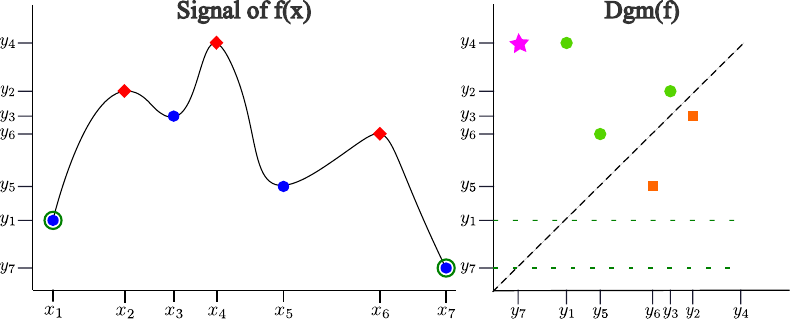}
    \caption{Sample signal (L) and persistence diagram (R). The odd number of critical points creates an asymmetry in the number of points in $\ord (f)$ and $\rel (f)$. The green, horizontal sparsely-dashed lines at $y_1$ and $y_7$ represent the invisibility of $x_1$ and $x_7$ to the superlevel relative persistence, as they are endpoints.}
    \label{fig:second_sample}
\end{figure}

When these three persistence diagrams are combined, we have our \textbf{extended persistence diagram}, which we will denote as $\dgm(f):=\ord(f) \sqcup \rel(f) \sqcup\ess(f)$, with the ordinary, relative, and essential diagrams labeled $\ord(f), \rel(f),$ and $\ess(f)$, respectively. Figures \ref{fig:first_sample} and \ref{fig:second_sample} illustrate a continuous function $f$ on the left, and its extended persistence diagram on the right. On the left side, the graph of $f$ labels local minima with blue dots, and local maxima with red squares, and $m^-$ and $m^+$ are circles in green. On the right, $\dgm (f)$ shows $\ord (f)$ with orange points, $\rel (f)$ with lime green squares, and $\ess (f)$ with pink stars. Finally, the horizontal green lines on the right show the values for which there is no point in $\ord (f)$ or $\rel (f)$, which are precisely the critical values which occur at the boundary points on the left. 

\section{Total Variation and Total Persistence}
\label{sec:Total Variation and Total Persistence}

The \textbf{variation} of a real-valued map is a way to quantify the amount a function fluctuates in a given section. In particular, when $f$ is a one dimensional, piecewise differentiable continuous map on a compact interval, it is absolutely continuous\footnote{A function is \textbf{absolutely continuous} if, for any $\epsilon > 0$ there exists a $\delta>0$ such that, for any finite sequence of disjoint intervals $\{(x_i,y_i)\}_{i=1}^N $, within $D$ the domain of $f$, whenever the sum of their length
\[\sum_{i=1}^N|y_i-x_i| < \delta,\]
the net displacement over those intervals 
\[\sum_{i=1}^N |f(y_i)-f(x_i)|<\epsilon.\]

}, and its variation is the integral of its absolute derivative: 

\begin{equation*}\text{Var}(f) = \int |f'(x)| \dee x.\end{equation*}

The \textbf{total variation} of a function $f$ is simply the variation of all of $M$ (the domain of $f$)\footnote{Formally, total variation is defined more generally as the supremum among all partitions, of the sum of the variation on each partition of the domain \cite[p.126]{Athreya2006}. However, the integral formulation holds whenever $f$ is absolutely continuous (see Section \ref{subsec:Non-Tame Functions}), which will be assumed unless specified, including for every signal outside Section \ref{subsec:Non-Tame Functions} and the corresponding Appendices \ref{subsec:app_Non-Tame Functions}-\ref{subsec:app_Open and Unbounded Domains}. \label{footnote:tVar_partn_definition}}:

\begin{equation*} \tvar (f) = \int_{m^-}^{m^+} |f'(x)| \dee x.\end{equation*}

The total variation can be thought of as the ``total distance traveled'' by a point moving along the $y$-axis according to the signal. This sum can be computed by partitioning the integral at the critical points, breaking it into monotonic segments which can be evaluated directly, as seen in Figure \ref{fig:total_variation}.\\

The \textbf{total persistence} of a function is the sum of the persistence of all of its features, as represented by points on the (extended) persistence diagram. In particular, the \textbf{persistence} of a point on $\dgm(f)$ is the absolute difference between its death and birth time-- that is, the duration for which it persists. This gives us a formula for the total persistence: 

\begin{equation*}\|\dgm(f)\|_1=\sum _{(A,B)\in \dgm(f)}|B-A|.\end{equation*}

\subsection{Computation and Comparison}
\label{subsec:Computation and Comparison}

A central result of Biswas et al.\ \cite[Corollary 4.6]{Biswas2023} was the equation of total variation and total persistence. Though this identity was previously known for functions whose domains were circles \cite{Dequeant2008}, Biswas et al.\ extend it to geometric networks. Based on their definition, the interval should be seen as a special case of embedded geometric network with two degree-one vertices.

We show that, given the sequence $\{y_i\}_{i=1}^n$ of critical values, the total variation can be computed as
\[\tvar(f) = \left|-y_1 + (-1)^ny_n +2\sum_{i=2}^{n-1} (-1)^iy_i \right|= \tper (f),\]
the weighted alternating sum of the critical values. When the first critical point is a local minimum, the sum will be nonnegative, and the absolute value signs are superfluous. This formula motivates the later Permutation Theorem (Theorem \ref{thm:PermThm}), which will confirm that the total persistence, and in fact the entire persistence diagram from which it is calculated, is completely determined by the sequence of critical values.

\subsection{Formula Derivation}
\label{subsec:Formula Derivation}

Consider a generic function $f$, with critical points $(x_1, y_1), (x_2, y_2), \ldots, (x_n, y_n)$.

Without loss of generality, assume that $y_2>y_1.$ This proof relies upon the computation of total variation and the total persistence, which are identical between $f$ and $-f$ (see Appendix \ref{app:negation}). Alternatively, if $y_1 > y_2$, switch ``odd'' for ``even'' throughout the following proofs.

\begin{figure}
    \centering
    \includegraphics[width=0.6\linewidth]{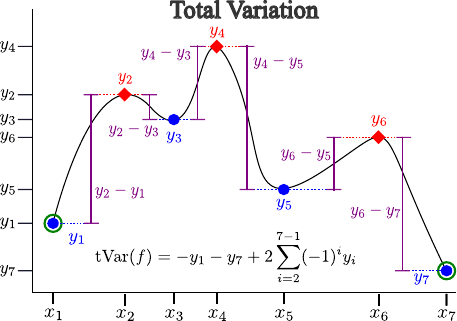}
    \caption{The total variation can be calculated by breaking the interval into monotonic segments at the critical points, and taking the difference of the endpoint values. }
    \label{fig:total_variation}
\end{figure}

First, observe that the function $f(x)$ is monotonic between any two consecutive\footnote{That is, two critical points with no other critical point between them in the $x$ direction.} critical points. Therefore, between any two such points  $c_i = (x_i, y_i)$, 
\begin{equation*}\int _{x_i}^{x_{i+1}}|f'(x)| \dee x = |y_{i+1}-y_i|.\end{equation*}
Thus, we can decompose our interval $M$ into the intervals between its critical points, as seen in Figure \ref{fig:total_variation}. With our assumption from earlier, since the critical points are indexed by their $x$ coordinate, we know that even-indexed points $y_{even}$ are local maxima, and odd-indexed points $y_{odd}$ are local minima. This means that every $|y_{i+1} - y_i|$ has the form $y_{even} - y_{odd}$, and so we have that
\begin{equation*} \tvar(f) = (y_2-y_1) + (y_2-y_3) + (y_4-y_3) + (y_4-y_5) + \ldots, \end{equation*}

which simplifies to \begin{equation*}
    \tvar (f) = -y_1 + (-1)^n y_n + 2\sum _{i=2}^{n-1} (-1)^i y_i.
\end{equation*}

\begin{figure}
    \centering
    \includegraphics[width=0.4\linewidth]{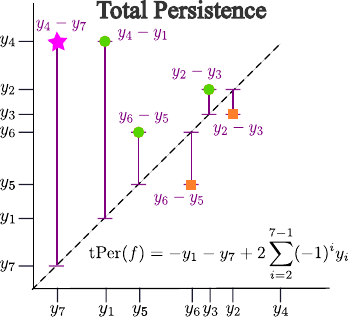}
    \caption{The persistence of each point in $\dgm (f)$ is shown, with $\ord (f)$ as green circles, $\rel (f)$ as orange squares, and $\ess (f)$ as a magenta star.  }
    \label{fig:total_persistence}
\end{figure}

\begin{table}[h!]
    \centering
    \begin{tabular}{c | c c c c c}
        \hline
        & $y_1$ & $y_{even}$ & $y_{odd}$ & $y_n$ & $y^+$\\
        \hline
        Sublevel: & Birth & Death & Birth & ifOdd & +Death\\
        Superlevel: & * & Birth & Death & ifEven & -Birth\\
        \hline Total:  & 1 & 2 & 2 & 1 & *
    \end{tabular}
    \caption{This table shows how many times a point in $\dgm (f)$ appears at each critical value of $f$. The right boundary critical value $y_n$ has one birth, located in the sublevel persistence if $n$ is odd, or in the superlevel persistence if $n$ is even. The parity of the first three columns depends on the $y_1<y_2$ assumption. More detail can be found in Appendix \ref{app:tVar=tPer_Proof_1_Elab}.}
    \label{tab:tPer_contribution_table}
\end{table}

As we see in Figure \ref{fig:total_persistence}, the total persistence of $f$ can be similarly broken down into differences of critical values $y_i$. Because we are still assuming $y_2 > y_1$, odd-indexed critical points are local minima, and even-indexed points are local maxima. For more explanation of this process, see Appendix \ref{app:tVar=tPer_Proof_1_Elab}.

Under this assumption, Table \ref{tab:tPer_contribution_table} summarizes the contributions of each critical point to the persistence diagram $\dgm (f)$. As before, since each local maximum is greater than the local minimum with which it is paired, the formula 

\begin{equation*}
    \tper = \sum_{(A,B)\in \dgm (f)} |A-B|
\end{equation*}
reduces to 
\begin{equation*}
    \tper = \sum_{(A,B)\in \dgm (f)} |A-B| =-y_1 + (-1)^ny_n +\sum_{i=2}^{n-1} (-1)^iy_i.
\end{equation*}

\vspace{5pt}

By harnessing the alternating nature of critical points, we can directly compute the total variation of $f$ to be 
\begin{equation*}
\tvar (f) = \tvar(f) = \int_{m^-}^{m^+} |f'(x)| \dee x =  -y_1 + (-1)^ny_n + 2\sum_{i=2}^{n-1} (-1)^i y_i.
\end{equation*}

By counting the birth times and death times, we get that the total persistence of $f$ is

\begin{equation*} \tper (f) = \| \dgm (f) \| _1 = -y_1 + (-1)^ny_n + 2\sum_{i=2}^{n-1} (-1)^i y_i. \end{equation*}

As these sums are identical, we have that  $\tvar (f) = \tper (f). \ \blacksquare$

\subsection{Circular Domains}
\label{subsec:Circular Domains}
Many signal applications are periodic or modeled with Fourier series, and are therefore often expressed more naturally as functions on the circle than on an (arbitrary) interval. When $f$ is a function on the circle $\bb S^1$, then the equality tVar($f$) = tPer($f$) holds, but the proof is different from the interval case.

Since there are non boundary points, every point is an interior point, and so each local minimum corresponds to a birth in the sublevel half of $\dgm (f)$, and a death in the superlevel half, and vice versa for each local maximum. Additionally, when we are on the circle, all non-global extrema will meet the essential component before the essential component ends. Therefore, the only points at $y^-$ and $y^+$ on $\dgm (f)$ are in the essential diagram; there are no ``additional'' points.

Furthermore, in superlevel persistence, there is a loop created when the global maximum is reached, because the line itself is connected. This loop persists for the entire duration of the superlevel persistence, and so there is a point in $\ess (f)$ created during the superlevel persistence, namely at $(y^+, y^- )$. 

Because every critical value corresponds to both a birth and a death, the persistence diagram is symmetric, as the only sources of asymmetry were the births associated with boundary points, and the consequent extra death at the global maximum. This symmetry is also a corollary of Biswas's Theorem 3.2 \cite{Biswas2023}. \\

We will model this function by considering a function on the interval, and identifying the boundary points $c_1$ and $c_n$ with one another.

Therefore, for a function $f : \bb S^1 \to \rr$, the total persistence \footnote{The absolute values account for the cases when the $y_1<y_2$ assumption does not hold; that is, when the left endpoint is a local maximum.}
\[ \| \dgm (f) \|_1 = \sum_{(A,B)\in \dgm (f)} |A-B| = 2\left|\sum_{i=1}^{n-1} (-1)^iy_i\right|.\]

Similarly, we can compute the total variation by simply setting $y_1 = y_n$ to yield 
\[\tvar (f) = \left |-y_1 + (-1)^n y_n + 2\sum _{i=2}^{n-1} (-1)^i y_i \right |. \]
If $n$ is even, then $y_1=y_n$ are not an actual critical point, and so their sum contributions cancel out.%
Conversely, if $n$ is odd, then $y_1 = y_n$ is a true critical point, and so is recovered in the sum. Combining these results, 
\[\tvar (f) =\begin{cases} 2\left |\sum _{i=2}^{n-1} (-1)^i y_i \right |, & \text{ if }n\text{ is even}\\ \\
    2\left |\sum _{i=1}^{n-1} (-1)^i y_i \right |, & \text{ if }n\text{ is odd}
\end{cases}. \]

\subsection{Non-Separating Functions}
\label{subsec:Non-Separating Functions}

As Biswas showed, $\tvar (f) = \tper (f)$ for any continuous, \textit{generic} function. Inspired by Belton et al.\ \cite{Belton2026}, we will extend this finding to nongeneric functions, and then to discrete domains. Specifically, we will first relax the separability requirement, and we will address non-tame functions in section \ref{subsec:Non-Tame Functions}. 

Recall that a continuous function $f: \bb I \to \rr$ is generic if there are finitely many critical values, and if all critical points have distinct values; that is, $f$ is \textit{tame} and \textit{separating}. One appeal of the separating restriction is that it guarantees that the elder rule is always unambiguous, because whenever two components meet, the distinct critical levels at which they were born cannot be the same. We will next examine functions which are tame but \textit{not} separating.

When we get rid of this restriction, if two components $k_{a}$ and $k_b$ meet, where $f(a) = f(b) = c$, then there is no way to adjudicate which point's component survives. It is equally possible that $k_a$ dies and $k_b$ survives, or that $k_a$ survives and $k_b$ dies. 

\begin{figure}
    \centering
    \includegraphics[width=0.9\linewidth]{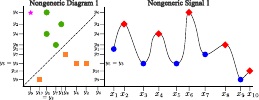}
    \caption{A simple nonseparating signal, where $y_3=y_5$. In the persistence diagram, the shared level has two births and two deaths, a consequence of the shared points both being interior points and nonglobal extrema. }
    \label{fig:non-separating1}
\end{figure}

\begin{figure}
    \centering
    \includegraphics[width=0.9\linewidth]{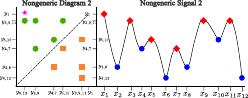}
    \caption{An example of a highly non-separating function. The critical points $x_2$ and $x_8$ share a critical value, as do the pairs $x_4$, $x_7$; $x_5$,$x_{10}$; and $x_6, x_{12}$. The three critical points at $x_3$, $x_9$, and $x_{11}$ all share a critical value as well. }
    \label{fig:non-separating2}
\end{figure}

However, in a persistence diagram, components are not graphed by their $x$-values, but by their $y$-values. Since $f(a)=f(b)=c$, the persistence diagram will graph two points in the column $c$: one at the level at which $k_a$ and $k_b$ meet, and one at which that merged component dies (or at $y^+$ if it is the essential component). Thus, the assignment of seniority to $k_a$ or to $k_b$ is arbitrary, and the persistence diagram produced either way is identical. A simple example of this can be seen in Figure \ref{fig:non-separating1}, and we can see in Figure \ref{fig:non-separating2} that even a multiplicity of overlaps, which creates many potential combinations of birth and death component choices, will only create one persistence diagram. 

\begin{figure}
    \centering
    \includegraphics[width=0.9\linewidth]{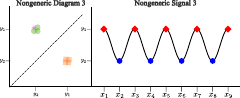}
    \caption{Another very non-separating function, specifically a (shifted) cosine wave. The repeated maxima and minima create multiple points in $\dgm(f)$, at $(y_1,y_2)$ and $(y_2, y_1)$, which are slightly dispersed for visibility.}
    \label{fig:non-separating3}
\end{figure}

One aspect of persistence diagrams of nongeneric functions is that they can feature multiple points at the same coordinates (or alternatively, points with multiplicity greater than one), as seen in Figure \ref{fig:non-separating3}. However, the equality of total persistence and total variation is still maintained, as one can easily see that the total variation and total persistence in Figure \ref{fig:non-separating3} is $8(y_1-y_2)$.

Another aspect present in non-generic functions is \textit{horizontal segments},  which Belton calls ``flats'' \cite{Belton2026}. Given $x_1<x_2 \in D$, the interval $[x_1,x_2]$ is a \textbf{horizontal segment} of $f$ if, for all $x \in [x_1,x_2]$, $f(x)$ is constant. Because every point in a flat horizontal segment is a critical point, and there are uncountably infinitely many points in any interval of positive width, any function with any flat segment is neither tame nor separating, and therefore is not generic. 

\begin{figure}
    \centering
    \includegraphics[width=0.9\linewidth]{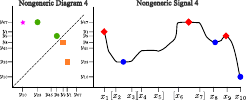}
    \caption{A non-separating signal with a horizontal component. The horizontal segment at $y_{45}$ is a \textit{monotone}, and so it has no extrema, and thus no contribution to $\dgm (f)$. The horizontals at $y_{23}$ and $y_{67}$, conversely, are a minimum and maximum, respectively, and function accordingly. }
    \label{fig:non-separating4}
\end{figure}

Horizontal segments are treated much the same as multiple isolated critical points: for each segment, one point in the segment is arbitrarily chosen to name the component $k$ which continues on from the shared critical value until its death, and a point is plotted at those coordinates irrespective of the choice of component, as seen in Figure \ref{fig:non-separating4}.

However, just like regular critical points, not all horizontal segments are plotted in the persistence diagram, but only those which contribute any persistence; that is, local maxima and minima. %
If the function is only increasing or only decreasing both before and after the horizontal interval, then we say the segment is \textit{monotone}, and it does not appear in $\dgm(f)$. Indeed, in Figure \ref{fig:non-separating4}, the horizontal segment at $y_{45}$ is monotone, as the function is increasing both before and after it, and so no points occur in the row or column of $y_{45}$ in $\dgm (f)$. More formally, a segment $[a,b]$ is called \textbf{monotone} on a function $f$ if $f$ is constant on the interval, and, for arbitrarily small $\epsilon >0$, 
\[\sign(f'(a-\epsilon)) = \sign(f'(b+\epsilon)) \neq 0.\]
A corollary of this definition is that a  monotone horizontal segment cannot contain a boundary point. Note that this definition is specific to the continuous case; with discrete functions, non-constant segments can also be monotone. \\

We can now prove that the total persistence and total variation of continuous functions remain equal even among non-separating functions, with the only limitation being a finite number of critical points and horizontal segments (as this is sufficient to guarantee finite variation); that is, $f$ must have a reduction $
\tilde f$ which is tame.

Let $f : \bb I \to \rr$ be a continuous function with finitely many critical points and horizontal segments, which is not necessarily generic. First, if our function has any horizontal segments, we can create another function $\tilde f : \bb I / \sim \to \rr$ by ``shrinking'' all horizontal segments to points; for more details, see Corollary \ref{cor:equiv_rel} of Lemma \ref{lem:hor_seg}. We can show that this function $\tilde f$ has the same total persistence and total variation, so we can show $\tvar (f) = \tper (f)$ by simply showing that $\tvar (\tilde f) = \tper (\tilde f)$. 

The motivation for using $\tilde f$ is that all critical points are isolated, and so we can apply the same counting methods to $\tilde f$ as we could to our generic functions in earlier sections. As a consequence, we can assume from here on out that $f$ has only isolated critical points, with the understanding that we will substitute $\tilde f$ for $f$ should this prove not to be the case. \\

Because there are finitely many critical points, there are in particular finitely many local extrema, which we will denote $c_1, \ldots, c_n$. That is, we are disregarding any \textit{degenerate} critical points, such as inflection points. We still know that $f$ is monotonic between all local extrema, so we can compute $\tvar(f)$ as before: 
\[ \tvar (f) = \int_{m^-}^{m^+} |f'(x)| \dee x = \sum _{i=1} ^{n-1}\int _{x_i}^{x_{i+1}}|f'(x)| \dee x =  \sum _{i=1} ^{n-1}|y_{i+1}-y_i| = -y_1+ (-1)^n y_n +2 \sum _{i=1} ^{n-1} (-1)^iy_i.\]

The counting argument summarized in Table \ref{tab:tPer_contribution_table} applies identically to our non-separating function, and for the same reason, our total persistence reduces to 

\begin{equation*}
    \tper = \sum_{(A,B)\in \dgm (f)} |A-B| =-y_1 + (-1)^ny_n +\sum_{i=2}^{n-1} (-1)^iy_i,
\end{equation*}
which is equal to $\tvar(f)$. $\blacksquare$ \\

\subsection{Discrete Functions}
\label{subsec:Discrete Functions}
Now, let $g : X \to \rr$ be a \textit{discrete} function, from a finite set $X \subsetneq \nn$ to the real numbers. We will show that $\tvar(g) = \tper(g)$, by way of showing that there exists a continuous function $f$ with identical total persistence and total variation to $g$, such that $\tvar (f) = \tper(f)$. 

In particular, let $f$ be a \textit{linear interpolation}\footnote{In fact, any interpolation which is monotonic between every pair of adjacent critical points will suffice; a linear interpolation was just chosen for ease of computation. } of $g$: that is, $\forall x \in [1,n]$, 
\[ f(x) = \begin{cases}
        g(x) & \text{ if } x \in \nn \\
        (\lceil x \rceil - x)g(\lfloor x \rfloor) + (x-\lfloor x \rfloor)g(\lceil x \rceil) & \text{ if } x \notin \nn
    \end{cases}.\]

    \begin{figure}
        \centering
        \includegraphics[width=0.9\linewidth]{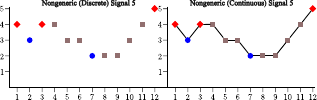}
        \caption{A discrete function $g$ (left), and its linear interpolation $f$ (right). The grey squares denote a non-extremum, or additional points in a horizontal flat extremum. For example, the local maximum at $x=3$ could be moved to 4, and the local maximum at $x=7$ could be moved to 8 or 9.}
        \label{fig:non-separating5}
    \end{figure}

First, observe that if $g(j)$ is not a local extremum, then $|g(j+1) - g(j-1)| = |g(j+1)-g(j)| - |g(j) - g(j-1)|$. Therefore, in the sum 
\[\tvar(g) = \sum _{i=1}^{n-1} |g(i+1)-g(i)| = |g(2)-g(1)| + \ldots + |g(n)-g(n-1)|,\]
we can replace each such $|g(j+1)-g(j)| - |g(j) - g(j-1)|$ with $|g(j+1) - g(j-1)| $, resulting in 
\[\tvar(g) = |g(2)-g(1)| + \ldots +|g(j-1)-g(j-2)| + |g(j+1) - g(j-1)| + |g(j+2) - g(j+1)| + \ldots + |g(n-1)-g(n)|.\] \\

\begin{figure}
    \centering
    \includegraphics[width=0.9\linewidth]{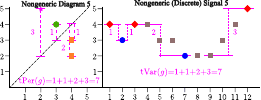}
    \caption{The total persistence (L) and variation (R) of functions $g$ (and equivalently, its linear interpolation $f$), from Figure \ref{fig:non-separating5}. The points in the persistence diagram (L) do not depend on the choice of maxima and minima in Figure \ref{fig:non-separating5}.}
    \label{fig:non-separating5.1}
\end{figure}

The result of this can be seen in the total variation shown in the right side of Figure \ref{fig:non-separating5.1}, with regard to the discrete function $g$ and its linear interpolation $f$ in Figure \ref{fig:non-separating5}. \\

Similarly, we observe that only the local extrema contribute any persistence to $\tper(g)$. Because the pairing of local extrema is identical between $f$ and $g$, their persistence diagrams will be identical, and $\tper(f) = \tper(g)$ follows. 

Because we know $f$ is a continuous function with finitely many isolated critical points and horizontal segments, $\tper(f) = \tvar(f)$. As we have equated each of these totals with their counterpart in $g$, we have
$\tper(g) = \tper(f) = \tvar(f)=\tvar(g)$, so $\tper(g) = \tvar(g)$, and we are done. 

\subsection{Non-Tame Functions}
\label{subsec:Non-Tame Functions}
All of the generalization so far have assumed that $f$ has a finite number of local extrema; that is, that $f$ is tame. However, this is not always the case; when $f$ is continuous, even some functions on a compact domain $D$ have infinite oscillations in $f(x)$ over a (bounded) sub-interval within $D$. For example, the function \[ f(x) = \begin{cases}
x^2\cos\left(\frac 1 x\right), & x \in (0,2] \\
0, & x=0
\end{cases}\]
is continuous, and has (nondegenerate) critical points whenever its derivative, 
\[\frac {\dee f} {\dee x} =  \sin \left( \frac 1 x \right) + 2 x\cos \left( \frac 1 x \right),\] equals zero, which happens with increasing frequency as $x \to 0^+$, and in fact an infinite number of times in the neighborhood of $x=0$.  However, because of the $x^2$ term, since \[\left|\sin\left(\frac 1 x \right)\right|,\left|\cos\left(\frac 1 x \right)\right| \leq 1,\]
we know that \[\int_0^2|f'(x)| \dee x = \int_0^2 \left| \sin \left( \frac 1 x \right) + 2 x\cos \left( \frac 1 x \right)\right| \dee x < \int_0^2 1 + 2x(1) \dee x = 6 < \infty.\]
In other words, because the contribution to $\tvar(f)$ of each critical point goes to zero as the points get closer to zero, the infinite sum of their total variation is convergent, and so the total variation is finite. \\

While compact continuous functions can fail to be tame, non-tame functions can also exist if the domain $D$ is not compact. When $D$ is an open (or clopen) interval, total variation and total persistence are only defined when the limits of $f$ exist at all endpoints. In this case, each is evaluated as though it were continuously extended to the $\bar D$, the closure of its domain, as total persistence and total variation of $f$ and the extended $\tilde f$ can be shown to be equal; see Proposition \ref{prop:cont_extn} in Appendix \ref{subsec:app_Open and Unbounded Domains} for more details. 

The final, most general, way a function can fail to be tame is if its domain is unbounded, and such functions can be either continuous or discrete. As in the case of the open interval, the functions must be asymptotically convergent; in any direction in which $D$ is unbounded, this means that $\lim_{x \to \pm \infty} f(x)$ must exist. 

Indeed, if such an unbounded function \textit{is} tame, it can be transformed into a (tame) function on an open interval by an infinite sequence of \textit{partial} positive horizontal scalings, under which total variation, total persistence, and the persistence diagram are invariant; for details and proof, see Corollary \ref{cor:horiz_scal_local} in Appendix \ref{sec:app_Diagonal Affine Transformation Theorem Cases}. Namely, for a function $f:\rr \to \rr$ and constants $0<a<<b$, we define 
\[g(x) = \begin{cases}
    f(x), & -b+a<x<b-a, \\
    f(\sign(x)\frac {18}a(b-a)(|x|-b+a+\frac a{18})) & b-a \leq |x| < b-\frac a2,\\
    f( \sign(x) \frac {360}a(b-a)(|x|-b+\frac a2 + \frac a{36})) & b-\frac a2 \leq |x| < b-\frac a4,\\
    f( \sign(x) \frac {7200}a(b-a)(|x|-b+\frac a4 + \frac a{72})) & b-\frac a4 \leq |x| < b-\frac a8,\\
    \vdots & \vdots
    
\end{cases} \]
\[= \begin{cases}
    f(x) & -b+a\leq x\leq b-a,\\
    f(\sign(x) \frac {18}a (20)^{\lfloor \log_2(\frac a{b-|x|})\rfloor}(b-a)(|x|-b+\frac{19}{18}a(\frac 12)^{\lfloor \log_2(\frac a{b-|x|})\rfloor})) & b-a < |x| < b
\end{cases}.\]

More generally, a function on an unbounded domain may have infinitely many critical points; that is, not be tame. If such a function is discrete, it must have bounded variation; that is, 
\[ \sum _{i=-\infty}^\infty |f(x_i)-f(x_{i-1})|<\infty.\] 
This condition is necessary and sufficient for the calculation of total variation and total persistence.

Conversely, so long as the function $f$ is \textit{absolutely continuous}, we can show that its total persistence is still equal to its total variation; see the proof of Theorem \ref{thm:abs_cont_eq} in Appendix \ref{subsec:app_Open and Unbounded Domains} for more explanation. 

Absolute continuity is important here is because it is necessary and sufficient \cite[Thm. 7.20, p. 148]{Rudin1987} to evaluate an arbitrary
\[\int_{k_1}^{k_2} f'(x) \dee x\]
as 
\[f(k_2)-f(k_1),\]
without which the formulation of
$\tvar(f)$ as 
\[ \int_D |f'(x)| \dee x\]
does not hold (see Theorem \ref{thm:abs_cont_eq}). Absolute continuity allows us to evaluate the general formulation of total variation \cite[p.126, def. 4.3.2]{Athreya2006}
\[ \tvar(f) = \sup_Q \sum_{i=1}^n |f(x_{i})-f(x_{i-1})|\]
as 
\[ \tvar(f) = \sup_Q \left [ \int_a^{x_1} |f'(x)|\dee x + \ldots +  \int_{x_{n-1}}^{b} |f'(x)|\dee x \right]= \int_a^b |f'(x)| \dee x. \]
However, absolute continuity is not required for $\tvar(f) = \tper(f)$ to hold; for example, the Cantor function (see \cite[sec. 4.5.3, p. 136]{Athreya2006}) has total persistence and total variation equal to one, but is not absolutely continuous, as 
\[ \int_0^1|\mathcal{C}'(x)| \dee x = 0 \neq 1 = \tvar(\mc C).\]
Thus, in the continuous case, bounded variation is necessary and sufficient for total variation and total persistence to exist, but it is not sufficient for the convenient formulation
\[\tvar (f) = \int_a^b|f'(x)| \dee x,\] merely necessary. Absolute continuity is necessary and sufficient for this formula to hold, but merely sufficient (not necessary) for total variation and total persistence to exist and to be calculable.

In Theorem \ref{thm:abs_cont_eq}, we generalize the formula for total variation and total persistence to 
\[\sum_{x \in Crit(f)}M_f(x)(1+I^{Interior}_f(x))f(x),\]
where
\[I_f^{Interior}(x) = \begin{cases}
    1 & \text{ if } x \in \text{Int} (f) \\
    0 & \text{ if } x \in \del(f)
\end{cases}\]
and 
\[M_f(x) = \begin{cases}
    1 & \text{ if $x$ is a local maximum of $f$} \\
    -1 & \text{ if $x$ is a local minimum of $f$}
\end{cases}.\]

We note that the absolute continuity condition does not require that $f$ be defined over a bounded domain. If $f$ is defined on the real line $\rr$ or a half-line, total variation can still be computed as long as $f(x)$ approaches a limit as $x\to \pm \infty$. For more detail, see Appendix \ref{subsec:app_Open and Unbounded Domains}. 

\section{Permutation Equivalence}
\label{sec:Permutation Equivalence}
\subsection{Permutation Theorem}
\label{subsec:Permutation Theorem}
For any two %
continuous functions $f,F : \bb I \to \rr$ on the interval, we let $c_{i,(j)} = (x_i,y_{(j)})$ denote the $i$th local extremum of $f$, and $C_{i,(j)} = (X_i, Y_{(j)})$ denote the $i$th local extremum of $F$. In this section, $x$ and $y$ are indexed separately, such that $x_i<x_j, y_{(i)}<y_{(j)} \iff i < j$. We may assume without loss of generality that critical points of $f$ and $F$ are isolated; if not, we can substitute $\tilde f$ and $\tilde F$ as explained in Corollary \ref{cor:equiv_rel} 
in Appendix \ref{sec:app_Generalizing Equivalence Domain}. Furthermore, although these proofs will use discrete functions for simplicity, as we showed in Section \ref{subsec:Discrete Functions} the persistence diagrams of continuous functions and the discrete sequence of its critical points are identical, so all of these results will hold for continuous functions on the interval as well. 

The object of this section is to show that the property shown in \eqref{lem:PermThm1} is sufficient to conclude that the \textit{ordering} of critical values in the persistence diagrams $\dgm(f)$ and $\dgm F$ are identical. That is, if, say, the second smallest and second largest critical values of $f$ are paired in $\dgm(f)$, then the second smallest and second largest critical values of $F$ will be paired in $\dgm(F)$. Although this is a sufficient condition, it is not necessary, as we will show in Section \ref{subsec:Persistence_Diagram_Preimage}.

First, we show that the condition in \eqref{lem:PermThm1} implies that the ordinal sequence of critical values matches between $f$ and $F$. That is, if, say, the critical value of the fourth critical point of $f$ is the global maximum of $f$, then the corresponding critical value of the fourth critical point of $F$ will be the global maximum of $F$.

\begin{lemma}
\label{lem:PermThm1}
    Let functions $f$ and $F$ have $n$ critical points. Suppose that, 
    \begin{equation} 
    \label{eqn:PermThmLemmaIf}
    \forall i,j \in I_n, \ f(x_i)<f(x_j) \iff F(X_i)<F(X_j)
    \end{equation} That is, there exists a single permutation $\pi$ on $n$ elements, such that permuting the sequences of critical values of $f$ and $F$ by $\pi$ will produce two monotonically increasing sequences. If $f$ and $F$ are generic, this permutation will be unique, and the monotonicity will be strict.

    Then \begin{equation*}
    \label{eqn:PermThmLemmaThen}
    F(X_a) = Y_b \iff f(x_a) = y_b \ \forall a,b \in I_n.
    \end{equation*}
\end{lemma}

\begin{proof}

We suppose first that $F(X_a) = Y_{(b)}$, and we want to show that $f(x_a) = y_{(b)}$. \\

First, we note \eqref{eqn:PermThmLemmaIf} directly implies that 
\begin{equation} 
\label{eqn:LmPrmThmEqlty} 
\forall i,j \in I_n, \ f(x_i)=f(x_j) \iff F(X_i)=F(X_j), \end{equation}
as $f(x_i) = f(x_j)$ is the only condition under which $f(x_i) \not<f(x_j)$ and $f(x_i) \not> f(x_j)$. 

Next, applying \eqref{eqn:PermThmLemmaIf} to each pair of points $X_i,X_j$ which includes $X_a$ in either position, we can produce the set $S_{b-1} \subset I_n$ of indices $k$ whose value $f(x_k) < f(x_a)$. Let us define
\[T_{b-1} = \{j : F(X_j)=Y_{(b-1)}\}, \] the set of indices of the preimage of the next-greatest critical level. We once again apply \eqref{eqn:PermThmLemmaIf} to $X_i,X_j$, where each $i \in S_{b-1} \setminus T_{b-1}$ and $j \in T_{b-1}$. Using \eqref{eqn:LmPrmThmEqlty}, we can see that $f(x_i) = y_{(b-1)} \iff i \in T_{b-1}$, by simply comparing all pairs of indices in $T_{b-1}$. 

We will use this as a base case, and construct an inductive argument to continue this process $b-2$ additional steps to construct $T_1$.

Suppose we are given a set $S_{b-k}$, containing the indices of all $X_i$ for which $F(X_i) < Y_{(b-k+1)}$. As above, we define
\[T_{b-k} = \{j : F(X_j)=Y_{(b-k)}\}, \] 
and apply \eqref{eqn:PermThmLemmaIf} to $X_i,X_j$, where each $i \in S_{b-k} \setminus T_{b-k} := S_{b-k-1}$ and $j \in T_{b-k}$. From this, we get that $f(x_i) = y_{(b-k)} \iff i \in T_{b-k}$.

When we iterate this process until $k = b-1$, we can partition $S_{b-1}$ into $b-1$ sets $T$, each of which correspond to a unique critical level less than $f(x_a)$. As $S_{b-1}$, the set of all indices with critical values less than $f(x_a)$, is partitioned into precisely $b-1$ sets $T$, which each correspond to a distinct critical level $y$ less than $f(x_a)$, we know that $f(x_a)$ is greater than precisely $b-1$ critical values $y$, and thus $f(x_a) = y_{(b)}$.

As there is no inherent distinction between $f$ and $F$, the same argument shows that if $f(x_a) = y_{(b)}$, then $F(X_a) = Y_{(b)}$, and we are done. 
\end{proof}

\begin{corollarylem}
\label{cor:PermThm1}
A critical point $(x_a,y_{(b)})$ is a local minimum [maximum] of $f$ if and only if $(X_a,Y_{(b)})$ is a local minimum [maximum] of $F$.  
\end{corollarylem}

\begin{proof}
    By Lemma \ref{lem:PermThm1}, when they exist, \[\sign(f(x_{a+1}) - y_{(b)}) = \sign (F(X_{a+1})-Y_{(b)}),\] and
    \[\sign(y_{(b)}-f(x_{a-1})) = \sign (Y_{(b)}-F(X_{a-1})).\]

As these determine whether a critical point is a local minimum or maximum, $(X_a, Y_{(b)})$ and $(x_a,y_{(b)})$ will both be either a local minimum or a local maximum. 
\end{proof}

\begin{lemma}
\label{lem:PermThmOpp}
    Let functions $f$ and $F$ have $n$ critical points. Suppose that, 
    \begin{equation} 
    \label{eqn:PermThmOppIf}
    \forall i,j \in I_n, \ f(x_i)<f(x_j) \iff F(X_{n+1-i})<F(X_{n+1-j})
    \end{equation} That is, there exists a permutation $\pi$ on $n$ elements, such that permuting the sequences of critical values of $f$ and $F$ by $\pi$ will produce one monotonically increasing sequence and one monotonically decreasing sequence. If $f$ and $F$ are generic, this permutation will be unique (up to which sequence is permuted to increase and which to decrease), and the monotonicity will be strict.

    Then \begin{equation*}
    \label{eqn:PermThmLemmaOppThen1}
    F(X_a) = Y_{(b)} \iff f(x_{n+1-a}) = y_{(b)}.
    \end{equation*} and, by the symmetry between $f$ and $F$,  \begin{equation*}
    \label{eqn:PermThmLemmaOppThen2}
    f(x_a) = y_{(b)} \iff F(X_{n+1-a}) = Y_{(b)} .
    \end{equation*}
\end{lemma}

\begin{proof}

We suppose first that $F(X_a) = Y_{(b)}$, and we want to show that $f(x_{n+1-a}) = y_{(b)}$. \\

we repeat the proof of Lemma \ref{lem:PermThm1} under the substitution $x_i \mapsto x_{{n+1-i}}$, which is valid since (\ref{eqn:PermThmOppIf}) plays the role of (\ref{eqn:PermThmLemmaIf}) under this substitution This yields that $f(x_{n+1-a}) = y_{(b)}$. 

We switch $f$ and $F$ and repeat to yield the forward $\implies$ implication for the second statement. Then for each statement, we substitute $a$ with $n+1-a$, which yields the backwards $\impliedby$ implication for each statement, as $n+1-(n+1-a) = a$. 
\end{proof}

Finally, we can show that the condition in \eqref{lem:PermThm1} implies this ordinal equality of $\dgm(f)$ and $\dgm(F)$. 

\begin{theorem}[Permutation Theorem] 
\label{thm:PermThm}
    Let functions $f$ and $F$ have $n$ critical points. Suppose that
    \begin{equation}
    \label{eqn:PermThmSuppose}
    \forall i,j \in I_n, \ f(x_i)<f(x_j) \iff F(X_i)<F(X_j).
    \end{equation} %

    Then, \begin{equation}
    \label{eqn:PermThmThen}
    (y_i, y_j) \in \dgm(f) \iff (Y_i, Y_j) \in \dgm(F).\end{equation} Moreover,
    \begin{enumerate}
        \item $(y_i, y_j) \in \ord(f) \iff (Y_i, Y_j) \in \ord(F)$
        \item $(y_i, y_j) \in \rel(f) \iff (Y_i, Y_j) \in \rel(F)$
        \item $(y_i, y_j) \in \ess(f) \iff (Y_i, Y_j) \in \ess(F)$
    \end{enumerate}
    However, the converse is false. 
\end{theorem}

\begin{proof}
First, to provide a counterexample for the converse implication, we can define two functions $f_1, f_2$ with three critical points, given by the mappings
\[f_1: \begin{pmatrix}
    x_1 & x_2 & x_3 \\
    y_{(1)} & y_{(3)} & y_{(2)}
\end{pmatrix}; \text{ and } f_2: \begin{pmatrix}
    x_1 & x_2 & x_3 \\
    y_{(2 )}& y_{(3)} & y_{(1)}
\end{pmatrix}\]
we have that 
\[\dgm(f_1) = \{(y_{(1)},y_{(3)}), (y_{(2)}, y_{(3)}) \} = \dgm(f_2).\]

However, clearly 
\[f_1(x_1) < f_1(x_3)\]
yet 
\[ f_2(x_1) \not< f_2(x_3),\]
contradicting our initial condition \eqref{eqn:PermThmSuppose}. Thus, merely having the same persistence diagram does not mean that two functions have the same ordering of critical values. \\

   \paragraph{$\ess(f)$:} First, let us suppose that a point $(y_{(i)},y_{(j)}) \in \ess(f)$, for two integers $i,j \in I_n$. Then, as $f$ is a function on the interval, $i$ must equal one, and $j$ must equal $m \leq n$, the number of unique critical levels of $f$, which equals $n$ if and only if $f$ is generic; this makes $(y_{(1)}, y_{(m)})$ the only point in $\ess(f)$. 
   
   As noted in \eqref{eqn:LmPrmThmEqlty}, the equality of any two indices \[f(x_i) = f(x_j) \iff F(X_i) = F(X_j)\] 
   follows from the assumed condition \eqref{eqn:PermThmSuppose}. 
   From here, it is straightforward to show that $F$ has the same number of unique critical levels as $f$, and therefore that $(Y_{(1)}, Y_{(m)}) $ is the sole point in $ \ess(F)$. Thus, \eqref{eqn:PermThmThen} holds for $\ess$.

 \paragraph{$\ord(f)$:} Next, let us suppose that $(y_{(i)},y_{(j)}) \in \ord(f)$. Because death times occur after birth times, we know that $i<j$.\\

 We will employ an inductive argument, to show that the sequence of subsets of critical values of $f$ and $F$ will be born and die at the same points. More precisely, we want to show that
 \begin{equation}
 \label{lab:induction_rule}
     \forall i \in I_n, \ \forall k \in \rr, \  \ y_{(i)} \in f^{-1}(-\infty, k] \iff Y_{(i)} \in F^{-1}(-\infty,k].
 \end{equation} 
 
 First, consider $f^{-1}(-\infty, y_{(1)}]$ and $F^{-1}(-\infty, Y_{(1)}]$; that is, when $k = y_{(1)}$. By Lemma \ref{lem:PermThm1}, the indices must be identical; that is, if $y_{(1 )}= f(x_a)$, then $Y_{(1)} = F(X_a)$. As these are the only points in the respective sublevel sets, the dual implication \eqref{lab:induction_rule} holds for $(-\infty, y_{(1)}]$. 

 Because the image under any function of the critical points includes any critical values, we know that the set $f^{-1}(y_{(i)}, y_{(i+1)})$ does not contain any critical points. Therefore, we only need to check $k$ as it reaches critical levels to verify that the dual implication \eqref{lab:induction_rule} holds over all $k$. \\

 Now, suppose that \eqref{lab:induction_rule} holds for all $k \in (-\infty, y_{(i-1)}]$. To assess whether it holds on $(-\infty, y_{(i)}]$, we need only look at what changes at the point $y_{(i)}$ itself. 

If $f$ is non-separating, then a level $y_{(i)}$ may have one or multiple local minima and one or more local maxima.

\subparagraph{Minimum} Suppose $y_{(i)}$ has at least one local minimum. Then, a new connected component of the sublevel set of $f$ will be born at every point in $f^{-1}(y_{(i)})$. Similarly, at $Y_{(i)}$ a new connected component will be born at $F^{-1}(Y_{(i)})$ in the sublevel set of $F$. By Lemma \ref{lem:PermThm1}, these levels will have the same inverse. That is, $x_j \in f^{-1}(y_{(i)}) \iff X_j \in F^{-1}(Y_{(i)})$. Therefore, the sublevel sets of each function will include all of the components in the previous (respective) sublevel set, plus one born at each $x_j$ or $X_j$ which is the inverse of $y_{(i)}$ or $Y_{(i)}$. As this will result in the same set of indices added to the respective previous sublevel sets, which were equal by assumption, the equality in \eqref{lab:induction_rule} will continue to hold over all of $(-\infty, y_{(i)}]$. 

\subparagraph{Maximum} Now, suppose $y_{(i)}$ contains at least one local maximum. 

Each local maximum $x_q \in f^{-1}(y_{(i)})$ adds $x_q$ to a connected component on one of its sides, and this causes the death of the younger connected component whenever there is a component on the other side; that is, when $x_q$ is an interior point. The connected components will be those born at the nearest points $x$ from each side of $x_q$. 

However, this can be observed directly from the set of critical points in the previous sublevel set $f^{-1}(-\infty, y_{(i-1)}]$: The rightmost component (when it exists, for $q < n$) can be observed as starting at $x_{q+1}$ and including every $x$ until one reaches the maximal $x_{q+a}$ such that every \[f(x_{q+1}), f(x_{q+2}), \ldots, f(x_{q+a}) < f(x_q) = y_{(i)}.\] 
This connected component shall be named for $\min  \{f(x_{q+1}), f(x_{q+2}), \ldots, f(x_{q+a})\}$. 

The leftmost component (for $q>1$) can similarly be shown to include a set of consecutively indexed critical points, starting from $x_{q-1}$, and continuing sequentially downwards until the smallest $q-b<q$ for which \[f(x_{q-1}), f(x_{q-2}), \ldots, f(x_{q-b}) < f(x_q),\] and will likewise be named for the minimal critical value in the set. 

Finally, for each point $x_q$ we will compare these two set minima. The component $k_{p_+}$ with the larger set minimum (if a second component exists), is said to die at $y_{(i)}$, and the component $k_{p_-}$ with the smaller minimum, subsumes the other and the $y_q$ which connects it. If these components have equal birth times, then the choice of which component dies and which lives is arbitrary, as they are not distinguished in $\dgm(f)$ apart from by birth and death time. This means, for each point $x_q$, we add a point to $\dgm(f)$ born at $y_{p_+}$, and which dies at $y_{(i)}$; that is, the point $(y_{p_+}, y_{(i)}) \in \dgm (f)$. \\

However, because \eqref{lab:induction_rule} held for $f^{-1}(-\infty, y_{(i-1)}]$, we know that the set of critical points in this set has precisely the same indices as the corresponding set for $F^{-1}(-\infty, Y_{(i)}]$, and so the same two components merge, and the component(s) which die at $Y_{(i)}$ has the same ind(ices) as the component(s) which die at $y_{(i)}$. Therefore,  if \eqref{lab:induction_rule} holds for the sublevel sets of $y_{(i-1)}$ and $Y_{(i-1)}$, it must also hold for $y_{(i)}$ and $Y_{(i)}$. 

Therefore, by induction, we have that \eqref{lab:induction_rule} holds for all $k \in \rr$.

In particular, whenever a local maximum is reached at some $y_{(i)}$ or $Y_{(i)}$, the points $(y_{p_+}, y_{(i)}) \in \dgm(f)$ and $(Y_{p_+}, Y_{(i)}) \in \dgm (F)$ in the respective persistence diagrams will have the same indices. 

\paragraph{$\rel (f)$:} Because the duality between ordinary and relative persistence provides that $\rel (f)$ can be computed by negating $\ord(-f)$, we apply the result above to $\ord(-f)$ to prove it for $\rel(f)$; see Lemma \ref{lem:Dgm_parts_equivalence} in Appendix \ref{app:negation} for details. \\

Thus, we have shown that the sublevel sets (and by duality, the superlevel sets as well) for $f$ and $F$ will contain critical points with precisely the same indices. Furthermore, we have shown that this process is sufficient to compute the locations of points on the persistence diagram. Therefore, if we compute the persistence diagrams of $f$ and $F$, they will be the same, the points will have the same pairs of indices, and so we show \eqref{eqn:PermThmThen}, as desired. \\
\end{proof}

\begin{corollarythm}
   Let functions $f$ and $F$ have $n$ critical points. Suppose that
    \begin{equation*}
    \forall i,j \in I_n, \ f(x_i)<f(x_j) \iff F(X_{n+1-i})<F(X_{n+1-j}).
    \end{equation*} 

    Then, \begin{equation}
    \label{eqn:PermThmCorThen}
    (y_i, y_j) \in \dgm(f) \iff (Y_i, Y_j) \in \dgm(F).\end{equation} Moreover,
    \begin{enumerate}
        \item $(y_i, y_j) \in \ord(f) \iff (Y_i, Y_j) \in \ord(F)$
        \item $(y_i, y_j) \in \rel(f) \iff (Y_i, Y_j) \in \rel(F)$
        \item $(y_i, y_j) \in \ess(f) \iff (Y_i, Y_j) \in \ess(F)$
    \end{enumerate}
    However, the converse is false.  
\end{corollarythm}

\begin{proof}
    As in the previous case, because $f$ and $F$ are continuous functions with $n$ isolated critical points, $\ess(f) = \{(y_1,y_m)\}$ and $\ess(F) = \{(Y_1,Y_m)\}$, where $m \leq n$ is again the number of distinct critical values of $f$, or equivalently of $F$. 

    The proof for $\ord(f)$, and thus also for $\rel(f)$, is identical to the $\ord$ section of the previous Theorem \ref{thm:PermThm}, with the single change that the index of each $X_i$ is replaced by $X_{n+1-i}$. As this does not affect the adjacency of any two points, the same argument will result in the same conclusion: As the sublevel sets are traced, each time a local maximum $y_i$ is reached in $f$ (or equivalently in $F$), the points $(y_{p_+}, y_i)$ and $(y_{p_+}, Y_i)$ are added to $\dgm(f)$ and $\dgm(F)$, respectively; as points added always have the same coordinates, \eqref{eqn:PermThmCorThen} is reached for $\ord$, and then $\rel$ by duality. 

    The (incorrect) converse statement, that two functions which have the same persistence diagram must have reversed critical value orderings, can be shown by identifying $f=F$ for any non-symmetric function $f$; that is, a function of $n$ critical points for which there exists some positive integer $i \leq n$ such that $f(x_i) \neq f(x_{n+1-i})$. 
\end{proof}

\subsection{Persistence Diagram Computation}
\label{subsec:Persistence Diagram Computation}

Because we have established that any two functions with the same sequence of ordinal critical values, will have the same ordinal pairings in their persistence diagram, one may want to extrapolate to compute the persistence diagram, again up to ordinal numbers, of an arbitrary permutation. 

Let a generic function $h$ on the interval have $n$ (distinct) isolated local extrema, $C = \{c_i : i \in I_n \} = \{(x_i, y_{(j)}) : i,j\in I_n\}$, where $x_1<x_2<\ldots x_n$ and $y_{(1)}<\ldots<y_{(n)}$. Define $X = \{x_i : i \in I_n\}$ and $Y = \{y_j : j \in I_n\}$-- that is, the preimage and image (respectively) of the nondegenerate critical points of $h$-- such that $C \subset X \times Y$. 

Then, we define a \textit{discrete} bijection $f : I_n \to I_n$ by $ f(i) = j \iff h(x_i) = y_j$, for all $i,j \in I_n$. For clarity, we will notate the discrete restriction of $h$ to its set of critical points as $F:= h|_{Crit(h)}$. We note, as shown in Section \ref{subsec:Discrete Functions}, that $\dgm(h) = \dgm(F)$, as $h$ is a (monotonic) interpolation of $F$ by definition. 

\begin{figure}[h]
    \centering
    \begin{tikzpicture}[>=stealth]

        \node (I) at (0, 3) {$I$};
        \node (X)  at (3, 3) {$X$};
        \node (Y)  at (3, 0) {$Y$};
        \node (R1) at (0,0) {$\rr$};
        \node (In1) at (6,3) {$I_n$};
        \node (In2) at (6,0) {$I_n$};
        \node (In-1) at (8,-2) {$I_{n-1}$};
        \node (In3) at (11,-2) {$I_n$};
        \node (i_j) at (9.5,2) {$\{(i,j)\}$};
        \node (yi_yj) at (9.5,5) {$\{(y_i,y_j)\}$};

        \node (R2) at (-1, 4) {$\rr$};
        \node (x_set)  at (3, 5) {$I_n$};
        \node (y_set)  at (3, -2) {$I_n$};
        \node (dgm_h) at (12.5,5) {$\dgm(h)$};

         \draw[->] (I) -- node[above] {CritPt()} (X);
         \draw[->] (I)  -- node[right] (h_label) {$h$} (R1);
         \draw[->] (R1) -- node[below] {CritVal()} (Y);
         \draw[->] (X) -- node[below] {$\chi^{-1}$} (In1);
         \draw[->] (Y) -- node[above] {$\psi^{-1}$} (In2);
         \draw[->] (In1) -- node[right] {$f$} (In2);
         \draw[->] (In2) --node[above right] {$\setminus \{n\}$} (In-1);
         \draw[->] (In-1) --node[above] {$g$} (In3);
         \draw[->] (i_j) --node[left] {$\psi$} (yi_yj);

         \draw[->, bend left=20, dotted, thick] (h_label) to (X);
         \draw[->, bend right=20, dotted, thick] (h_label) to (Y);
         \draw[->, dotted, thick, bend right=20] (In-1) to node[left] {$i$} (i_j);
         \draw[->, dotted, thick, bend left=20] (In3) to node[left] {$\times$} node[right] {$j$}  (i_j);
         \draw[->, dotted, thick] (X) -- node[right] {$F$} (Y);

         \draw[dotted, thick] (R2) -- node[above right] {$\supsetneq$} (I);
         \draw[<-] (Y) -- node[left] {$\psi$} (y_set);
         \draw[<-] (X) -- node[right] {$\chi$} (x_set);
         \draw[->] (yi_yj) -- node[above] {=} (dgm_h);

    \end{tikzpicture}
    \caption{The continuous function $h$ is filtered down to its critical points, and the indices of the paired $x_i$ and $y_j$ are used to create the discrete function $f$. The pairing map $g$ is used to generate the set of indices of the paired critical values in $\dgm(h)$, for which an index $g(i)$ exists for every $i<n$. $F$ is the restriction of $h$ to its critical points, and $\psi^{-1} \circ F$ = $f \circ \chi^{-1}$. }
    \label{fig:PerDiagComp}
\end{figure}
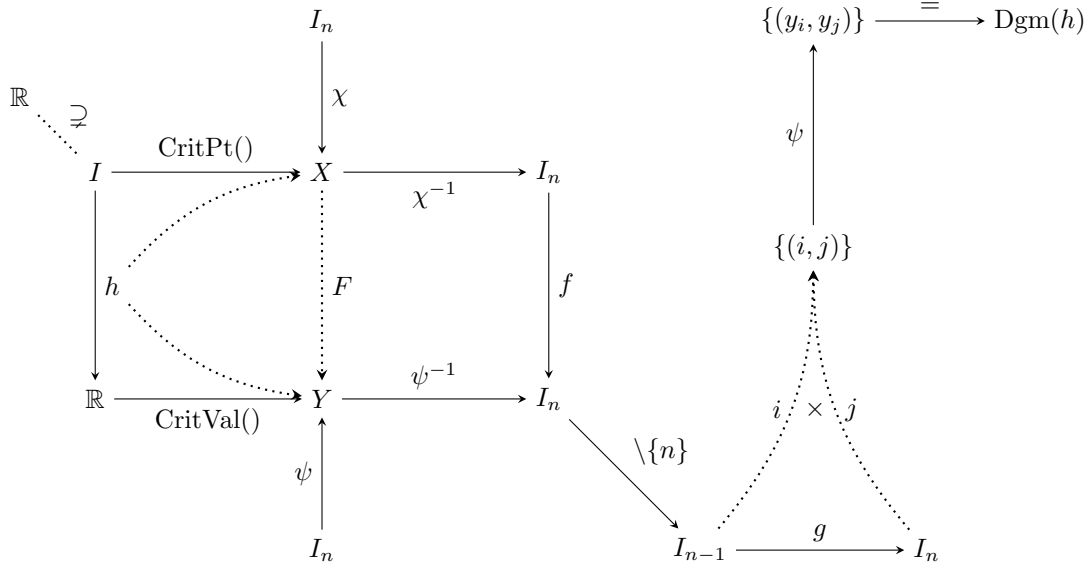

Finally, since we wish to move between the sets $X,Y$ and their indices within the set, we will define $\chi: I_n \to X$ by $i \mapsto x_i$, and $\psi : I_n \to Y$ by $j \mapsto y_j$. Because each critical value is unique in $Y$ (and obviously in $X$ too, as $h$ is a function), we know that $\chi$ and $\psi$ are bijections. The relationships between these sets and functions is illustrated in Figure \ref{fig:PerDiagComp}.\\

We wish to characterize a function $g : I_{n-1} \to Y$, such that, for all $i \in I_{n-1}$, the point $( y_i, y_{g(i)})) \in \dgm(h)$. 

Because $h$ is generic, there is exactly one point in each column of $\dgm(h)$ for every point in $I_{n-1}$, and so $g$ is well-defined. As $y_n$ is the global maximum level, there is no point in $\dgm(h)$ in its column, because it is not recorded by the relative superlevel persistence. 

The image of $g$ will always be \[Y \setminus\{ y_i : i \neq n, \text{ and }(F)^{-1}(y_i) = x_1 \text{ or } x_n, \};\]
that is, the set of critical levels which are the global maximum, or which correspond to interior critical points of $h$. This set will have $n-2$ elements if the local maximum is an interior point, and $n-1$ elements if it is a boundary point. \\

The essential class is identical for all functions on the interval, namely $\ess(h) = \{( y_1,y_n) \}$. It remains to characterize $\ord(h)$ and $\rel(h)$: \\

We let $A:= \{F(x_1), F(x_3), \ldots, F(x_{n\text{ or } n-1})\}$ and $B:=\{F(x_2), F(x_4),\ldots,F(x_{n\text{ or } n-1})\}$, indexed such that $a_1<a_2<\ldots<a_{\lfloor\frac n2 \rfloor-1}$ and $b_1<b_2<\ldots < b_{\lceil \frac n2 \rceil-1}$; that is, the sets of local minimum and maximum values (not necessarily respectively). We define constants $k_- := \lfloor \frac n2 \rfloor -1 = |A|-1$, and $k_+ := \lceil \frac n2 \rceil -1 = |B|-1$.  Finally, we define sets 
\begin{itemize}
\item $A^i := A \setminus \{F(x_1),F(x_n)\}$,
\item $A^\ell := A \setminus \{y_{(1)}, y_{(n)}\}$;
\item $B^i := B \setminus \{F(x_n)\}$, and 
\item $B^\ell := B \setminus \{y_{(1)}, y_{(n)}\}$,
\end{itemize}as the restrictions of $A$ and $B$ to interior $(i)$ and nonglobal $(\ell)$ critical points, respectively.  

Suppose that the first critical point of $h$ is a local minimum, meaning that $A$ is the set of interior local minima, and $B$ the set of interior local maxima. Then, there exist (bijective) permutations maps $\pi_A: A^i \to A^i$ and $ \pi_B: B^i \to B^i$ specific to $h$ (or equivalently, $F$) such that 
\begin{equation}
\label{eq:comp_pers_minA}\ord (h) = \{ ( a^\ell_1, \pi_B \circ h(x_2)), ( a^\ell_2, \pi_B\circ h(x_4)), \ldots, ( a^\ell_{k_+}, \pi_B\circ h(x_{2k_+}))\}\end{equation} and 
\begin{equation}
\label{eq:comp_pers_minB}\rel (h) = \{ ( b^\ell_1, \pi_A \circ h(x_3)), ( b^\ell_2, \pi_A\circ h(x_5)), \ldots, ( b^\ell_{k_-}, \pi_A\circ h(x_{2k_-+1}))\}.\end{equation}

Conversely, if the first critical point of $f$ is a local maximum, then the formulas are reversed, because of the duality of $\ord$ and $\rel$. In other words, we have 
\begin{equation}
\label{eq:comp_pers_maxA}\ord (h) = \{ ( b^\ell_1, \pi_A \circ h(x_3)), ( b^\ell_2, \pi_A\circ h(x_5)), \ldots, ( b^\ell_{k_-}, \pi_A\circ h(x_{2k_-+1}))\}\end{equation}and 
\begin{equation} 
\label{eq:comp_pers_maxB}\rel (h) = \{ ( a^\ell_1, \pi_B \circ h(x_2)), ( a^\ell_2, \pi_B\circ h(x_4)), \ldots, (a^\ell_{k_+}, \pi_B\circ h(x_{2k_+}))\}.\end{equation}

There is not currently a known formula for $\pi_A$ and $\pi_B$ in general, but it has been computed for at least $n \leq 6$, as shown in Appendix \ref{sec:app_Persistence Diagram Computed Examples}. The current computations of $\pi_A$ and $\pi_B$ apply to all generic discrete (or continuous) functions on the interval with finitely many critical points. Future research may attempt to calculate cases for $n \geq 7$, or else to generalize to functions on the circle, or on other domains.

\subsection{Persistence Diagram Preimage}
\label{subsec:Persistence_Diagram_Preimage}

Another way to attempt to characterize the types of functions which share a persistence diagram, is to try to find connections from the bottom up. In general, when $n$ is not small ($n >> 5$ or so), a persistence diagram will correspond to multiple ordinal permutations. 

\begin{proposition}[Swap]
\label{prop:Swap}
    
Let $f : I_n \to \rr$ be a discrete time series. 
Suppose there exists an $i \in \{2,\ldots, n-3\}$ such that 
\begin{itemize}
    \item if $x_{i+1}$ is a local maximum, then $f(x_{i+1})<f(x_{i-1}),f(x_{i+3})$.
    \item if $x_{i+1}$ is a local minimum, then $f(x_{i+1})>f(x_{i-1}),f(x_{i+3})$.
\end{itemize}
Then we define a ``persistence transposition'' $\tau _i$ such that $\tau_i(f)(x_j) := \begin{cases} f(x_j) \text{ if } j \notin \{i,i+2\}, \\ f(x_{2i+2-j}) \text{ if } j \in \{i,i+2\} \end{cases}$. 

Let $h = \tau_i (f)$. Then, $\dgm(f) = \dgm(h)$.
\end{proposition}

\begin{figure}
    \centering
    \includegraphics[width=\linewidth]{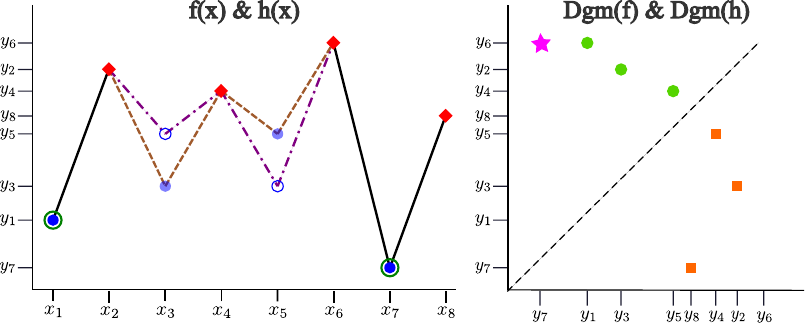}
    \caption{Functions $f$ and $h$ have switched values for $f(x_3)$ and $f(x_5)$, but are otherwise identical. Per Proposition \ref{prop:Swap}, this swap does not change the persistence diagram: Neither $x_3$ nor $x_5$ are boundary points, and $f(x_4)<f(x_2),f(x_6)$ and likewise for $h$. }
    \label{fig:Swap}
\end{figure}

\begin{proof}
    We will prove the equality of these persistence diagrams by construction. Because the critical values themselves are not changed, $\dgm(f)$ and $\dgm(h)$ each have one essential point at $(y_1,y_n)$. 

    We assume, without loss of generality, that $f(x_i) \neq f(x_{i+2})$, as otherwise the functions are identical, and so therefore must their persistence diagrams be.

    We further assume, without loss of generality, that $f(x_i)$ is a local minimum. If it is a local maximum, then we consider instead $-f$ and $-h$, for which $x_i$ is a local minimum; see Lemma \ref{lem:Dgm_parts_equivalence} in Appendix \ref{app:negation} for details. 

    To construct $\ord(f)$ and $\ord(h)$, we gradually increase our sublevel set threshold $\alpha$. Let us call $k := f(x_{i+1})$, $c^-:=\min\{f(x_i),f(x_{i+2})\}$ and $c^+ := \max \{f(x_i),f(x_{i+2})\}$. For $\alpha < c^-$, $f$ and $h$ are identical, so their persistence diagrams are, too. Each has a component born at $c^-$ and $c^+$, which do not die before $k$, because it is less than $f(x_{i-1})$ and $f(x_{i+3})$. In each case, as the components meet at $k$, the component born at $c^+$ must die there, because each component born at $c^+$ and $c^-$ has not merged with any external components, and so the component born at $c^+$ is younger. 
    
    Then, the component born at $c^-$, which at $\alpha = k$ contains all of $[i,i+2]$, continues until it dies. Because the diagrams so far of $f$ and $h$ are identical, and the function signals are identical for $\alpha > k$, the component born at $c^-$ will die at the same point in $\dgm(f)$ and $\dgm(h)$. Thus, each diagram will have points at $(c^+,k)$, and $(c^-,k_1)$, for some $k_1 > k$. If $c^-$ is the global minimum (after accounting for any tiebreaking rules, if $f$ and $h$ are non-separating), then $(c^-,k_1)$ is the essential point; otherwise both points are in $\ord$.

Finally, we will examine $\rel$. As before, for $\alpha >k$, the superlevel relative persistence will be identical between $f$ and $h$, and so too will their diagrams be. Upon reaching $k$, a new component will be born; we know $k$ is not the global maximum because $k < f(x_{i-1}), f(x_{i+3})$. 

The superlevel component born at $k$ will die at $c^+$, and then the components which include $x_{i-1}$ and $x_{i+3}$ will meet at $c^-$. Because $f$ and $h$ are identical outside of $[x_i, x_{i+2}]$, the age of each of the two components at $c^-$ will be identical between $f$ and $h$. As before, because they cannot merge with any older components, the younger of the two components will die at $c^-$. Thus, both $\dgm(f)$ and $\dgm(h)$ will have points at $(k,c^+)$ and $(k_2,c^-)$, for some $k_2>k$. 

Since the two functions are otherwise identical, all other points will be the same in the two persistence diagrams. Thus, we have shown that $\dgm(f) = \dgm(h)$. 
\end{proof}

\begin{proposition}
Let functions $f, h : I_n \to \rr$ be discrete time series.

Then, $\dgm(f) = \dgm (h)$ if %
$f$ and $h$ can be connected by a finite combination of persistence transpositions and reversals. 
    
\end{proposition}

\begin{proof}
 $\tau_i$ is a \textit{persistence transposition} if it meets the conditions of Proposition \ref{prop:Swap}.

If $f=h$, the equivalence of their persistence diagrams is trivial. If $f$ is the reversal of $h$, their persistence diagrams are identical, as shown in Proposition \ref{prop:hor_reflection}. Finally, if $f$ and $h$ can be connected by a finite number of persistence transpositions or reversals, we can create a sequence of functions $g_1,\ldots, g_z, z\in \nn$, such that $g_1 = f, g_n = h$, and every $g_i, g_{i+1}$ are separated by one persistence transpositions or reversal. By Propositions \ref{prop:Swap} and \ref{prop:hor_reflection}, each $\dgm(g_i)=\dgm(g_{i+1})$, meaning $\dgm(f) = \dgm(h)$ by transitivity.

\end{proof} 
\begin{figure}
    \centering
    \includegraphics[width=\linewidth]{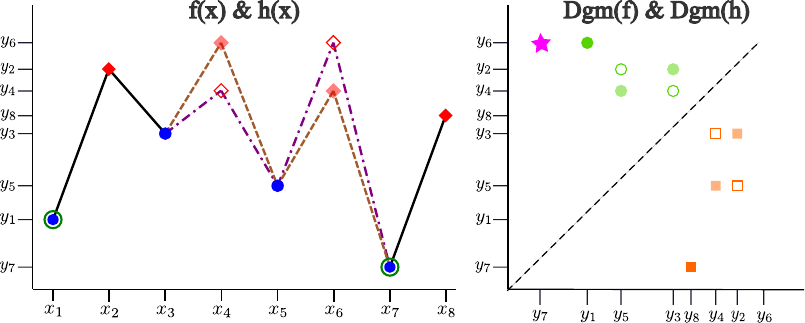}
    \caption{A swap $\tau_4^*$ does not preserve the persistence diagram, because the local maximum $f(x_3) > f(x_5)$, meaning $\tau^*_4$ is not a persistence transposition on $f$ or $h$. The brown evenly [purple unevenly] dotted lines and translucent [outlined] local maxima correspond to the translucent [outlined] points in $\dgm$.}
    \label{fig:Bad_Swap}
\end{figure}
\begin{conjecture}[Persistence Diagram Preimage]
\label{conj:Pers_Diag_Preimage}
    If $f,h$ are discrete time series such that $\dgm(f) = \dgm(h)$, then $f$ and $h$ can be connected by a finite combination of persistence transpositions and reversals. 
\end{conjecture}

  We must show that if $\dgm(f) = \dgm(h)$, then $f$ and $h$ can be connected by a finite number of persistence transpositions and reversals. 

To show this condition is implied, we will prove the contrapositive. We suppose that $f$ and $h$ cannot be expressed as a combination of persistence transpositions and reversals. We must show that this implies that $\dgm(f) \neq \dgm(h)$. 

First, we can assume that the (multi)sets of critical values of $f$ and $h$ are the same, as if these multisets differ, then the persistence diagrams will have different multisets of birth and death times, and therefore will not be equal. 

We can also assume that the (multi)sets $\{f(x_1), f(x_n)\}$ and $\{h(x_1), h(x_n)\}$ are equal, as the critical value associated with each boundary point only appears a single time (as a birth, unless it is the global maximum) in the persistence diagram, unlike interior points which appear twice. If the function is non-separating, the boundary point will still be noticed because each appears once, while every other point appears twice. 

We further assume, without loss of generality, that $f(x_1) = h(x_1)$, as when this is not the case, we substitute $\tilde h$, the reversal of $h$, as we know $\dgm(\tilde h) = \dgm (h)$. Thus we know that, for each $x_i$, $f(x_i)$ and $h(x_i)$ are both local maxima, or both local minima. 

Finally, we assume, again without loss of generality, that the (multi)sets of local maxima and local minima are the same between $f$ and $h$. If this were not the case, then the sets of births in $\ord$ or $\rel$ would differ, meaning $\dgm(f) \neq \dgm(h)$. \\

From here, we observe that, because $f$ and $h$ have the same sets of local minima and local maxima, there exists a permutation $\pi$ which permutes the sequence of local extrema of $f$ to the sequence of local extrema of $h$. Furthermore, as permutations, these can be decomposed into a product adjacency transpositions.

Because $f$ and $h$ cannot be written as a finite combination of persistence transpositions and reversals, every adjacency decomposition $d$ of $\pi$ will contain at least one swap $\tau^*$ which is not a persistence transposition (or a reversal).

If a decomposition $d$ contains either $\tau_1$ or $\tau_{n-2}$, then the set of boundary points will change, and therefore so will the persistence diagram (assuming that $f(x_1) \neq f(x_3)$ or $f(x_{n-2}) \neq f(x_n)$, understanding that otherwise the swap would not change anything). 

If a decomposition $d$ does not contain either $\tau_1$ or $\tau_{n-2}$, it must contain at least one $\tau_i^*$, for which $2 \leq i \leq n-3$, yet which is not a persistence transposition. Since the only additional requirement for a persistence transposition is that $f(x_{i+1})$ be less prominent (extreme) than $f(x_{i-1})$ and $f(x_{i+3})$, it must be that $f(x_{i+1})$ is more prominent than at least one of these. That is, if $f(x_{i+1})$ is a local minimum [maximum], then $f(x_{i+1})$ is less [greater] than at least one of $f(x_{i-1})$ and $f(x_{i+3})$. \\

However, the mere presence of a non-persistence transposition between two functions does not guarantee that their persistence diagrams are not equal. For example, letting $\vec x := [1,\ldots,n]^T$, consider $f(\vec x) = [2,5,4,6,3,8,1,7]^T$ and $h(\vec x) = [2,5,4,8,3,6,1,7]^T$. Though the swap that separates them, $\tau_4^*$, is not a persistence transposition, $\dgm(f) = \dgm(h)$, and indeed 
\[f = (\tau_2 \circ \tau_3 \circ \tau_4 \circ \tau_5 \circ \tau_2 \circ \tau_5 \circ \tau_4 \circ \tau_3 \circ \tau_2)
(h).\]
 Several example preimages are illustrated in Appendix \ref{subsec:app_PD Preimage Illustrations}, with their corresponding persistence transposition connections.

\section{Signal Transformations}
\label{sec:Signal Transformations}

It is not obvious that any formula, such as those developed in Section \ref{subsec:Persistence Diagram Computation}, would be more computationally efficient than simply computing the critical values. However, the fact that the critical points depend solely on $F$ (the function restricted to local extremum points) and not on $h$ (the original function) can be used to show that the persistence diagram $\dgm(F)$ is generally invariant to horizontal transformations, and qualitatively invariant to global vertical transformations which respect ordering. A more complete summary of the effects of such transformations is the subject of the next section. 

\subsection{Diagonal Affine Transformation Theorem}
\label{subsec:Diagonal Affine Transformation Theorem}
The simplest transformations to characterize are \textit{diagonal affine} transformations: that is, transformations that can be decomposed into some number of (horizontal or vertical) translations, scalings, or reflections. When a function undergoes a diagonal affine transformation, we can completely characterize any resulting change to the function's persistence diagram:

\begin{theorem_countless}[Diagonal Affine Transformation]
\label{thm:sgnl_trans_thm} The effect of diagonal affine transformations upon total persistence, total variation, and the persistence diagram of a function $f$ on the interval can be summarized as follows: 
\begin{enumerate}
    \item  Total variation, total persistence, and persistence diagrams are invariant under horizontal shifts, horizontal scaling, and horizontal reflection. 

    \item Total variation and total persistence are invariant under vertical shifts, while the persistence diagram will be shifted diagonally.

    \item Vertical scaling will scale the persistence diagram, the total variation, and the total persistence, by the (positive) scaling constant.
    
    \item Total persistence and total variation are invariant under vertical reflection. Furthermore, each point in Ord and Rel is reflected through the origin, while the point in Ess is reflected across the anti-diagonal line. 
\end{enumerate}
\end{theorem_countless}

The proof of each of these cases can be found in Appendix \ref{sec:app_Diagonal Affine Transformation Theorem Cases}, where it is labeled as Theorem \ref{thm:sgnl_trans_thm}.

\subsection{Function Addition}
\label{subsec:Function Addition}
Many function transformations benefit from the characteristic of \textit{linearity}-- that is, the sum of the transformed functions is the transformation of the sum of the functions. This makes tools like differentiation, integration, and Fourier series immensely powerful. However, the transformation from signal to persistence diagram does not have this property. In general, the persistence diagram of the sum of two functions $f+g$ does not bear any resemblance to the persistence diagrams of either $f$ or $g$, because there is no additive structure on the space of persistence diagrams. Indeed, even the sets of local extrema of $f+g$ are different in general from either $f$ or $g$, because $(f+g)' = 0 \implies f'=-g'$, which neither requires nor is guaranteed by either $f'$ or $g'$ alone being zero. 

If we assume $f$ and $g$ have critical points in the same places\footnote{This is not always a requirement to characterize the persistence diagram, but it generally is for nonconstant $g$.}, then we know that the set of critical points CritPoint($f+g$) $\supset$ CritPoint($f$) = CritPoint($g$), but CritPoint($f+g$) may also have more critical points if $\sign(f'(x)) = -\sign(g'(x))$ at a given $x$; since the critical points of $f$ and $g$ are the same, the signs will switch in the same places, and so either $\sign(f'(x)) = -\sign(g'(x))$ or $\sign(f'(x)) = \sign(g'(x))$ for all $x$. 

If we further assume that $\sign(f'(x)) = \sign(g'(x))$, then at last we know that CritPoint($f+g$) = CritPoint($f$) = CritPoint($g$). However, we will still have, in general, no relation between the sets of critical \textit{values} of $f$, $g$, and $f+g$. \\

There are still some conditions under which persistence diagrams can be `added'' in some sense. For example, if $g$ is a constant function or a scalar multiple of $f$, the persistence diagram of $f+g$ can be computed according to the corresponding transformation in Theorem \ref{thm:sgnl_trans_thm} (Propositions \ref{prop:vert_shift} through \ref{prop:vert_reflect}). 

Furthermore, persistence diagrams can be "added" for functions $g$ with a small enough range, as a consequence of the stability condition (see \cite{CohenSteiner2005}). That is, the bottleneck distance
\[d_B(X,Y) = \inf_\gamma \sum _x \|x-\gamma(x)\|_\infty\]
is bounded from above by the $L_\infty$ norm of the difference between two functions, $\|f-g\|_\infty$. \cite{CohenSteiner2005}

In particular, given the set LE($f$) of local extrema (points) of $f$, define $\eta_f := \min_{a \neq b \in \text{CV}(f)}|f(a)-f(b)|$. Note that $\eta_f = 0 \iff f$ is a non-separating function. 

\begin{proposition}[Addition and Stability]
\label{prop:Addition_and_Stability}
    
Let $f$ be a function, and $F := f|_{\text{LE}(f)}$ be a discrete function which takes the values of the function $f$ restricted to the set  LE($f$) of local extrema of $f$.  

Then, if there exists a function $g$ on the same domain, such that
\[\mathrm{sign} (g'(x_i)) = \mathrm{sign} (f'(x_i)) \ \forall x_i \footnote{This implies that $f$ and $g$ have the same sets of critical points, in addition to having identical parity between them.}\] and
\[\|g\|_\infty < \frac {\eta_f}2,\] 
then
\[ \dgm(f+g) = \{(A + g\circ F^{-1}(A), B+g\circ F^{-1}(B)), \ : \ (A,B) \in \dgm(f)\}.\]
\end{proposition}

\begin{proof}
    If $f$ is not generic, then $\eta_f=0$, and so no such $g$ exists, and the implication is vacuously true. 
    Thus, when $f$ is generic, $F$ is a bijection, and so $F^{-1}$ is well-defined. Moreover, by the Permutation Theorem (Theorem \ref{thm:PermThm}), $f$ and $F$ will have the same persistence diagram.

    Because \[|(f+g)(x_i)-f(x_i)| = |g(x_i)| < \frac {\eta_f} 2,\]
    if 
    \[f(x_i) - f(x_j)>0,\]
     then
     \[(f+g)(x_i) -(f+g)(x_j) = [f(x_i)-f(x_j)] + [g(x_i) - g(x_j)] > [f(x_i)-f(x_j)] - \frac {\eta_f} 2 - \frac{\eta_f}2\]\[ =  [f(x_i)-f(x_j)] - \eta_f \geq 0,\]
     by the definition of $\eta_f$. 

By genericity, we know that if $i\neq j$, $f(x_i) \neq f(x_j)$. Therefore, as the same implication above holds for $f(x_j) - f(x_i)$, we have that
\[f(x_i) < f(x_j) \iff (f+g)(x_i) < (f+g)(x_j),\]
and therefore, by Theorem \ref{thm:PermThm}, $(y_i,y_j) \in \dgm(f) \iff (Y_i,Y_j) \in \dgm(f+g)$. 

However, because we know that the signs of $f$ and $g$ are the same, we know that they have the same set of critical points, and as they have the same sign between critical points, $f+g$ has every critical point of $f$ and $g$ and no more. 

As we showed for Theorem \ref{thm:PermThm}, $f$ and $f+g$ have the same ordering of critical points. Therefore, we know that for every critical point $x_i$, if $f(x_i) = y_j$, then $Y_j = (f+g)(x_i) = y_j + g(x_i)$. As a result, for every $(f(x_i),f(x_j)) \in \dgm(f)$, $(f(x_i) + g(x_i), f(x_j) + g(x_j)) \in \dgm(f+g)$. The unique critical point $x_i$ such that $f(x_i) = y_j$ can be written as $F^{-1}(y_j)$, giving us our decomposition of $\dgm(f+g)$, as desired: 
\[ \dgm(f+g) = \{(A + g\circ F^{-1}(A), B+g\circ F^{-1}(B)), \ : \ (A,B) \in \dgm(f)\}.\]
     
\end{proof}

\subsection{Integration and Differentiation}
\label{subsec:Integration and Differentiation}
Another very common transformation we wish to examine is differentiation and integration. In some applications, it is easier to indirectly measure some desired characteristic (for example, position or velocity of an object) by measuring another quantity, and calculating the desired characteristic by computing derivatives or antiderivatives (in this example, acceleration). 

We will show that there cannot exist a nontrivial invariant of a function $f$, which is invariant under summation / integration, difference / differentiation, the persistence diagram replacement function $D_{h}(f) = (\dgm_h^{-1} \circ \dgm)(f) = \{ h : \dgm(f) = \dgm(h)\}$. In particular, for any two functions $f,g$, there exists a sequence $O_1,O_2,\ldots, O_n$ of these three functions such that $g = (O_n \circ \ldots \circ O_2 \circ O_1) (f)$. 

We start with the mathematically simpler case of discrete data. In discrete calculus, the analogs of differentiation and integration are differences and sums. Then, because any continuous function can be formulated as an interpolation of discrete points, we will build the differential and integral cases from the discrete case.

\begin{figure}[h]
\centering
\begin{tikzpicture}[>=Stealth]

\node (f) at (0,0) {$f$};
\node (F) at (4,0) {$F$};

\node (dgmf) at (0,3) {$\dgm(f)$};
\node (dgmF) at (4,3) {$\dgm(F)$};

\draw[->] (0.3,0.1) -- node[above] {$\sum$ / $\int$} (3.7,0.1);
\draw[->] (3.7,-0.1) -- node[below] {$\Delta$ / $\frac \dee {\dee x}$} (0.3,-0.1);

\draw[->] (0.1,0.3) -- node[right] {$\dgm$} (0.1,2.7);
\draw[dashed, ->] (-0.1,2.7) -- node[left] {$\dgm^{-1}$} (-0.1,0.3);

\draw[->] (4.1,0.3) -- node[right] {$\dgm$} (4.1,2.7);
\draw[dashed, ->] (3.9,2.7) -- node[left] {$\dgm^{-1}$} (3.9,0.3);

\draw[dotted, ->] (0.8,3.1) -- node[above] {?} (3.2,3.1);
\draw[dotted, ->] (3.2,2.9) -- node[below] {?} (0.8,2.9);

\end{tikzpicture}
\caption{The relationship between the persistence diagrams of a function $f(x)$ and its sum or integral $F(x) = \sum_{k=1}^if(x_k)$. The dashed lines represent the preimage of the non-injective function $\dgm$, as in Section \ref{subsec:Persistence_Diagram_Preimage}.}
\label{fig:Disc_Int/Diff_Dgm}
\end{figure}
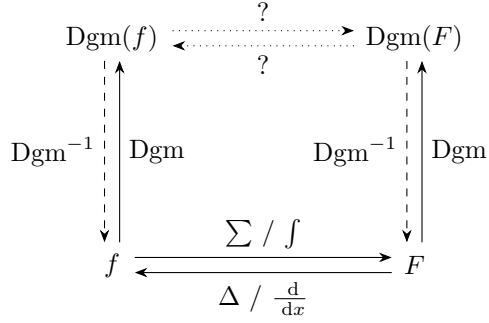

Let $f(x_i), \ i \in I_n$ be a (uniformly sampled) discrete function, with a summation function \[\Sigma : f \to F \ |
\ f(x_i) \mapsto \sum_{k=1}^{i-1} f(x_k).\] We will examine the effects of summation and differences simultaneously, by comparing $f$ and $F$ to functions which share persistence diagrams with them. 

The principal difficulty of characterizing any relationship between $\dgm(f)$ and $\dgm(F)$ is the many-to-one (non-injective) characteristic of the $\dgm$ function, and that its inverse $\dgm^{-1}$ is therefore one-to-many. Thus, na\"ively writing
\[\dgm(F) = \dgm \circ \Sigma \circ \dgm^{-1} (f)\] is not meaningful, as it is set-valued, not any specific function.\\

\begin{figure}
    \centering
    \includegraphics[width=0.7\linewidth]{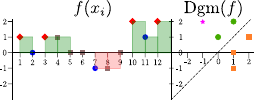}
    \caption{An example discrete function $f$, with its persistence diagram. The colored boxes represent the contribution of each point to the running total, which is calculated as $F(x_i)$ in Figure \ref{fig:IntDiff_FF}.}
    \label{fig:IntDiff_f}
\end{figure}

Generally, when working with non-injective functions and their inverses, it can be helpful to try to work with \textit{equivalence classes}, which are sets that are built up from \textit{equivalence relations}. Given a set $A$ of elements (in this case, discrete 1-D signal functions like $f$), we define an equivalence relation $\sim$ on each pair of elements $a_1,a_2 \in A$, such that the following axioms are satisfied:
\begin{enumerate}
    \item Identity: $a \sim a \ \forall a \in A$,
    \item Symmetry: $a_1 \sim a_2 \iff a_2 \sim a_1 \ \forall a_1,a_2 \in A$, and
    \item Transitivity: if $a_1 \sim a_2$ and $a_2 \sim a_3$, then $a_1 \sim a_3$, $\forall a_1, a_2, a_3 \in A$. 
\end{enumerate}
From this equivalence relation, we can construct equivalence classes $B_1, B_2,\ldots$, which collectively partition the set $A$, and for which, given $b_1 \in B_i$, $b_1 \sim b_2 \iff b_2 \in B_i$. \\

One obvious candidate for an equivalence relation would be to define the equivalence relation $\sim$ such that discrete functions $g_1 \sim g_2 \iff \dgm(g_1) = \dgm(g_2)$. In this case, we can label each equivalence class $\tilde D$ by the persistence diagram $D$ which all of its elements share. 

\begin{figure}
    \centering
    \includegraphics[width=0.7\linewidth]{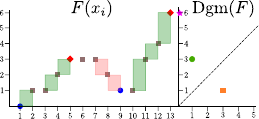}
    \caption{$F(x)$ is the running summation of the function $f$ in Figure \ref{fig:IntDiff_f}. The colored boxes represent the difference function between consecutive values, which traces out the function $f$. }
    \label{fig:IntDiff_FF}
\end{figure}

For example, we consider an arbitrary function $f$, as shown in Figure \ref{fig:IntDiff_f}. The highlighted boxes show the contribution to the running sum $F$, seen in Figure \ref{fig:IntDiff_FF} with $\dgm (F)$. The colored boxes in Figure \ref{fig:IntDiff_FF} represent the difference function of $F$, which is shown in Figure \ref{fig:IntDiff_f} as $f$. 

If we define $D = \dgm(F)$, then we can identify $\tilde D$ as the class of functions $g$ for which $\dgm(g) = D$. The next step is to see whether these equivalence classes are robust under summations or differences. That is, if $\dgm(f) =D = \dgm(g)$, do we know whether, taking the summations $F$ and $G$ of $f$ and $g$ respectively, it follows that $\dgm(F) = \dgm(G)$? What about the differences?

\begin{figure}
    \centering
    \includegraphics[width=0.7\linewidth]{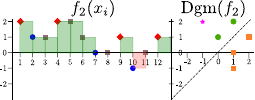}
    \caption{The function $f_2$ was created such that $\dgm(f) = \dgm(f_2)$ (see Figure \ref{fig:IntDiff_f}). As with $f$, the boxes denote the running sum, shown in Figure \ref{fig:IntDiff_FF2}.}
    \label{fig:IntDiff_f2}
\end{figure}

As it turns out, the persistence diagram is preserved under neither sequence of operations. To show summation, we define a function $f_2(x)$ such that $\dgm(f_2) = \dgm(f)$, as seen in Figure \ref{fig:IntDiff_f2}. Inspection confirms that the persistence diagrams are identical between Figures \ref{fig:IntDiff_f2} and \ref{fig:IntDiff_f}.

\begin{figure}
    \centering
    \includegraphics[width=0.7\linewidth]{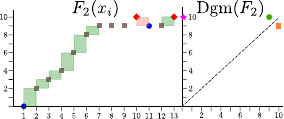}
    \caption{The function $F_2$ is the running sum of the function $f_2$, as seen in Figure \ref{fig:IntDiff_f2}. Contrast the persistence diagram $\dgm(F_2)$ with $\dgm(F)$ in Figure \ref{fig:IntDiff_FF}.}
    \label{fig:IntDiff_FF2}
\end{figure}

Next, we wish to calculate $F_2$ by taking the running sum of the values in $f_2$. We compute this running sum and inspect the persistence diagram, and see in Figure \ref{fig:IntDiff_FF2} that the persistence diagram $\dgm(F_2)$ is different from $\dgm(F)$ in Figure \ref{fig:IntDiff_FF}; note especially the difference in scaling in the former diagram.

\begin{figure}
    \centering
    \includegraphics[width=0.7\linewidth]{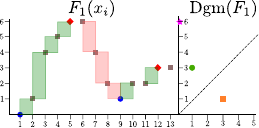}
    \caption{Similarly, $F_1$ is a function created such that $\dgm(F_1) = \dgm(F)$. The colored boxes trace out the function $f_1$, created by taking the first difference of $F_1$ (see Figure \ref{fig:IntDiff_f1}).}
    \label{fig:IntDiff_FF1}
\end{figure}

A similar process can be performed for differences. We define function $F_1$ such that $\dgm(F_1) = \dgm(F)$, as seen in Figure \ref{fig:IntDiff_FF1}. 

\begin{figure}
    \centering
    \includegraphics[width=0.7\linewidth]{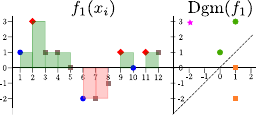}
    \caption{The function $f_1$ is created as the difference of $F_1$ (see Figure \ref{fig:IntDiff_FF1}), and the colored boxes of $f_1$ trace out $F_1$. Contrast $\dgm(f_1)$ with $\dgm(f)$, seen in Figure \ref{fig:IntDiff_f}.}
    \label{fig:IntDiff_f1}
\end{figure}

In Figure \ref{fig:IntDiff_FF1}, we take the first difference between consecutive entries of $F_1$, and compute $f_1$. We compute the persistence diagram $\dgm(f_1)$ and compare to $\dgm(f)$. Though some of the features are the same, such as the points at $(1,0)$ and $(0,1)$, other points have moved, as seen in Figure \ref{fig:IntDiff_f1}. The persistence diagrams do not even have the same number of points. 

The key takeaway from this exercise is that neither summation nor differences preserve the equivalence relation $\sim$, as originally defined, viz. $g_1 \sim g_2 \iff \dgm(g_1) = \dgm(g_2)$. 

Still, this does not mean that these equivalence classes must correspond to individual persistence diagrams. The equivalence classes could, in principle, partition the set of persistence diagrams into a number of classes, with some defining characteristic which any function $f$ shares with every function $g$ which it can reach through a sequence of the above operations. If this defining characteristic were invariant under both the summation / difference operations, and the persistence diagram operations, then two functions for which these characteristics differed would never be able to be transformed into one another by any sequence of the above operations. \\

However, no such characteristic exists with more than one possible class (that is, which is non-trivial). We can show that no non-trivial equivalence relation is invariant under both $\dgm()$ and either summation or differences. In particular, we will show that, for any two discrete signals $f,g$, there exists a sequence of operations among $\Sigma, \Delta, \dgm(),$ and $\dgm^{-1}()$ which connects $f$ to $g$, where the function $\dgm^{-1} \circ \dgm(h_1)$ takes on the value of an arbitrary function $h_2$ such that $\dgm(h_1) = \dgm(h_2)$. If we can connect any two specific functions $f$ and $g$ by a sequence of these four functions, we know at the minimum that they are identified under any equivalence relation $\sim_E$ which would be invariant under these four operations. Since we will show this for all discrete functions, then we have shown that there are no non-trivial functions $f,g$ for which $f \not \sim_E g$, and therefore that $\sim_E$ identifies every function $f$ in the same equivalence class, and is therefore a trivial equivalence relation.

\subsubsection{Calculus Connection Theorem}
\begin{theorem}[Calculus Connection]
\label{thm:Calc_connection}
    Let $f : I_{n_1} \to \zz$ and $g : I_{n_2} \to \zz$ be two discrete functions, expressed as row vectors $[f(i)]_{i=1}^{n_1}$ and $[g(i)]_{i=1}^{n_2}$, $n_1,n_2 \in \zz^+$. 

    Then, there exists a sequence of operations of the following types, such that their composition $O_1 \circ \ldots \circ O_k (f) = g$. 
\begin{itemize}
    \item $ \Sigma_c : [x_1, \ldots, x_k] \mapsto \left[c, c+x_1, c+x_1 + x_2,\ldots,c+\sum_{i=1}^k x_i\right] \in \zz^{k+1}$
    \item $\Delta: [x_1, \ldots, x_k] \mapsto \left[x_2-x_1, \ldots, x_k-x_{k-1}\right] \in \zz^{k-1}$
    \item $D_{\vec y}(\vec x) : \vec x \mapsto  y \in \dgm^{-1} \circ \dgm(\vec x)$. The dimension of $\vec y$ can be arbitrary, but the multiset equality $Crit(\vec y) = Crit(\vec x)$ must hold, as $\dgm(\vec x) = \dgm(\vec y)$.  
\end{itemize}
\end{theorem}

\begin{proof}
    Let $f_1 : I_{n_1} \to \zz$ and $g_1 : I_{n_2} \to \zz$ be arbitrary functions, $n_1,n_2>1 \in \zz^+$. 
    
    We note that each of these operations is invertible: $D_{\vec x} \circ D_{\vec y}(\vec x) = \vec x$, $\Delta \circ \Sigma_c(\vec x) = \vec x$, and $\Sigma_{x_1} \circ \Delta (\vec x) = \vec x$. Therefore, it suffices to find a function $h$ to which we connect both $f_1$ and $g_1$ to show that $f_1$ and $g_1$ are connected to each other. 

First, suppose one (or both) of $f_1(1)$ and $g_1(1)$ are negative or zero. Then, we define 
\[f_2:=\begin{cases} \Sigma_c (f_1), c > |f_1(1)| & \text{ if } f_1(1) \leq 0 \\
f_1 & \text{ if } f_1(1) >0     
\end{cases} \]
and 
\[g_2:=\begin{cases}\Sigma_c(g_1), c > |g_1(1)| & \text{ if } g_1(1) \leq 0 \\
g_1 & \text{ if } g_1(1) > 0\end{cases}.\]

Next, suppose $f_2(1) \neq g_2(1)$. Again, we replace $f_2$ with \[f_3 :=\begin{cases}\Sigma_{g_2(1)} f_2 & \text{ if } f_2(1) \neq g_2(1),\\
f_2 & \text{ if } f_2(1) = g_2(1)
\end{cases}.\] and so $f_3(1) = g_2(1)$. 

We will define $f = f_3$ and $g=g_2$.  Because $g_2(1) > 0$, it follows that $f_3(1) > 0$ as well, meaning that $f(1)=g(1) > 0$ holds.

Because repeated iterations of $\Sigma$ on a function $f$ or $g$ is just a sum of sums, we can compute the $m$th sum $\Sigma_0^m$ of $f = [f(i)]_{i=1}^{n}$ as
\[\Sigma_0^m (f) = \left[\sum_{k=m+1}^i \begin{pmatrix}
    i-k+m-1 \\ m-1
\end{pmatrix} f(k-m)\right]_{i=1}^{n+m},\]
where the sum $\sum_a^b$ of anything evaluates to zero if $a>b$.

In particular, for an arbitrarily large $m$, we know that $\Sigma_0^m f(m+1) = f(1)$ is the first non-zero entry of $\Sigma_0^mf$, and the $n-1$ terms after that will be polynomial in $m$, whose leading coefficient will be proportional to $f(1)$:
\[\Sigma_0^m f(m+2) = mf(1) + f(2); \Sigma_0^m f(m+3) = \frac 12 m(m+1) f(1) + mf(2) + f(3), \ldots.\]
Therefore, since $f(1) > 0$, we can choose a sufficiently large $m$ such that $\Sigma_0^m f(m+z) > 0 \forall z \in I_n$.

Let us take the maximum $m$ of such $m_f,m_g$ defined for the functions $f$ and $g$, respectively. Then, $\Sigma^m f(i)$ and $\Sigma ^m g(i)$ are all positive for any $i$ which comes after the first nonzero entry of each. We denote $F = \Sigma^m f$ and $G = \Sigma ^m g$.

For the next step, we want to ensure that the sum of all entries of $f$ and $g$ are equal. Let $S(h)$ denote the \textit{grand sum} of $h$: that is, the sum of every value $h(i)$ at every $i$ in the domain of $h$. Then, we compare $S(F)$ and $S(G)$. We will add additional monotone points (that is, those which are not local extrema) to at least one of the functions $F$, $G$ to bring their grand sums to be equal. 

We suppose that $S(F) \neq S(G)$, as otherwise we are done with this step. Let $d := S(F) - S(G).$ We wish to add points to $F$ and $G$ such that the sum of the points added to $G$ exceeds that of the points to $F$ by precisely $d$.

Because $S(F) \neq S(G)$, we know there is at least one index $i$ for which $F(i) \neq G(i)$. Because $F(m+1) = G(m+1)$, this means that there is at least one point $j \in [m+1,i)$ for which either $F(j) \neq F(j+1)$ or $G(j) \neq G(j+1)$. 

We detour into number theory. Without loss of generality, we declare that $F(j) \neq F(j+1)$. Let $a=F(j)$, $b=F(j) + \sign (F(j+1)-F(j))$, and $c = G(j)$ (though any index would work for $c$). 

Because $|a-b| = 1$, we know that $a$ and $b$ are relatively prime, since $a \geq 1$ and $b \geq \min(F(j),F(j+1)) \geq 1$. Thus, any integer greater than $ab-a-b$ can be written as a non-negative linear combination of $a$ and $b$; that is, as $k_1a + k_2b$, for $k_1,k_2 \geq 0$. In particular, regardless of the value of $d$, since $c > 0$, we know that, for a sufficiently large $k_3 > 0$, that $k_3 c + d > ab-a-b$, and therefore 
$\exists k_1,k_2,k_3 \in \zz^{\geq 0}$ such that $k_1a+k_2b = k_3c + d$. 

With these numbers in mind, we can create
\[y_F(i) = \begin{cases}
    F(i) &\text{ if } i \leq j \\
    F(j) &\text{ if } j < i \leq j+k_1 \\
    F(j) + \sign(F(j+1)-F(j)) & \text{ if } j+k_1 < i \leq j+k_1+k_2 \\
    F(i-k_1-k_2) & \text{ if } i > j+k_1+k_2\end{cases}\]
and 
\[y_G = \begin{cases}
    G(i) & \text{ if } i < j\\
    G(j) & \text{ if } j \leq i \leq j+k_3 \\
    G(i-k_3) & \text{ if } i > j+k_3
\end{cases}.\]

Because each function $y_F$ and $y_G$ is formed by simply inserting points into $F$ and $G$ (respectively) which respect monotonicity, and therefore do not create any additional local extrema, we know that $\dgm(y_F) = \dgm(F)$ and $\dgm(y_G) = \dgm(G)$, by Theorem \ref{thm:PermThm}. 

Finally, because $S(y_F) = S(y_G)$ by construction, the final point of $\Sigma_0 y_F$ and $\Sigma_0 y_G$ will be the same. Furthermore, because $y_F$ and $y_G$ are strictly nonnegative, $\Sigma_0 y_F$ and $\Sigma_0 y_G$ will be monotonically increasing, so 
\[\dgm(\Sigma_0 y_F) = \{(0, S(y_F))\} = \{(0,S(y_G))\} = \dgm(\Sigma_0 y_G).\]

Thus, we can construct a chain

\begin{figure}[h]
    \centering
    \begin{tikzpicture}[>=stealth] 
        \node (f1)   at (0, 0)   {$f_1$};
        \node (f)    at (2.5, 0) {$f$};
        \node (F)    at (5, 0)   {$F$};
        \node (yF)   at (7.5, 0) {$y_F$};
        \node (S0yF) at (10, 0)  {$\Sigma_0 y_F$};

        \node (g1)   at (0, -2)   {$g_1$};
        \node (g)    at (2.5, -2) {$g$};
        \node (G)    at (5, -2)   {$G$};
        \node (yG)   at (7.5, -2) {$y_G$};
        \node (S0yG) at (10, -2)  {$\Sigma_0 y_G$};

        \draw[->] (f1) -- node[above] {$\Sigma^{0,1,2}$} (f);
        \draw[->] (f)  -- node[above] {$\Sigma_0^m$} (F);
        \draw[->] (F)  -- node[above] {$D_{y_F}$} (yF);
        \draw[->] (yF) -- node[above] {$\Sigma_0$} (S0yF);
        
        \draw[->] (S0yF) -- node[right] {$D_{\Sigma_0 y_G}$} (S0yG);
        
        \draw[->] (g1) -- node[below] {$\Sigma^{0,1}$} (g);
        \draw[->] (g)  -- node[below] {$\Sigma_0^m$} (G);
        \draw[->] (G)  -- node[below] {$D_{y_G}$} (yG);
        \draw[->] (yG) -- node[below] {$\Sigma_0$} (S0yG);
        
    \end{tikzpicture}
    \caption{Because a chain of invertible operations connects $f$ and $g$ to a common function, eg. $\Sigma_0y_G$, we know that $f$ and $g$ are connected to each other. }
    \label{fig:Calc_Con_Thm}
\end{figure}

of invertible operations $\Sigma, \Delta, D$   which connects any $f_1$ and $g_1$. 
\end{proof}

\begin{corollarythm}
    \label{cor:calc_connection1}
    Let $f:I_{n_1} \to \rr$ and $g : I_{n_2} \to \rr$. 
    Then, there exists a sequence of operations of the following types, such that their composition $O_1 \circ \ldots \circ O_k (f) = g$. 
\end{corollarythm}

\begin{proof}
    Let $f_1$ and $g_1$ be discrete functions with real-valued outputs. We define $f_3$ and $g_2$, and thus $f$ and $g$, identically to the proof for Theorem \ref{thm:Calc_connection}. 
    $F$ and $G$ are calculated in the same way. We define $a,b,$ and $c$ identically to above, but because we are working with real numbers, $a,b,c$ and $d$ are not necessarily integers. 
    However, because $|a-b| = 1$, we can choose a sufficiently large $k_3$ such that $\exists k_1,k_2$ such that $k_1a+k_2b < k_3 + d < k_1a + k_2b + 1$, for $k_1,k_2 > 1$. We define $k_4 := k_3c + d - (k_1-1)a - k_2b$, such that $k_4$ is located between $a$ and $b$. 

Then, we create
    \[y_F(i) = \begin{cases}
    F(i) &\text{ if } i \leq j \\
    F(j) &\text{ if } j < i \leq j+k_1 -1 \\
    k_4 &\text{ if } i = j+k_1 \\
    F(j) + \sign(F(j+1)) & \text{ if } j+k_1 < i \leq j+k_1+k_2 \\
    F(i-k_1-k_2) & \text{ if } i > j+k_1+k_2\end{cases}\]
and 
\[y_G = \begin{cases}
    G(i) & \text{ if } i < j\\
    G(j) & \text{ if } j \leq i \leq j+k_3 \\
    G(i-k_3) & \text{ if } i > j+k_3
\end{cases}.\]

This is sufficient to show that $S(y_F) = S(y_G)$, and the rest of the proof continues identically. 
\end{proof}

\begin{corollarythm}
    \label{cor:calc_connection2}
 Let there be continuous functions $f,g:\bb I \to \rr$, with finitely many critical points on $I$. 
    Then, there exists a sequence of operations of the following types, such that their composition $O_1 \circ \ldots \circ O_k (f) = g$. 
\end{corollarythm}

\begin{proof}
    By Theorem \ref{thm:PermThm}, if $f$ has $n_1$ local extrema, and $g$ has $n_2$ local extrema, then we can construct discrete functions $f':I_{N_1} \to \rr$ and $g':I_{n_2}$ such that $f'(i)$ is the $i$th critical value of $f$, and $g'(i)$ is the $i$th critical value of $g$. In particular, these discrete functions will have the same persistence diagram as the continuous functions from which they were created. Thus, the continuous case reduces to the equivalence shown in \ref{cor:calc_connection1}.
\end{proof}

Future work can examine whether a meaningful invariant characteristic can exist among a restriction of the above operations; for example, by requiring that summations / integrations alternate with differences / derivatives, or by requiring summations / integrals to start at zero, rather than another constant. However, merely bounding from above the number of operations cannot serve as an invariant, because transitivity is required within every equivalence class created.

\section{Conclusions and Further Directions}
\label{sec:Conclusions and Further Directions}

The equation of the total variation and the total persistence of a function gives us another way to view the pairing of points in the persistence diagram algorithm: We are divvying up the finite amount of total variation between pairs of local minima and local maxima, and we are searching for those pairs which are allotted a significant proportion of this fixed quantity; these pairs in turn correspond to the large-scale features in the signal, whereas other pairs with less persistence correspond to much less total variation in the signal, and are therefore much more likely to be noise, or otherwise negligible. 

While not a novel realization, the generalization of this equivalence beyond Biswas' generic continuous functions implies that this relation is fundamental to persistence, and not merely an artifact of certain ``niceness" assumed in the function or its domain. 

This generalization extends even beyond the realm of continuous functions, to discrete signals, or time series. In the existing literature, papers tend to focus on one or the other of these settings. However, we have shown that many results and tools can be used in either setting with equal facility. 

In this vein, we have further abstracted the process of computing the persistence diagram to the minimal information required to compute it; that is, the sequence of critical values. Additionally, any such sequences with the same relative orderings will have the persistence diagrams that pair the same ordinal critical values. Our work on signal transformations likewise unified the continuous and discrete settings. \\

Further directions of research in this field are numerous. For example, one could further interrogate Conjecture \ref{conj:Pers_Diag_Preimage}, which claims that any discrete time series with identical persistence diagrams can be connected by a sequence of persistence transpositions and possibly a reversal. While Section \ref{subsec:Persistence Diagram Computation} calculated the permutations $\pi_A$ and $\pi_B$ only up to $n=6$, further work could calculate the permutations for longer signals, or prove the existence or non-existence of a general form for larger $n$. Finally, one could examine other signal transformations and their changes under persistence diagrams. Could it be that only affine transformations result in a consistent change in persistence diagrams?

\section{References}

\bibliographystyle{plainnat}
\bibliography{Sources}

\appendix
\section{Lemmas}
\label{sec:app_Lemmas}

\subsection{Negation Lemmas}
\label{app:negation}
\begin{lemma}
\label{lem:Dgm_parts_equivalence}
For a generic function $f : [m^-, m^+] \to \mathbb{R}$, and critical values $A<B$,
\begin{itemize}
\item $(A,B )\in \ess (f) \iff (-B, -A )\in \ess (-f)$
\item $(A,B )\in \ord(f) \iff (-A,-B )\in \rel(-f)$  
\item $(B,A)\in \rel(f) \iff (-B, -A)\in \ord(-f)$
\end{itemize}
\end{lemma}

\begin{proof}
First, as we are on the interval, $\ess (f) = \{(y^-, y^+)\}$, and $\ess (-f) = \{(-y^+, -y^-)\}$.

In the discrete case, the statements for $\ord(f)$ and $\rel(f)$ are precisely Belton's \cite[Theorem 3.14]{Belton2026}, and the proof for the continuous case is essentially the same:  

 Negation reverses the filtration order and swaps minima with maxima, so ordinary persistence of $f$ becomes the relative persistence of $-f$. Substituting $-f$ for $f$ then yields the final equivalence.
\end{proof}

\begin{lemma}
\label{lem:diagram_negation}
    For any generic function $f : [m^-, m^+] \to \rr$,
\[\tvar (f) = \tvar(-f)\]
and
\[\tper (f) = \tper (-f).\] 
\end{lemma}
\begin{proof}
The first equality, 
\[ \tvar (f) =  \int_{m^-}^{m^+} |f'(x)| \dee x = \int_{m^-}^{m^+} |-f'(x)| \dee x = \tvar(-f)\]
is apparent, and
\[ \tper (f) = \sum _{(A,B)\in \dgm(f)}|B-A| = \sum _{(-B, -A)\in \dgm(-f)}|-B - (-A)| = \tper(-f) \]
 follows from the equivalences established in Lemma \ref{lem:Dgm_parts_equivalence}: 

 \[\tper (f) = \|\dgm(f)\|_1 = \|\ord(f)\|_1 + \| \rel(f) \|_1 + \|\ess(f)\|_1 = \]
 \[  = \| \rel(-f) \|_1 + \|\ord(-f)\|_1 + \|\ess(-f)\|_1 = \tper(-f)\]
\end{proof}

\subsection{Other Lemmas}
\begin{lemma}
\label{lem:C1+BV=AC}
    Let $f:\rr \to \rr$ be a piecewise $C^1$ function with bounded variation. Then, $f$ is absolutely continuous.
\end{lemma}
\begin{proof}
Let $\epsilon > 0$.  Because $f$ has bounded variation, we know that its limits
 \[\lim_{x \to -\infty}f(x) := A\]
 and
 \[ \lim_{x \to \infty}f(x) =: B\]
 exist. Furthermore,
\[\int_{-\infty}^\infty|f'(x)| \dee x \leq \tvar(f) < \infty,\]
so $f' \in L^1$ by definition.

Therefore, there exists a $j \in (0,\tvar(f))$ such that 
$j=\int_{-\infty}^\infty|f'(x)|\dee x$, and thus
\[\lim _{i \to \infty} \int_{-i}^i|f'(x)|\dee x = j.\]

Thus, there exists $K\in (0,j)$ such that 
\[ \int_{-\infty}^{-K}|f'(t)|\dee t+\int_K^\infty|f'(t)|\dee t < \frac \epsilon 2.\]

Because $f$ is piecewise $C^1$, the Extreme Value Theorem tells us that $C_1$ is bounded on the compact interval $[-K,K]$. Let us define 
\[\mu = \sup_{t \in [-K,K]}f'(t).\]

Then, we can choose $\delta = \frac \epsilon{2\mu}$. 
Let $R := \{(x_i,y_i)\}_{i=1}^N$ be a collection of disjoint intervals within $\rr$, such that the sum of the intervals' lengths \[\sum_{i=1}^N|x_i-y_i|<\delta.\] 

We will split the set of intervals in $R$ into those which are within $[-K,K]$ and those outside of $[-K,K]$, splitting any intervals which straddle $\pm K$. Thus, we define
\[R_I = \{(x_i,y_i) \cap [-K,K] : (x_i,y_i) \in R, (x_i,y_i) \cap [-K,K]\neq \emptyset\}\]
and
\[R_O = \{(x_i,y_i) \cap \{\rr \setminus[-K,K]\} : (x_i,y_i) \in R, (x_i,y_i) \cap \{\rr \setminus [-K,K]\}\neq \emptyset\},\]
such that $\cup R = \cup(R_I \sqcup R_O)$, and \[R = R_I \sqcup R_O \iff \{\cup R \cap \{-K,K\}\} = \emptyset.\]

First, we consider $R_I$. We consider 
\[\sum_{(x_i,y_i) \in R_I}|f(y_i)-f(x_i)| \leq \sum_{(x_i,y_i) \in R_I} \int_{x_i}^{y_i} |f'(t)| \dee t \leq \sum_{(x_i,y_i) \in R_I}\int_{x_i}^{y_i} \mu \dee t.\]
Because the intervals in $R_I$ are disjoint, we continue
\[\sum_{(x_i,y_i) \in R_I}\int_{x_i}^{y_i} \mu \dee t = \int_{t \in \cup R_I} \mu \dee t = \mu \sum_{(x_i,y_i) \in R_I}|y_i-x_i| < \mu \delta = \frac \epsilon 2.\]

Furthermore, as $K$ is chosen so that 
\[ \int_{-\infty}^{-K}|f'(x)|\dee x+\int_K^\infty|f'(x)|\dee x < \frac \epsilon 2,\]
it follows that integrating over $R_O \subsetneq \{\rr \setminus[-K,K]\}$ that 
\[\int_{t \in K_O} |f'(t)| \dee t \leq\int_{-\infty}^{-K}|f'(t)|\dee t+\int_K^\infty|f'(t)|\dee t < \frac \epsilon 2.\]

Thus, since $\cup R = \cup(R_I \sqcup R_O)$, we know
\[\sum_{(x_i,y_i) \in R}|f(y_i)-f(x_i)| = \sum_{(x_i,y_i) \in R_I}(|f(y_i)-f(x_i)|) + \sum_{(x_i,y_i) \in R_O}(|f(y_i)-f(x_i)|) < \frac \epsilon 2 + \frac \epsilon 2 = \epsilon.\]
\end{proof}

\begin{lemma}
    \label{lem:accum_point_not_dense}
    Let $f:\rr \to \rr$ be a piecewise $C^1$ function, and $I\subseteq \rr$ be an interval (possibly unbounded) within which there is no nontrivial subinterval $I'\subset I$ such that $f$ is constant on $I'$.  
    
    Then the set $C \subset \rr$ of local extrema has a set of accumulation points $A$ for which $A \cap I$ is nowhere dense in $I$. 
\end{lemma}

\begin{proof}
Clearly, if $C$ does not have any accumulation points, the statement is vacuously true. Thus, we assume $A$ is nonempty.

Suppose, for sake of contradiction, that there was an interval $[a,b] \subset I$ over which $f$ is nowhere constant, such that $[a,b] \subset A$. Because the set of points at which $f'$ is discontinuous is finite, we can assume without loss of generality that this interval is fully $C^1$. 

For every point $t \in [a,b]$, there is an infinite sequence of critical points $\{s_k(t)\}_{k=1}^\infty$ such that every $s_k(t)$ is a critical point; that is, $f'(s_k(t)) = 0$. However, because $f'(x)$ is continuous, because $\{s_k(t)\}$ converges to $t$ as $k \to \infty$, then \[\lim_{k \to \infty}f'(s_k(t)) = 0 \implies \lim_{k \to \infty} f'(t) = 0 \implies f'(t) = 0.\]

Thus, $f'(t)=0$ for all $t \in [a,b]$. But then $f$ must be constant on $[a,b]$, because if there were two points $t_1<t_2$ such that $f(t_1) \neq f(t_2)$, then the Mean Value Theorem would imply the existence of a point $t_3 \in (t_1,t_2)$ such that $f'(t_3) = \frac{f(t_2)-f(t_1)}{t_2-t_1} \neq 0$, which contradicts the fact that $f'(t) = 0 \forall t \in [a,b] \supsetneq (t_1,t_2)$. 

But $[a,b] \subset I$, so the fact that $f$ is constant on $[a,b]$ contradicts our assumption. Therefore, no such interval of accumulation points of local extrema can exist. \\

This implies that $A^c$, the complement of $A$ in $\rr$, is dense in $\rr$, because if it were not, that would imply the existence of an $\epsilon$-neighborhood around some point $t$ within which there were no non-accumulation points, which we just showed is impossible.

To show that $A$ is nowhere dense in $I$, we must prove that \[\forall B_\epsilon(x) \subset I \  \exists B_\delta(y) \subseteq B_\epsilon(x) : B_\delta(y) \cap A = \emptyset.\] 

Because $A$ is the set of accumulation points of $C$, it must be a closed set in $\rr$. 

Therefore its complement $A^c$ is an open set in $\rr$, meaning that around every point $b \in A^c$ there exists a $\delta$-ball contained within $A^c$. 

Thus, we consider $B_\epsilon(x) \subset I$. Because no interval contains just points in $A$, there exists some point $y \in A^c \cap B_\epsilon(x)$. As $A^c$ is an open set in $\rr$, $A^c \cap B_\epsilon(x)$ is an open set in $B_\epsilon(x)$, and so there exists some $\delta > 0$ such that $B_\delta(y) \subset A^c$. Finally, because $B_\delta(y) \subset A^c$, it follows that $B_\delta(y) \cap A = \emptyset$, as desired.

\end{proof}

\begin{lemma}
\label{lem:BV_implies_limits}
    Let $f$ be defined on an open interval $(a,b)$. If either  $\lim_{x\to a^+}f(x)$ and $\lim_{x\to b^-}f(x)$ does not exist, then $f$ does not have bounded variation.
\end{lemma}

\begin{proof}
    If one limit (WLOG $\lim_{x\to a^+}f(x)$) did not exist, there exists a monotonic sequence $\{x_i\}\to a^+$ for which each $x_i \in (a,b)$ and $\lim_{x \to a^+} \sup f(x_i) > \lim_{x \to a^+} \inf f(x_i)$. This is because, if there  were no such sequence, then for every monotonic sequence $\{x_i\}$ in $(a,b)$ converging to $a$, $\{f(x_i)\}$ would be convergent to some limit $L$; if two sequences converged to different limits $L$ then a nonconvergent sequence could be made by interleaving them, and if all converged to some limit $L$, then $\lim_{x\to a} f(x)$ would exist and equal $L$.
    
    This means that $\forall \epsilon > 0$, there are infinite subsequences $\{y_j\}$ and $\{z_k\}$ of $\{x_i\}$ such that every $f(y_j) > \limsup f(x_i) - \epsilon$ and every $f(z_k) < \liminf f(x_i)+ \epsilon$.

    Because we $\{y_j\}$ and $\{z_k\}$ are (infinite) subsequences of $\{x_i\}$, we can construct a final sequence $\{w_l\}$ such that $\{w_l\}$ is monotonically decreasing towards $a$, and such that $w_l \in \{y_j\}$ whenever $l$ is odd but $w_l \in \{z_k\}$ whenever $l$ is even. Because $\{w_l\}$ is a subsequence of $\{x_i\}$, we know that 
    \[\sum_{i=1}^\infty |f(x_{i+1})-f(x_i)| \geq \sum_{i=1}^\infty |f(w_{i+1})-f(w_i)| > \sum_{i=1}^\infty\limsup f(x_i) - \liminf f(x_i) - 2\epsilon. \]
    Since, for small enough $\epsilon$, $\limsup f(x_i) - \liminf f(x_i) - 2\epsilon > 0$, we have that the partial sum \[\sum_{i=1}^N |f(x_{i+1})-f(x_i)|\] is unbounded as $N \to \infty$, and therefore that the sum \[\sum_{n=1}^\infty |f(x_{i+1})-f(x_i)|\] is infinite. 

    However, as $f(x)$ passes through every point $f(x)$ in $\{x_i\}$, we know that \[\tvar(f) = \int_a^b |f'(x)| \dee x = \int_{x_1}^b|f'(x)| \dee x + \int_{x_2}^{x_1}|f'(x)| \dee x + \ldots  \geq \sum_{n=1}^\infty |f(x_{i+1})-f(x_i)| ,\] and therefore  that $\tvar(f)$ is unbounded, as desired. 
\end{proof}

\section{Generalizing Equivalence Domain}
\label{sec:app_Generalizing Equivalence Domain}
\subsection{Formula Derivation Proof}
\label{subsec:app_Formula Derivation Proof}

\begin{proposition}
\label{app:tVar=tPer_Proof_1_Elab}
    Let $f$ be a generic continuous function with $n$ critical points. Then, the total variation and total persistence are equal, and can be calculated to be
    \[\tvar(f) = -y_1 + (-1)^ny_n +\sum_{i=2}^{n-1} (-1)^iy_i = \tper (f)\]
\end{proposition}

\begin{proof}
The calculation of total variation proceeds straightforwardly from 

\[\tvar (f) = -y_1 + (-1)^n y_n + 2\sum _{i=2}^{n-1} (-1)^i y_i.\]

Recall Table \ref{tab:tPer_contribution_table} summarizes the contributions of each critical point to the total persistence of $f$ as follows: 
\begin{table}[h!]
    \centering
    \begin{tabular}{c | c c c c c}
        \hline
        & $y_1$ & $y_{even}$ & $y_{odd}$ & $y_n$ & $y^+$\\
        \hline
        Sublevel: & Birth & Death & Birth & ifOdd & +Death\\
        Superlevel: & * & Birth & Death & ifEven & -Birth\\
        \hline Total:  & 1 & 2 & 2 & 1 & *
    \end{tabular}
    \label{table:1_repeat}
\end{table} 
As we vary $f_\alpha$ from $y^-$ to $y^+$, every time we reach a local minimum, a new feature is born. Because $y_1 < y_2$ [by assumption], local minima exist at $y_1$ and $y_i$ for every odd $i$, since local minima and maxima alternate. Therefore, a feature is born in either $\ord (f)$ or $\ess (f)$ (in the case of $y^-$) whenever an odd-indexed critical value $y_{odd}$ is reached. 

Now, let us consider $y_{even}$, the even-indexed local maxima. If a point $y_{even}$ is an \textit{interior point} (that is, not a boundary point $y_n$), then as $\alpha \to y_{even}$, two components in $f_\alpha$ approach one another, and at $f_{y_{even}}$, the younger component dies. Therefore, there is a death in $\ord (f)$ at every interior even-indexed critical value $f_{even}$. 

Furthermore, when the global maximum is reached, an additional death occurs in $\ess (f)$. If the global maximum is at a boundary point, then this will be the only death at the critical value, but if the global maximum is at an interior point, then there will be two deaths which occur at the critical value: the second youngest point will die at this critical value, and a corresponding point will be created in $\ord (f)$, and the youngest point $y^-$ will die, creating a point in $\ess (f)$. \\

The same argument works in the other direction, when considering superlevel persistence and points in $\rel (f)$. A point is born in $\rel (f)$ or $\ess (f)$ (in the case of the global maximum) at every $y_{even}$ \textit{besides the local maximum}; a point dies in $\rel (f)$ at every interior $y_{odd}$. However, because the interval is not itself a loop, there is no point in $\ess (f)$ during the superlevel persistence; this is why there is no point associated with the global maximum in the superlevel persistence. This balances out the ``extra'' death at the global maximum which occurs during the sublevel persistence. \\

From our definition, we have that the total persistence is
\begin{equation*}
    \tper = \sum_{(A,B)\in \dgm (f)} |A-B|.
\end{equation*}

However, we know that every point in $\dgm (f)$ contains exactly one local minimum and one local maximum. Furthermore, the critical value of a local minimum is always less than the critical value of the local maximum with which it is paired, as a point must be born before it dies. Therefore every persistence contribution $|A-B|$ has the form $y_{even} - y_{odd}>0$. 

We know that every interior point is associated with one point in the sublevel persistence, and every \textit{local} interior point (that is, not a global extremum) is associated with one point in the superlevel persistence. We know that boundary point is associated with only one point in $\dgm (f)$: either in sublevel persistence if the boundary point is a minimum or the global maximum, or in superlevel if it is a nonglobal maximum. %

As all local maxima are added and all local minima subtracted, by commutativity we have that

\begin{equation*}
    \tper = \sum_{(A,B)\in \dgm (f)} |A-B| =-y_1 + (-1)^ny_n +\sum_{i=2}^{n-1} (-1)^iy_i.
\end{equation*}
\end{proof}

As Biswas showed, $\tvar (f) = \tper (f)$ for any \textbf{continuous, generic} function. We first wish to show that this identity still holds for continuous non-separating functions. 

Recall that a continuous function $f: \bb I \to \rr$ is generic if are critical points have distinct values (and therefore there are no nonisolated critical points or horizontal segments). One appeal of this restriction is that it guarantees that the elder rule is always unambiguous, because there must always be an elder of two (or more) critical points. 

When we get rid of this restriction, if we have a meeting between two components $k_{a}$ and $k_b$, where $f(a) = f(b) = c$, then there is no way to adjudicate which point survives the meeting. It is equally possible that $k_a$ dies and $k_b$ survives, or that $k_a$ survives and $k_b$ dies. 

However, in a persistence diagram, components are not graphed by their $x$-values, but by their $y$-values. Since $f(a)=f(b)=c$, the persistence diagram will graph two points in the column $c$: one at the level at which $k_a$ and $k_b$ meet, and one at which that merged component dies (or at $y^+$ if it is the essential component). Thus, though the assignment of seniority to $k_a$ and $k_b$ is arbitrary, the persistence diagram produced is in any event identical. 

When more than two critical points share a $y$-level, it can be helpful to think of there being created a (total) ``ordering'' for each shared $y$ level, which determines which component survives a meeting between any two points. 

 Let $\mc O(f)$ denote a strict total ordering of the values $y$ of each critical point of $f$, which is a refinement of the (not strict) ordering imposed on the critical values by the ordering of $\rr$. We will show that the persistence diagram created for $f$ is completely independent of any of these ``persistence orderings'' $\mc O(f)$.
\subsection{Non-Separating Functions}
\label{subsec:app_Non-separating Functions}

\begin{lemma}
\label{lem:two_pts}
    Let $f:\bb I \to \rr$ be a continuous, non-separating function, with two critical points $c_i,c_j$ such that $f(x_i) = f(x_j)$. Suppose the components born at these critical points $k_{c_i}$ and $k_{c_j}$ meet at a level $\ell$. Then the persistence diagram $\dgm (f)$ does not depend on which component survives the meeting and continues past $\ell$. 
\end{lemma}

\begin{proof}
    Since we know that  $k_{c_i}$ and $k_{c_j}$ meet, we know that $c_i$ and $c_j$ must be the same type of extremum.
    Suppose first they are local minima. If $c_i$ is determined to survive past $c_j$, then $\dgm(f)$ will have one point at $(y_j, \ell)$, and one at $(y_i, \ell_2)$, where $\ell_2 \geq \ell$ is the point at which the merged component created from $k_{c_i}$ and $k_{c_j}$ dies.
    
    Conversely, if $k_{c_j}$ is determined to survive past $k_{c_i}$, then $\dgm(f)$ will have points at $(y_i, \ell )$ and $(y_j, \ell_2)$. However, since $y_i=y_j$, these are the same pairs of points. The column at $y=y_i=y_j$ will have points at levels $\ell$ and $\ell_2$; if $\ell=\ell_2$, then these points will be coincident. \\

    The argument repeats in much the same manner if $c_i$ and $c_j$ are instead local maxima. If $y_i=y_j=y^+$, then there will only be a point at $(y, \ell )$, as the number of loops created at $y^+$ will always be one less than the number of points there, so long as the function $f$ is defined on the interval $\bb I$, and not the circle $\bb S^1$. 
\end{proof}

\begin{lemma}
\label{lem:hor_seg}
    If a non-separating continuous function $f$ has a horizontal segment on the interval $[a,b]$ at level $l$, the persistence diagram $\dgm(f)$ does not depend on which point $c$ is chosen to represent a connected component born there. In other words, the persistence diagram does not depend on the persistence ordering of the points in $[a,b]$. 
\end{lemma}

\begin{proof}
    First, suppose the horizontal segment is monotone: that is, neither $a$ nor $b$ is an endpoint of $\bb I$, and $\sign( f'(a-\epsilon)) \neq \sign(f'(b+\epsilon))$ for sufficiently small $\epsilon > 0$. Then, in neither sublevel nor superlevel persistence does any point in the horizontal segment birth a new topological feature. Therefore, there can be no ambiguity in the persistence diagram, and the equality is trivial. \\

    Next, suppose the horizontal segment is a minimum. That is, $\sign(f'(a-\epsilon)) = -1 = -\sign(f'(b+\epsilon))$, if they exist. Then, one point $c \in [a,b]$ in the interval will continue beyond $l$. However, as the component born $k_c$ includes all of $[a,b]$ immediately from the point of its birth, regardless of choice of $c$, the local maximum $\ell_2$ with which it pairs cannot vary with choice of $c$. 

    Therefore, the interval will be represented in $\dgm(f)$ by one point at $(l, \ell_2 )$, regardless of the choice of $c$. 

    Finally, suppose the horizontal segment is a maximum. That is, $\sign(f'(a-\epsilon)) = 1 = -\sign(f'(b+\epsilon))$, if they exist. Then none of these points contribute to $\ord (f)$, or indeed $\ess (f)$, but rather $\rel (f)$. 

    In the case that $l < y^+$, our proof mirrors the minimum case exactly. 

    If $l=y^+$, then exactly one point at this level must be chosen to survive to the end. This is similar to the essential component in sublevel persistence, except that it does not contribute any persistence, because superlevel persistence tracks 1-D loops, not connected components. 

    If our interval $[a,b]$ is the entirety of $f^{-1}(l)$, then the entire interval will contribute no persistence. If the levelset $f^{-1}(l)$ has at least one other point $c'\notin [a,b]$, then the interval $[a,b]$ will contribute no persistence if and only if none of the other points in $f^{-1}(l)$ are chosen to be the ``original component'' which survives until the end. Otherwise, one point in $[a,b]$ will contribute $l-\ell_2$ persistence, where $\ell_2$ is the level at which the component born at the interval merges with the component born at $l=y^+$ which was decided to survive it. 
\end{proof}

\begin{corollarylem}
\label{cor:equiv_rel}
    Let us define an equivalence relation $\sim$ on $\bb I$, the domain of $f$. Let $c_1 \sim c_2$ if $f$ is constant on $[c_1,c_2]$. 

    Let $\tilde f : \bb I / \sim \  \to \rr $ be the function created by collapsing all flat intervals in $f$ to a single point. Then, $\dgm(f) = \dgm (\tilde f)$. 

    Furthermore, $\tvar(f) = \tvar (\tilde f)$ follows from the fact that the variation of a horizontal segment is zero. 
\end{corollarylem}

\begin{lemma}
\label{lem:one_level_arb}
    Let a non-separating continuous function $f$ have finitely many critical points at a level $l$. Suppose $\mc O_1(f)$ and $\mc O_2(f)$ are two persistence orderings of the critical points of $f$, which differ only at $l$. Then the persistence diagrams $\dgm_1(f)$ and $\dgm_2(f)$ generated by the respective orderings are the same. 
\end{lemma}

\begin{proof}
 This follows from Lemma \ref{lem:two_pts} by induction. In particular, any two (finite) orderings are separated by a finite number of swaps, each of which preserve a unique $\dgm(f)$ by Lemma \ref{lem:two_pts}.
\end{proof}

\begin{theorem}
\label{thm:dgm1=dgm2}
   A (not necessarily generic) continuous function $f : \bb I \to \rr$ with finitely many isolated critical points and horizontal segments, has a persistence diagram $\dgm(f)$ which does not depend on the persistence ordering of points at any level. 
\end{theorem}
\begin{proof}
    Let $\mc O_1(f)$ and $\mc O_2(f)$ be two persistence orderings of $f$. We must show that the persistence diagrams $\dgm_1(f)$ and $\dgm_2(f)$ are the same. 
    
    First, we will compute $\tilde f$, as shown in Corollary \ref{cor:equiv_rel}. This gives us that \[\dgm _1(f) = \dgm _1(\tilde f) \text{ and } \dgm _2 (\tilde f) = \dgm _2 (f),\]

    and so it remains to show that 

    \[\dgm _1(\tilde f) = \dgm _2 (\tilde f). \]

    Because all horizontal segments have been collapsed to isolated critical points (or non-critical points), we know that $\tilde f$ has finitely many critical points. It follows that there are finitely many levels $l$ with multiple critical points, each of which has finitely many critical points.
    
    Therefore, a finite number of applications of the procedure in Lemma \ref{lem:one_level_arb} shows that $\dgm_1(f) = \dgm_2(f)$, and recovers the uniqueness of $\dgm(f)$. 
    
\end{proof}

\begin{corollarythm}

    \label{cor:tVar=tPer_arb_cont}
    Let $f: \bb I \to \rr$ be a continuous function with finitely many critical points, which is not necessarily generic. Then the total variation $\tvar (f) = \tper (f)$, the total persistence, which can be computed as
    \[\left |-y_1+(-1)^ny_n +2\sum_{i=2}^{n-1}(-1)^iy_i\right|.\]
\end{corollarythm}

\begin{proof}
First, as in the generic case, we will assume without loss of generality that $y_1 < y_2$. We know that $y_1 \neq y_2$ because there are finitely many points, and therefore they are isolated.

As above, $f(x)$ is monotonic between any two consecutively indexed critical points. Thus, our computation of variation is unchanged, and we find
\[\tvar (f) = -y_1+(-1)^ny_n +2\sum_{i=2}^{n-1}(-1)^iy_i.\]

To compute $\tper(f)$ from $\dgm(f)$, we refer to the earlier Table \ref{tab:tPer_contribution_table}. 
\begin{table}[h!]
    \centering
    \begin{tabular}{c | c c c c c}
        \hline
        & $y_1$ & $y_{even}$ & $y_{odd}$ & $y_n$ & $y^+$\\
        \hline
        Sublevel: & Birth & Death & Birth & ifOdd & +Death\\
        Superlevel: & * & Birth & Death & ifEven & -Birth\\
        \hline Total:  & 1 & 2 & 2 & 1 & *
    \end{tabular}
    \label{tab:1_repeat_2}
\end{table} 
Because the difference is computed identically, regardless of the pairings, we recover the same summation:
\[\tper (f) = \left |-y_1+(-1)^ny_n +2\sum_{i=2}^{n-1}(-1)^iy_i\right|.\]

Thus, $\tvar (f) = \|\dgm(f)\|_1=\tper(f)$. 
\end{proof}

\begin{corollarythm}
\label{cor:dgm1=dgm2_flats}
    Let $f: \bb I \to \rr$ be a function with finitely many critical points and horizontal segments, which is not necessarily generic. Then the total variation $\tvar (f)$ and total persistence $ \tper (f)$ are equal and can be calculated as
    \[\left |-y_1+(-1)^ny_n +2\sum_{i=2}^{n-1}(-1)^iy_i\right|.\]
\end{corollarythm}

\begin{proof}
    Compute $\tilde f$, as shown in Corollary \ref{cor:equiv_rel}. As $\tper(f) = \tper (\tilde f)$ and $\tvar(f) = \tvar (\tilde f)$, we need only compute $\tper(\tilde f) = \tper(\tilde f)$. However, $\tilde f$ has finitely many critical points, and so the equality follows from Theorem \ref{thm:dgm1=dgm2}.
\end{proof}

\subsection{Discrete Functions}
\label{subsec:app_Discrete Functions}
\begin{theorem}
\label{thm:f_and_g}
    Let $g: X \to \rr$ be a \textit{discrete} function, which maps from a finite set of natural numbers $X = \{1,2,\ldots,n\}$ to the real numbers. Let a piecewise differentiable function $f:[1,n] \to \rr$ be the \textit{linear interpolation} of the points of $x$: that is, $\forall x \in [1,n]$,
    \[ f(x) = \begin{cases}
        g(x) & \text{ if } x \in \nn \\
        (\lceil x \rceil - x)g(\lfloor x \rfloor) + (x-\lfloor x \rfloor)g(\lceil x \rceil) & \text{ if } x \notin \nn
    \end{cases}.\]
    Then, $\tvar (f) = \tvar (g)$ and $\tper (f) = \tper (g)$. 
\end{theorem}

\begin{proof}
    By definition, \[\tvar(g) = \sum _{i=1}^{n-1} |g(i+1)-g(i)|.\] As $f$ is monotonic between critical points, of which the set of points in $g$ is a superset, we know that 
    \[\tvar(f) = \int _1^n |f'(x)| \dee x = \sum_{i=1}^{n-1}|f(x_{i+1})-f(x_i)|, \]
where $\{ y_i\}$ is the sequence of \textit{consecutively distinct} critical points: that is, the set of critical points which excludes any additional points in a constant interval, because the total variation within a constant interval is zero.

Furthermore, note that, on any interval over which $f$ is monotonic, for any $x^* \in [x_i,x_{i+1}]$, \[|f(x_{i+1}) - f(x_i)| = |f(x_{i+1})-f(x^*)| + |f(x^*)-f(x_i)|.\]

Therefore, as $f$ is indeed monotonic over every interval between critical points, we can take every point in $X \setminus \{x_i\}$, and add it to the sum 
\[\sum_{i=1}^{n-1}|f(x_{i+1})-f(x_i)|\]
without changing the overall total. Eventually, we will arrive at \[\tvar(g) = \sum _{i=1}^{n-1} |g(i+1)-g(i)|,\]
and as the sum is unchanged, we will have shown that 
\[\tvar (f) = \tvar (g).\] 

To show that $\tper (f) = \tper (g)$, we claim that $\dgm(f)=\dgm(g)$.

To show this, note that the persistence diagram depends only on the locations of local extrema. Interpolation does not create or destroy any isolated extrema, and any discrete horizontal segment just becomes a continuous horizontal segment, in both cases represented by an arbitrary point within the range, which in any event has an identical critical value, and so is displayed identically in the persistence diagram.

Because the discrete function $g$ and its linear interpolation $f$ have identical extrema (modulo the arbitrary choice given by horizontal segments), their persistence diagrams are identical, and their total persistence are thus equal. 

Note that this is not restricted to a linear interpolation $f$ of $g$, and would work equally for any interpolation of $g$ which remained (strictly) monotonic between all critical points. 
\end{proof}

\begin{corollarythm}
    \label{cor:f_and_g}
    Let $g:X\to \rr$ be a finite discrete function. Then $\tvar (g) = \tper(g)$, which can be computed as
    \[\left |-y_1+(-1)^ny_n +2\sum_{i=2}^{n-1}(-1)^iy_i\right|.\]
\end{corollarythm}

\begin{proof}
    Apply Theorem \ref{thm:f_and_g} and then Corollary \ref{cor:dgm1=dgm2_flats} to the resulting interpolation function. 
\end{proof}

\subsection{Functions on the Circle}
\label{subsec:app_Functions on the Circle}
\begin{lemma}
\label{lem:cont_circle}
    Let $f : \bb S^1 \to \rr$ be a continuous function (not necessarily generic) on the circle, with finite isolated critical points. Then, $\tvar (f) = \tper (f)$ are simply
   \[ 2\left|\sum_{i=1}^{n-1} (-1)^iy_i\right|.\]
\end{lemma}

\begin{proof}
   Suppose $f$ has $n-1$ distinct critical points (such that, when written out, $c_1 = c_n$).
   
   Recalling Section \ref{subsec:Circular Domains}, we have that 
   \[\tvar (f) = 2\left |\sum _{i=1}^{n-1} (-1)^i y_i \right |, \]
twice the alternating sum of all critical points.

Furthermore, as there are no boundary points and no horizontal segments, every critical point is a local minimum or local maximum. Thus, every critical point is represented in any representation of $\dgm(f)$ as one birth and one death, with one point born at each of $y^-$ and $y^+$ being arbitrarily chosen to represent the essential class, which is now symmetric. 

Therefore, $\tper(f)$ is simply 
\[ \| \dgm (f) \|_1 = \sum_{( A,B) \in \dgm (f)} |A-B| = 2\left|\sum_{i=1}^{n-1} (-1)^iy_i\right|,\]
   the alternating sum of all critical points. 
\end{proof}

\begin{corollarylem}
\label{cor:cont_circle_flats}
    Let us define an equivalence relation $\sim$ on $\bb S^1$. For a function $f:\bb S^1\to\rr$, let $c_1 \sim c_2$ if $f$ is constant on $[c_1,c_2]$. 

    Let $\tilde f : \bb I / \sim \  \to \rr$ be the function created by collapsing all flat intervals in $f$ to a single point. Then, $\tper (f) = \tper(\tilde f)$ and $\tvar (f) = \tvar (\tilde f)$.
\end{corollarylem}

\begin{proof}
    Because all horizontal segments in $f$ are interior segments, this is a special case proven in Lemma \ref{lem:hor_seg}, and a direct analog to Corollary \ref{cor:dgm1=dgm2_flats}.
\end{proof}

\begin{corollarylem}
\label{cor:cont_circle_combined}
    Let $f:\bb S^1 \to \rr$ be a continuous function (not necessarily generic) on the circle, with finite isolated critical points and horizontal segments. Then, $\tvar (f) = \tper (f)$.
\end{corollarylem}

\begin{proof}
    From $f$, construct $\tilde f$ as in Corollary \ref{cor:cont_circle_flats}. Then apply Lemma \ref{lem:cont_circle}. 
\end{proof}

\begin{theorem}
      Let $g : X \to \rr$ be a discrete function (not necessarily generic) on the circle (that is, $X \subset \bb S^1)$, with finite isolated critical points and horizontal segments. Then, $\tvar (g) = \tper (g)$ are computed as 
      \[ 2\left|\sum_{i=1}^{n-1} (-1)^iy_i\right|. \]
\end{theorem}

\begin{proof}
    Let us define a linear interpolation $f: \bb S^1 \to \rr$ of the points in $g$. The same argument in Theorem \ref{thm:f_and_g} demonstrates that $\tvar(f) = \tvar(g)$ and $\tper(f) = \tper (g)$. Finally, Corollary \ref{cor:cont_circle_combined} of Lemma \ref{lem:cont_circle} gives us that $\tvar(f) = \tper(f)$, completing the proof. 
\end{proof}

\subsection{Non-Tame Functions}
\label{subsec:app_Non-Tame Functions}
\begin{theorem}[Section \ref{subsec:Non-Tame Functions}]
\label{thm:abs_cont_eq}
    If $f$ is an absolutely continuous function defined on a compact interval, then $\tper(f) =\tvar(f)$, which can be computed as a weighted sum of critical values. 
\end{theorem}

\begin{proof}
 As before, we can assume WLOG that all local extrema of $f$ are isolated, by the process in Corollary \ref{cor:equiv_rel}, as $\dgm(f)$ (and thus $\tper(f)$) and $\tvar(f)$ are unaffected by ``shrinking'' the horizontal sections. 

Suppose that $f$ has finitely many local extrema. Then equality of $\tper(f)$ and $\tvar(f)$ is proven in Section \ref{subsec:Formula Derivation} if $f$ is an interval, or else in Section \ref{sec:app_Generalizing Equivalence Domain} if $f$ is defined on a circle. \\

Thus, we suppose next that $f$ has infinitely many local extrema. We first consider the case when $f$ defined on the interval,
\[f:[a,b] \to \rr.\]

First, as $f$ is absolutely continuous, $f' \in L^1$ \cite[Thm. 7.20, p. 148]{Rudin1987}, meaning 
\[ \int_a^b |f'(x)| \dee x < \infty.\]

For any partition of $f$ we have that 
\[|f(x_i) - f(y_i)| \leq \left |\int_{y_i}^{x_i} \left|f'(x)\right| \dee x \right|,\]

and so we know that
\begin{equation} 
\label{eq:bound_variation_inf_crit_points} \sum_{i=1}^N|f(x_i)-f(y_i)| \leq \int_a^b|f'(x)| \dee x< \infty,\end{equation}
and therefore that $f$ has bounded variation. 

Furthermore, the absolute continuity of $f$ means that an arbitrary
\[\int_{k_1}^{k_2} f'(x) \dee x\]
evaluates to
\[f(k_2)-f(k_1).\]

Because the total variation is bounded \eqref{eq:bound_variation_inf_crit_points}, we can take the sum
\[ \tvar(f) = \sup_Q \sum_{i=1}^n |f(x_{i})-f(x_{i-1})| = \sup_Q \sum_{i=1}^{n} a_n,\]

where $Q \subsetneq [a,b]$ is a finite partition of the domain. 

Let us define $Q_\infty$ as a partition of $[a,b]$ which contains an $x_i$ at every local extremum of $f$ within $[a,b]$. 
Let us define $Q_1$ as the partition $Q = \{a<x_i<x_{i+1}<b\}$, with $i$ chosen such that $|f(x_i)-f(x_{i-1})|$ is minimized. For every $j \in \nn$, let us define define the partition \[Q_{j+1} := Q_j \cup\{x_{i-1},x_i\},\]
where \[(x_{i-1}, x_i) = \argmax_{(i-1,i) \notin Q_j} |f(x_i)-f(x_{i-1})|.\]

For every $N \in \bb N$, let $Q_N$ be a refinement of $Q_\infty$ [that is, the set of points in the boundaries of the partition of $Q_N$ is a proper subset of the boundary points of $Q_\infty$]. Then, the total variation over $Q_N$,
\begin{equation}  \label{eq:Q_N}
\sum_{i=1}^N |f(x_i)-f(x_{i-1})|\end{equation}
will be less than that over $Q_\infty$,
\begin{equation} \label{eq:Q_infty} \sum_{i=1}^\infty |f(x_i)-f(x_{i-1})|,\end{equation}
by no more than the sum of the variation over all non-monotonic intervals, which for sufficiently high $N$ includes only intervals with infinitely many critical points. However, because the sum \eqref{eq:Q_infty} is convergent, the sum \eqref{eq:Q_N} will increase arbitrarily close to \eqref{eq:Q_infty}, meaning that $\tvar(f)$, the supremum of the partial sums $Q$, is at least \eqref{eq:Q_infty}.

However, since $f$ is absolutely continuous and monotone on every $(x_{i-1}, x_i)$ in $Q_{\infty}$, by \cite[Thm. 7.20, p. 148]{Rudin1987} we can rewrite each of the absolute differences $|f(x_i)-f(x_{i-1})|$ as 
\[\int_{x_{i-1}}^{x_i}|f'(x)|\dee x\]

\[\sum_{i=1}^{n_Q}|f(x_{Q,i})-f(x_{Q,i-1})| =\sum_{i=1}^{n_Q}|a_{Q,n}| \] and reorder the terms $a_n$ such that 
\[|a_1| \geq |a_2| \geq \ldots \geq |a_n|.\]

Then, for all $\epsilon > 0$, all but finitely many monotonic intervals contribute less than $\epsilon$ to the total variation $\tvar(f)$ on $[a,b]$; that is, we can take a partition $Q$ of intervals for which $|f(x_i)-f(x_{i-1})|\geq \epsilon$, which will always have a finite number of terms, and the sum
\[\tvar_Q(f) = \sum_{i=1}^{n_Q} a_{Q,n}\]
can be made arbitrarily close to 
\[\int_a^b|f'(x)|\dee x\]
by choosing a sufficiently small $\epsilon$. Thus, the supremum 
\[\sup_Q \sum_{i=1}^{n} a_n = \tvar(f)\]
of this set will be precisely
\[\int_a^b|f'(x)|\dee x,\]
just as in the finite case. 

Because $f$ is monotonic between local extrema, we can partition this sum into segments between each critical point $x_i$:
\[\tvar(f) = \sum_{i=1}^\infty |f(x_{i+1})-f(x_i)|,\]
which we know is a convergent sum because of bounded variation. As before, we know that each local maximum is greater than the adjacent local minima, so the local maximal values are added each time they show up, and the local minimal values are subtracted. Recalling the derivations before, we note that the boundary points appear once and all interior points appear twice, and so letting $x_1 = m^- = a$ and $x_\infty = m^+ = b$, we have that the total variation is
\[ \sign[f(a)-f(a+\epsilon)]f(a) + \sign[f(b)-f(b-\epsilon)]f(b) + 2\sum_{x\in IMax(f)}[f(x)] + 2\sum_{x \in IMin(f)}[f(x)].\]
More succinctly, if we define
\[M_f(x) = \begin{cases}
    1 & \text{ if $x$ is a local maximum of $f$} \\
    -1 & \text{ if $x$ is a local minimum of $f$}
\end{cases}\]
and $ICrit$ as the set of interior (nondegenerate) critical points,
then \[\tper(f) = M_f(a)f(a) + M_f(b)f(b) + 2 \sum_{x \in ICrit(f)}M_f(x)f(x),\]
consistent with the formula in \ref{subsec:Formula Derivation}. 

Similarly, even in the infinite case the pairing of critical points continues as shown in Table \ref{tab:tPer_contribution_table}:
\begin{table}[h!]
    \centering
    \begin{tabular}{c | c c c c c}
        \hline
        & $y_1$ & $y_{even}$ & $y_{odd}$ & $y_n$ & $y^+$\\
        \hline
        Sublevel: & Birth & Death & Birth & ifOdd & +Death\\
        Superlevel: & * & Birth & Death & ifEven & -Birth\\
        \hline Total:  & 1 & 2 & 2 & 1 & *
    \end{tabular}
    \label{tab:tPer_contribution_table_3}
\end{table} 

As before, boundary points appear once and interior points appear twice, every absolute value taken is the difference between a local maximum and a local minimum, the local maximum always has the greater value and is therefore added against the subtracted local minimum, and we evaluate the convergent infinite sum to arrive at 
\[\tvar(f) = M_f(a)f(a) + M_f(b)f(b) + 2 \sum_{x \in ICrit(f)}M_f(x)f(x),\]
and so $\tper(f) = \tvar(f)$, as desired. \\

When $f:\bb S^1 \to \rr$ is defined on the circle, there are no boundary points $a$ and $b$, and the same pairing of (interior) critical points as for the interval case yields 

\[\tvar(f) = 2 \sum_{x \in ICrit(f)}M_f(x)f(x) = \tper(f),\]

where the local minimum is treated analogously to the local maximum, just as in Section \ref{subsec:Circular Domains}. 

Therefore, if we let
\[I_f^{Interior}(x) = \begin{cases}
    1 & \text{ if } x \in \text{Int} (f) \\
    0 & \text{ if } x \in \del(f)
\end{cases}\]
we can express the total persistence and total variation of $f$ as
\[\tvar(f) =\sum_{x \in Crit(f)}M_f(x)(1+I^{Interior}_f(x))f(x) = \tper(f).\]
\end{proof}

Though sufficient, absolute continuity is \textit{not} necessary for total persistence to equal total variation; for example, the Cantor function \cite[4.5.3, p. 136]{Athreya2006} is not absolutely continuous, yet has total variation and total persistence of 1.

\subsection{Open and Unbounded Domains}
\label{subsec:app_Open and Unbounded Domains}
\begin{proposition}
\label{prop:cont_extn}
    Let $f$ be a bounded absolutely continuous function defined on the open interval $(a,b)$, such that $\lim_{x\to a} f(x) := A$ and $\lim_{x\to b} f(x) := B$ exist. Then, define $\dgm(f) := \lim_{\epsilon \to 0}\dgm(f|_{[a+\epsilon, b-\epsilon]})$, the limit of the persistence diagrams of the closed intervals contained in the domain of $f$, $(a,b)$.  Then, $\tvar(f) = \tvar(F)$ and $\dgm(f) = \dgm(F)$, where the continuous extension $F$ of $f$ is defined as
    \[F := \begin{cases}
        A & \text{ if } x=a  \\
        f(x) & \text{ if } x\in (a,b) \\
        B & \text{ if } x=b.
    \end{cases}\]
\end{proposition}
\begin{proof}
First, as $f$ is absolutely continuous, \[\tvar(f) =  \int_a^b|f'(x)|\dee x = \tvar(F),\] as they differ only by the integral of two constants, which is zero. As $F$ is an absolutely continuous function on a compact interval, Theorem \ref{thm:abs_cont_eq} characterizes $\tvar(F)$ as
\[\sum_{x \in Crit(f)}M_f(x)(1+I^{Interior}_f(x))f(x),\]
which is therefore $\tvar(f)$ and thus also $\tper(f)$.\\

    It remains to show that $\dgm(f) = \dgm(F)$. For sufficiently small $\epsilon > 0$, let us define $f_\epsilon := f|_{[a+\epsilon, b-\epsilon]}$. By construction, $f_\epsilon$ agrees with $f$ on $[a+\epsilon, b-\epsilon]$.

    Let us define $F_\epsilon$ as an extension of $f_\epsilon$ to $[a,b]$: 
    \[ F_\epsilon := \begin{cases}
      F(a+\epsilon) & \text{ if } x \in [a,a+\epsilon)]  \\
        f_\epsilon(x) & \text{ if } x\in [a+\epsilon,b-\epsilon]\\
        F(b-\epsilon) & \text{ if } x\in (b-\epsilon, b].  
    \end{cases}\]
    
    As $F_\epsilon(a) = A$ and $F_\epsilon(b) = B$, we have that 
    \[\|F-F_\epsilon\|_\infty = \max \{\sup_{[a,a+\epsilon]}|F-F(a+\epsilon)|, \sup_{[b-\epsilon, b]}|F-F(b-\epsilon)|\}\] which goes to zero as $\epsilon \to 0$ by continuity. 
    
    \citet{CohenSteiner2005} gives us the \textit{stability condition} of persistence diagrams. That is, sufficiently close functions will have arbitrarily close persistence diagrams. Symbolically, we get 
    \[d_B(\dgm(f), \dgm(g)) \leq \|f-g\|_\infty,\]
    where the bottleneck distance 
    \[d_B(X,Y) = \inf _\gamma \sup _x \|x - \gamma(x) \|_\infty.\]
    represents the maximum distance between a matched pair of points in the optimal pairing of points in the two persistence diagrams. That this distance goes to zero as $\epsilon$ goes to zero means that the persistence diagrams of $F$ and $F_\epsilon$ converge; that is, $\dgm(F_\epsilon) \to \dgm(F)$ as $\epsilon \to 0$. Therefore, as $\dgm(f)$ was defined as the persistence diagram of $f_\epsilon$ as $\epsilon \to 0$, it remains to show that $\dgm(f_\epsilon) = \dgm(F_\epsilon)$. 

    However, $F_\epsilon$ is an extension of $f_\epsilon$ by a constant line (on each side). As discussed in Corollary \ref{cor:equiv_rel}, a horizontal segment is invisible to persistence, and so $f_\epsilon = \tilde F_\epsilon$, and $\dgm(f_\epsilon) = \dgm(F_\epsilon)$. 
\end{proof}

\begin{proposition}
\label{prop:Unbounded_Domain}
    Let $f(x):\rr \to \rr$ be a piecewise $C^1$ (continuously differentiable) function on $\rr$ with bounded variation.
Then we can compute a sequence of functions
\[f_{M,N} : [M,N] \to \rr,\] such that
\[\tvar(f) = \lim_{M \to -\infty, N \to \infty} \tvar(f_{M,N}) =  \sum_{x \in Crit(f)}M_{f}(x)(1+I^{Interior}_{f}(x))f(x),\] and  \[\dgm(f) = \lim_{M \to -\infty, N \to \infty} \dgm(f_{M,N}).\] 
\end{proposition}

\begin{proof}

First, the absolute continuity of $f$ follows from its bounded variation and the fact that it is piecewise $C^1$; for the proof of this, see Lemma \ref{lem:C1+BV=AC}. 

 For every $M<N \in \rr$, we can construct a piecewise continuous function $f_{M,N} : \rr \to \rr$ such that
 \[f_{M,N}(x)  =\begin{cases}
     A & \text{ if } x \leq M \\
     f(x) & \text{ if } M<x < N \\
     B & \text{ if } x \geq N
 \end{cases}.\]

By an analogous argument to Lemma \ref{lem:hor_seg}, just with a function defined on $\rr$ rather than $\bb I$, we can invoke Corollary \ref{cor:equiv_rel} to create a piecewise continuous function $\tilde f_{M,N} : [M,N] \to \rr$, such that $\dgm(\tilde f_{M,N}) = \dgm(f_{M,N})$. 

Although both total variation and persistence diagrams can be computed on piecewise (absolutely) continuous functions, we can also monotonically perturb $\tilde f_{M,N}$ to make it absolutely continuous.

First, we distinguish between two types of asymptotes of such functions $f$. For each limit (for example, at $\infty$), we call the asymptote of $f$ \textit{(positive) oscillatory} if it crosses the asymptote ($y= B$) infinitely many times; that is, if for any threshold $n \in \rr$, there exists some $x>n$ such that $f(x) = B$. Conversely, we call the asymptote \textit{(positive) non-oscillatory} if there exists $n \in \rr$ such that $\forall x > n, \ f(x) \neq B$. The same characterizations apply analogously for the asymptote at negative $\infty$. Clearly, if $f : \rr \to \rr$, each asymptote is either oscillatory or non-oscillatory. 

In particular, $f$ is either negative oscillatory or negative non-oscillatory. If $f$ is negative oscillatory, then by continuity we can choose a sufficiently small $M$ such that $f(M) = B$ and 
\[\max_{x<M}|f(x) - B| < \delta,\]
for an arbitrarily small $\delta > 0$. 
Conversely, if $f$ is negative non-oscillatory, then we choose a sufficiently small $M$ such that
\begin{enumerate}
    \item $\max_{x < M}|f(x)-M| < \delta,$
    \item $M$ is not an accumulation point of local extrema, and 
    \item $\sign(f'(M)) = \sign(A-f(M))$. 
\end{enumerate}

First, by the definition of a limit there exist infinitely many points that fulfill the first criterion, including all points on the interval $(\infty, \mu)$ for some sufficiently small $\mu$.
Because the set of accumulation points of local extrema is nowhere dense wherever $f$ is non-constant\footnote{See Lemma \ref{lem:accum_point_not_dense}.}, it is easy to find a point which fulfills the second criterion.

Finally, given a point $M_1$ which fulfills the first two criteria, we can find a point $M$ that fulfills the third by applying the Intermediate Value Theorem to $f$ on $[M_1,M_2]$ for some $M_1<M_2$. Because $M < M_2$, the first criterion is also fulfilled for $M$. If $M$ is an accumulation point of local extrema, we can choose instead any point $M'$ within the $\epsilon$-ball $B_\epsilon(M)$ centered around $M$. For small enough $\epsilon$, the first and third criteria are preserved from $M$ (as $f'$ is assumed continuous for sufficiently small $M$), and there will exist many choices of $M'$ which are not accumulation points of local extrema, per Lemma \ref{lem:accum_point_not_dense}.

We choose $N$ with the same (or symmetric) criteria, depending on whether $f$ is \textit{positive} oscillatory or positive non-oscillatory. 

If $M$ (or equivalently, $N$) is oscillatory, we do not need to do anything. If $M$ is non-oscillatory, we choose some $\gamma_- > 0$ sufficiently small that $\forall x \in [M,M+\gamma_-]$, $\sign (f'(x)) = \sign (A-f(M))$. We make a similar choice of $\gamma_+$ if $f$ is positive non-oscillatory, and we define 
\[\hat f_{M,N}(x) = \begin{cases}
    
    \frac 1 {\gamma_-}( (x-M)f(M+\gamma_-) + (M + \gamma_--x)A ) & \text{ if } f \text{ is negative non-oscillatory and }x \in [M,M+\gamma_-] \\
    \frac 1 {\gamma_+}( (x-N-\gamma_+)B + (N-x)f(N-\gamma_+) ) & \text{ if } f \text{ is positive non-oscillatory and }x \in [N-\gamma_+,N] \\
    f_{M,N}(x) &\text{ otherwise }
\end{cases}.\]
 Because $f$ (and thus $f_{M,N}$) is monotonic on $[M,M+\gamma_-]$ and $[N-\gamma_+,N]$ by construction, $\hat f_{M,N}$ is a monotonic perturbation of $f_{M,N}$ which neither creates nor destroys critical points (as it is a linear interpolation), and therefore does not change either $\tvar(f_{M,N})$ [using the supremum formulation if $f_{M,N}$ is not continuous] or $\tper(f_{M,N})$. 

 Finally, because $\hat f_{M,N}$ is constant on $(-\infty,M]$ and again on $[N,\infty)$, we replicate the arguments of Lemma \ref{lem:hor_seg} and Corollary \ref{cor:equiv_rel}, but applied to the setting where the domain of $f$ is $\rr$. 
 We define a function $\tilde f_{M,N} : [M,N] \to \rr$ as in that example, and argue that $\tvar(\tilde f_{M,N}) = \tvar(\hat f_{M,N})$ and $\dgm(\tilde f_{M,N}) = \dgm(\hat f_{M,N})$. Tracing back the chain of equivalences, we get that \[\tvar(f_{M,N}) = \tvar(f_{M,N}) = \tvar(\tilde f_{M,N}) =\sum_{x \in Crit(\tilde f_{M,N})}M_{f_{M,N}}(x)(1+I^{Interior}_{f_{M,N}}(x))f_{M,N}(x) = \tper(f_{M,N}),\] and that $\dgm (f_{M,N})$ can likewise be calculated from $\dgm(\tilde f_{M,N})$, an absolutely continuous function on a compact interval. \\

 It remains to equate the total variation and persistence diagrams of $f_{M,N}$ with those of $f$. 
First, because $f$ has bounded variation, we know that 
 \[ \int_{-\infty}^\infty |f'(x)| \dee x =: j\geq 0,\]
 meaning that 
 \[ \lim _{H \to \infty} \int_{-H}^H |f'(x)| \dee x =j\]
 and thus that 
  \[ \lim _{H \to \infty} \int_{-\infty}^{-H} |f'(x)| \dee x + \int_{H}^\infty |f'(x)| \dee x=0.\]
  Because \[\int_{-\infty}^\infty \hat f_{M,N}(x) \dee x = \int_M^N f(x) \dee x,\]
  we get that \[\int_{M}^N \tilde f_{M,N}(x) \dee x = \int_M^N f(x) \dee x \to \tvar(f)\]
  as $M \to -\infty$ and $N \to \infty$. 

Similarly, by the definition of a limit, we get that
\[\|f_{M,N}-f\|_\infty \to 0\]
as $M \to -\infty$ and $N \to \infty$. Thus, by the stability condition of Cohen-Steiner et al.\ \cite{CohenSteiner2005}, as
\[d_B(\dgm(f_{M,N}),\dgm(f)) \leq \|f_{M,N}-f\|_\infty,\]
the maximum distance between the persistence diagrams goes to zero, by the Squeeze Theorem, as $M \to -\infty$ and $N \to \infty$. 

Because the difference between the total variations of $f$ and $f_{M,N}$, and the difference between the persistence diagrams of $f$ and $f_{M,N}$, vanish as $M$ and $N$ go to their respective infinities, we can write
\[\tvar(f) = \lim_{M \to \infty, N \to -\infty} \tvar(f_{M,N}) = \sum_{x \in Crit(f)}M_{f}(x)(1+I^{Interior}_{f}(x))f(x)\]

and \[\dgm(f) = \lim_{M \to -\infty, N \to \infty} \dgm(f_{M,N}),\]
which is a function of the ordering of all of the (possibly infinite) critical points of $f$ over $\rr$. 

\end{proof}

\begin{corollaryprop}
    \label{cor:half_unbounded_domain}
    Let $f$ be a piecewise $C^1$ function defined on a real half-line (ray); that is, either $(-\infty, \alpha]$ or $[\alpha,\infty)$, for some $\alpha \in \rr$. 

Then we can compute a sequence of functions
\[f_{M} : [-M,\alpha] \to \rr,\] or \[f_{M} : [\alpha,M] \to \rr\]such that
\[\tvar(f) = \lim_{M \to \infty} \tvar(f_{M}) =  \sum_{x \in Crit(f)}M_{f}(x)(1+I^{Interior}_{f}(x))f(x),\] and  \[\dgm(f) = \lim_{M \to -\infty} \dgm(f_{M}).\] 
    
\end{corollaryprop}

\begin{proof}
    We repeat the same proof as for Proposition \ref{prop:Unbounded_Domain}, with focus only on the respective positive or negative asymptote. 
\end{proof}

\begin{theorem}
    If $f$ is an absolutely continuous 1-D signal on a bounded domain, or a piecewise $C^1$ 1-D signal with bounded variation on an unbounded domain, then $\tper = \tvar$, and can be expressed as 
    \[\tvar(f) =\sum_{x \in Crit(f)}M_f(x)(1+I^{Interior}_f(x))f(x) = \tper(f).\]
\end{theorem}
\begin{proof}

We know that the limits $\lim_{x\to a^+}f(x)$ and $\lim_{x\to b^-}f(x)$ both exist, by the contrapositive of Lemma \ref{lem:BV_implies_limits}, and so we can meaningfully calculate $\dgm(f)$ for the open interval $(a,b)$, as per Lemma \ref{prop:cont_extn}, by substituting $F:[a,b] \to \rr$, where the boundary points are defined as the respective one-sided limits. After this substitution, we can assume WLOG that $f$ is not defined on a bounded open interval, but on the closed interval, a circle, or an unbounded domain.

The unbounded real line case is the conclusion of \ref{prop:Unbounded_Domain}, and closed rays are \ref{cor:half_unbounded_domain}. For open rays, we know that the limit point exists on the side of the finite boundary (by Lemma \ref{lem:BV_implies_limits}, where BV is known by assumption), and so the finite open boundary is equated to the corresponding closed ray defined over the open ray's closure in the exact fashion of Proposition \ref{prop:cont_extn}. 

Finally, this equation is known for the compact domain (interval or circle) by Theorem \ref{thm:abs_cont_eq}. 
\end{proof}

\section{Persistence Diagram Computed Examples}
\label{sec:app_Persistence Diagram Computed Examples}
In Section \ref{subsec:Persistence Diagram Computation}, we express the persistence diagram of any discrete function by pairing its critical values together. In particular, we show that, for all discrete functions $h$ with a particular $n$ critical points, there are permutations $\pi_A$ and $\pi_B$ on $n-1$ elements such that $\ord(h)$ and $\rel(h)$ are given by the formulas in Equations \eqref{eq:comp_pers_minA} - \eqref{eq:comp_pers_maxB}.

Recalling the values
\[A:= \{F(x_1), F(x_3), \ldots, F(x_{n\text{ or } n-1})\},\] 
\[B:=\{F(x_2), F(x_4),\ldots,F(x_{n\text{ or } n-1})\},\] 
\[k_- := \lfloor \frac n2 \rfloor -1 = |A|-1,\] and 
\[k_+ := \lceil \frac n2 \rceil -1 = |B|-1;\]
and the sets 
\begin{itemize}
\item $A^i := A \setminus \{F(x_1),F(x_n)\}$,
\item $A^\ell := A \setminus \{y_{(1)}, y_{(n)}\}$,
\item $B^i := B \setminus \{F(x_n)\}$, and 
\item $B^\ell := B \setminus \{y_{(1)}, y_{(n)}\}$
\end{itemize}
as the restrictions of $A$ and $B$ to interior $(i)$ and nonglobal $(\ell)$ critical points, respectively, the following equations compute the persistence diagrams of $h$, depending on permutations $\pi_A: A^i \to A^i$ and $ \pi_B: B^i \to B^i$:

If the first critical point of $h$ is a local minimum, then we have
\begin{equation*}
\ord (h) = \{ ( a^\ell_1, \pi_B \circ h(x_2)), ( a^\ell_2, \pi_B\circ h(x_4)), \ldots, ( a^\ell_{k_+}, \pi_B\circ h(x_{2k_+}))\}\end{equation*} and 
\begin{equation*}
\rel (h) = \{ ( b^\ell_1, \pi_A \circ h(x_3)), ( b^\ell_2, \pi_A\circ h(x_5)), \ldots, ( b^\ell_{k_-}, \pi_A\circ h(x_{2k_-+1}))\}.\end{equation*}

Conversely, if the first critical point of $h$ is a local maximum, then the formulas are reversed, because of the duality of $\ord$ and $\rel$, giving us 
\begin{equation*}
\ord (h) = \{ ( b^\ell_1, \pi_A \circ h(x_3)), ( b^\ell_2, \pi_A\circ h(x_5)), \ldots, ( b^\ell_{k_-}, \pi_A\circ h(x_{2k_-+1}))\}\end{equation*} and 
\begin{equation*} 
\rel (h) = \{ ( a^\ell_1, \pi_B \circ h(x_2)), ( a^\ell_2, \pi_B\circ h(x_4)), \ldots, (a^\ell_{k_+}, \pi_B\circ h(x_{2k_+}))\}.\end{equation*} \\

These formulas hold for all functions $h$ with $n$ critical points, where the first critical point is a local minimum or local maximum respectively. For the following computations, we will combine the maps $\pi_A$ and $\pi_B$ into one function $g$, as defined above. That is, for an image $Y$ of a signal described by a continuous function $h$ (or equivalently, a discrete function $F$), we wish to compute 
\[ g : \psi ^{-1}(Y \setminus \{y_n\}) \to \psi^{-1}(Y \setminus\{ y_{(i)} : i \neq n, \text{ and }H^{-1}(y_{(i)}) = x_1 \text{ or } x_n, \}),\]
or equivalently, letting $G:= \psi \circ g$, 
\[ G : (Y \setminus \{y_n\}) \to (Y \setminus\{ y_{(i)} : i \neq n, \text{ and }H^{-1}(y_{(i)}) = x_1 \text{ or } x_n, \}),\]
such that 
\[G(y_i) = \begin{cases}
    y_{(n)} & \text{ if } i=1, \\
    \pi_A(y_i) & \text{ if } y_i \text{ is a nonglobal maximum,} \\
    \pi_B(y_i) & \text{ if } y_i \text{ is a nonglobal minimum.}
\end{cases}.\]
That is,  
\[ G(y_i) = y_j \iff ( y_i,y_j)  \in \dgm(h).\]
\paragraph{n $\leq$ 4}
\begin{itemize}
    \item  When $n=2$, $G(y_{(1)}) = y_{(2)}$. 
    \item When $n=3$, $G(y_{(1)}) = y_{(3)}$ and $G(y_{(2)}) = F(x_2)$. 
    \item When $n=4$, 
$G(y_{(1)}) = y_4$. If $x_1$ is a local minimum, then $G(y_{(2)}) = F(x_2)$ and $g(y_{(3)}) = F(x_3)$. If $x_1$ is a local maximum, then $g(y_{(2)}) = F(x_3)$ and $g(y_{(3)}) = F(x_2)$. 
\begin{itemize}
    \item Because each set $|A|=|B| = 1$ has one element, both $\pi_A = \pi_B = e$ are the identity map.
\end{itemize}
\end{itemize}

\paragraph{n=5}
If $x_1$ is a local minimum, we have that $\pi_A$ has only one element, and is therefore trivial. However, $\pi_B$ is not, and can take one of two values, the identity $e :=\begin{pmatrix}
    y_{(i)} & y_{(j)} \\
    y_{(i)} & y_{(j)}
\end{pmatrix}$ or the nonidentity $d := \begin{pmatrix}
    y_{(i)} & y_{(j)} \\
    y_{(j)} & y_{(i)}
\end{pmatrix}$. \footnote{This is a permutation, not a matrix. Thus $e$ maps $y_{(i)}$ and $y_{(j)}$ each to themselves, and $d$ maps them to each other. } In particular, we can show
\[ \pi_B = \begin{cases}
    d & \text{ if } F(x_5)=y_{(3)} \text{ or } F(x_4) = y_{(3)},\\
    e & \text{ if } F(x_1) = y_{(3)} \text{ or } F(x_2) = y_{(3)}, \\
   d & \text{ if } F(x_3) = y_{(3)} \text{ and } F(x_4) = y_{(5)}, \\
   e & \text{ if } F(x_3) = y_{(3)} \text{ and } F(x_2) = y_{(5)}.
\end{cases}\]

The latter two conditions reflect that, when $F(x_3) = y_{(3)}$, one will always have that $g(y_{(2)}) = y_{(5)}$ and $g(y_{(3)}) = y_{(4)}$. The latter part of the former conditions reflect that when $y_{(3)}$ is a local maximum, $g(y_{(2)}) = y_{(3)}$ and $g(y_{(4)}) = y_{(5)}$. \\

If $x_1$ is a local maximum, then 
\[ \pi_B = \begin{cases}
    d & \text{ if } F(x_5)=y_{(3)} \text{ or } F(x_2) = y_{(3)},\\
    e & \text{ if } F(x_1) = y_{(3)} \text{ or } F(x_4) = y_{(3)}, \\
   d & \text{ if } F(x_3) = y_{(3)} \text{ and } F(x_4) = y_{(1)}, \\
   e & \text{ if } F(x_3) = y_{(3)} \text{ and } F(x_2) = y_{(1)}.
\end{cases}\]

\paragraph{n=6}

Suppose first that $x_1$ is a local minimum. Then, $A$ is the set of local minima, and $B$ is the set of local maxima. Then, 
\[ \pi_A = \begin{cases}
    e & \text{ if }F(x_2) = y_{(2)}, \\
    d & \text{ if } F(x_6) = y_{(2)}, \\
    e & \text{ if } F(x_4) = y_{(2)} \text{ and } F(x_3) < F(x_5), \\
    d & \text{ if } F(x_4) = y_{(2)} \text{ and } F(x_3) > F(x_5) \\
\end{cases},\]
and 
\[ \pi_B = \begin{cases}
    e & \text{ if }F(x_5) = y_{(5)}, \\
    d & \text{ if } F(x_1) = y_{(5)}, \\
    e & \text{ if } F(x_3) = y_{(5)} \text{ and } F(x_4) > F(x_2), \\
    d & \text{ if } F(x_3) = y_{(5)} \text{ and } F(x_4) < F(x_2) \\
\end{cases}.\]

When $x_1$ is a local maximum, we can simply reindex the $x$ coordinates by $\tilde F(x_i) = F(x_{7-i})$, creating a new function $\tilde F$ with an identical persistence diagram, but a first local minimum. The new $\pi_A$ and $\pi_B$ are computed accordingly. 

Flowcharts summarizing the calculation of persistence diagrams in the $n=5$ and $n=6$ cases are shown in Appendix \ref{subsec:app_Computed PD Flowcharts}.

\paragraph{n $\geq$ 7} When $n \geq 7$, at least one of the $\pi$s will have at least $3! = 6$ possible values, so the complexity grows on the order of $O( (\frac n2 -1)!)$, and quickly becomes too large to concisely enumerate.

\section{Diagonal Affine Transformation Cases}
\label{sec:app_Diagonal Affine Transformation Theorem Cases}

This section gives the full proof of the Diagonal Affine Transformation Theorem, from Section \ref{subsec:Diagonal Affine Transformation Theorem}. 

\begin{theorem}[Diagonal Affine Transformation Theorem]
\label{thm:app_sgnl_trans_thm} The effect of affine transformations upon total persistence, total variation, and the persistence diagram of a function $f$ on the interval can be summarized as follows: 
\begin{enumerate}
    \item  Total variation, total persistence, and persistence diagrams are invariant under horizontal shifts, horizontal scaling, and horizontal reflection. 

    \item Total variation and total persistence are invariant under vertical shifts, while the persistence diagram will be shifted diagonally.

    \item Vertical scaling will scale the persistence diagram, the total variation, and the total persistence, by the (positive) scaling constant.
    
    \item Total persistence and total variation are invariant under vertical reflection. Furthermore, each point in Ord and Rel is reflected through the origin, while the point in Ess is reflected across the anti-diagonal line. 
\end{enumerate}
\end{theorem}

    Let $f : X \to \rr$ be a continuous function with finitely many isolated critical points, on the set $X = [a,b] \subset \rr$. Remembering that equality of $\dgm$ implies equality of $\tper$, we will address each of these cases in turn: 

\begin{propositionthm}
\label{prop:hor_shift}
    Total variation, total persistence, and persistence diagrams are invariant under a horizontal shift.
\end{propositionthm}

\begin{proof}
    Let us define $h(x) := f(x-c)$, for some constant $c \in \rr$. $f'(x)=h'(x+c)$, so
    \[ \tvar (f) = \int_a^b|f'(x)| \dee x = \int_{a+c}^{b+c}|h'(x)| \dee x = \tvar(h). \]
    From the definition of $h$, we know that the $i$th critical point of $h$, $X_i= x_i + c$. Thus, since the ordering and value of critical points are unchanged, we have that $f(x_i) = h(X_i)$, which means
    \[ \forall i,j \in I_n, \ f(x_i)<f(x_j) \iff h(X_i) < h(X_j),\] and therefore that $\dgm(f) = \dgm(h)$. 
    
\end{proof}
\begin{propositionthm}
\label{prop:hor_scaling}
    Total variation, total persistence, and the persistence diagram are invariant under positive horizontal scaling. 
\end{propositionthm}
\begin{proof}
    Let $h(x) := f(\frac xc)$, for some $c > 0$. Then $|h'(x)| = \frac 1c|f'(\frac xc)|$, and so we can take
    \[ \tvar (h) = \int_{ca}^{cb}|h'(x)| \dee x = \int_{ca}^{cb}\frac 1c\left|f'\left(\frac xc\right)\right| \dee x.\]
    Next, let us substitute $u = \frac xc$, meaning $c\dee u = $ d$x$. This gives us 

    \[\int_{a}^{b}\frac 1c\left|f'\left(u\right)\right| c\dee u =\int_{a}^{b}|f'(u)| \dee u = \tvar(f).\]

    Because horizontal scaling does not change the order, we know that the $i$th critical point of $h$ is \[X_i = cx_i,\] the scaled $i$th critical point of $f$. As $h(X_i) = f(x_i)$, we know that
    \[ \forall i,j \in I_n, \ f(x_i)<f(x_j) \iff h(X_i) < h(X_j),\] and thus, by Theorem \ref{thm:PermThm}, that $\dgm(f) = \dgm(h)$.
\end{proof}

\begin{corollarypropthm}
    \label{cor:horiz_scal_local}

    Total variation, total persistence, and the persistence diagram are invariant under \textbf{partial} positive horizontal scaling.
\end{corollarypropthm}

\begin{proof}
    A \textit{partial} horizontal scaling is defined as a horizontal scaling that affects only a sub-interval of the domain. 

    Given constants $a\leq b<c\leq d \in \rr$, consider functions $f : [a,d] \to \rr$. Because $a=b, c=d$ is precisely the case of Proposition \ref{prop:hor_scaling}, we assume WLOG that $c<d$; if this is not the case, we can reflect our signal by Proposition \ref{prop:hor_reflection}.
    
    We will scale the interval $[b,c]$ by a factor $k>0$. Because both quantities and the persistence diagram are invariant under a horizontal shift, per Proposition \ref{prop:hor_shift}, we can assume without loss of generality that the left endpoint $a$ is held constant. 

    Thus, we define $g:[a,d+\frac{k-1}{c-b}]$ as
    \[g(x):= \begin{cases}
        f(x) & \text{ if } a \leq x <b, \\
        f(b+k(x-b)) &\text{ if }b \leq x \leq c,\\
        f(x+(k-1)(c-b)) & \text{ if } c<x\leq d
    \end{cases}.\]
    By the chain rule, we compute
     \[g'(x):= \begin{cases}
        f'(x) & \text{ if } a \leq x <b, \\
        kf'(b+k(x-b)) &\text{ if }b \leq x \leq c,\\
        f'(x+(k-1)(c-b)) & \text{ if } c<x\leq d
    \end{cases},\]
and so
    \[\tvar(g) = \int_a^d |g'(x)| \dee x = \int_a^b|g'(x)| \dee x + \int_b^c|g'(x)| \dee x + \int_c^d|g'(x)| \dee x =\]
    \[= \int_a^b|f'(x)| \dee x + \int_b^{b+k(c-b)}k\left|f'\left(\frac xk\right)\right| \dee x + \int_{b+k(c-b)}^{d+(k-1)(c-b)}|f'(x)| \dee x.  \]
    Substituting $u = \frac xk$, and therefore $\dee u = \frac{\dee x}k$, we get
    \[= \int_a^b|f'(x)| \dee x + \int_b^{b+k(c-b)}k\left|f'\left(u\right)\right| \dee u + \int_{b+k(c-b)}^{d+(k-1)(c-b)}|f'(x)| \dee x = \tvar(f).\]

Because $k>0$, the ordering of critical points is unchanged. As the critical values are also unchanged, the sequence of critical values of $f$ and $g$ are the same, and the equivalence of $\dgm(f)$ and $\dgm(g)$, and thus also of $\tper(f) = \tper(g)$, follow from Theorem \ref{thm:PermThm}. 
\end{proof}
\begin{propositionthm}
\label{prop:hor_reflection}
    Total variation, total persistence, and the persistence diagram are invariant under horizontal reflection. 
\end{propositionthm}

\begin{proof}
    Because any horizontal reflection can be decomposed into a horizontal reflection about the $y$-axis and a horizontal translation, we will suppose, without loss of generality, that $h(x) := f(-x)$.

Using the substitution $u = -x$, we evaluate
\[ \tvar (f) = \int_a^b|f'(x)| \dee x = \int_{-b}^{-a}|f'(-u)| \dee u = \int_{-b}^{-a}|h'(u)| \dee u =\tvar(h). \]

Furthermore, as the order of critical points is reversed, we have that  
\begin{equation*}
    \forall i,j \in I_n, \ f(x_i)<f(x_j) \iff F(X_{n+1-i})<F(X_{n+1-j}),
    \end{equation*} 
so Corollary \ref{cor:PermThm1} gives us that
\[(y_i, y_j) \in \dgm(f) \iff (Y_i, Y_j) \in \dgm(h);\]
the equality of the persistence diagrams follows then from the fact that $f$ and $h$ have the same set of critical values. 
\end{proof}

\begin{propositionthm}
\label{prop:vert_shift}
    Total variation and total persistence are invariant under a vertical shift, and the persistence diagram is translated diagonally. 
\end{propositionthm}

\begin{proof}
    Let us define $h(x) := f(x) + c$, for some constant $c \in \rr$. First, we can directly evaluate total variation to get that
    \[ \tvar (f) = \int_a^b|f'(x)| \dee x = \int_a^b|h'(x)| \dee x = \tvar(h),\]
    as $f' = h'$. 
    Furthermore, as $c$ is a constant added to $f$, neither the position nor the relative ordering of critical points are changed. Thus, $x_i = X_i$, and thus $h(X_i) = f(x_i) + c$, and
     \[ \forall i,j \in I_n, \ f(x_i)<f(x_j) \iff f(x_i) + c <f(x_j) + c\]
     is trivial. By Theorem \ref{thm:PermThm}, 
     \[(y_i,y_j) \in \dgm(f) \iff (Y_i, Y_j) \in \dgm(h),\] meaning
     \[(f(x_i), f(x_j)) \in \dgm(f) \iff (f(x_i) + c, f(x_j)+c)) \in \dgm(h),\] and so $\dgm(h)$ is a diagonal translation of $\dgm(f)$ by $c$ in each direction, as both birth and death times are translated by $c$. As $\dgm$ is merely translated, 
     \[\tper(h) = \sum_{(Y_i, Y_j) \in \dgm(h)}|Y_i - Y_j| = \sum_{(y_i+c, y_j+c) \in \dgm(h)}|y_i +c- (y_j+c)| = \sum_{(y_i, y_j) \in \dgm(f)}|y_i - y_j| = \tper(f).\]
\end{proof}

\begin{propositionthm}
\label{prop:vert_scale}
    Vertical scaling will scale the persistence diagram, the total variation, and the total persistence by the (positive) scaling constant. 
\end{propositionthm}

\begin{proof}
    Let $h(x):=cf(x)$, for some constant $c>0$. As $|h'(x)| = c|f'(x)|$, we find
    \[\tvar(h) = \int_a^b |h'(x)| \dee x = c\int_a^b |f'(x)| \dee x = c* \tvar(f).\]

    As $c>0$, the ordering is preserved, and so the implication \eqref{eqn:PermThmSuppose} holds. Therefore, Theorem \ref{thm:PermThm} provides that 
    \[(y_i, y_j) \in \dgm(f) \iff (Y_i, Y_j) \in \dgm(h).\] By the definition of $h$, we get that each $Y_i = cy_i$, and so 
    \[(y_i, y_j) \in \dgm(f) \iff (cy_i, cy_j) \in \dgm(h), \]
    meaning $\dgm(h)$ is equal to $\dgm(f)$ scaled by $c$. 
    Finally, 
    \[\tper(h) = \sum_{(Y_i, Y_j) \in \dgm(h)}|Y_i - Y_j| = \sum_{(cy_i, cy_j) \in \dgm(h)}|cy_i - cy_j| = \sum_{(y_i, y_j) \in \dgm(f)}c|y_i - y_j| = c* \tper(f).\]
\end{proof}

\begin{propositionthm}
\label{prop:vert_reflect}
    Total persistence and total variation are invariant under vertical reflection. Points in the ordinary and relative persistence diagrams are reflected through the origin, while the essential point $(y_{(1)},y_{(n)})$ is mapped instead to $(-y_{(n)},-y_{(1)})$. 
\end{propositionthm}

\begin{proof}
    Because any vertical reflection can be decomposed into a vertical reflection about the $x$-axis and a vertical translation, we will consider $h(x) = -f(x)$. Because $|h'(x)| = |f'(x)|$, 
    \[\tvar(h) = \tvar(f)\] is apparent. The persistence diagram changes according to Lemma \ref{lem:Dgm_parts_equivalence} in Appendix 
    \ref{app:negation}. Thus, 
    \[\tper(h) = \sum_{(Y_i, Y_j) \in \dgm(h)}|Y_i - Y_j| = |Y_{(n)}-Y_{(1)}| + \sum_{(Y_i, Y_j) \in \ord(h)\sqcup \rel(h)}|Y_i - Y_j| = \]
    \[|Y_{(n)}-Y_{(1)}| + \sum_{(-Y_i, -Y_j) \in \ord(f)\sqcup \rel(f)}|-Y_j - (-Y_i)|  = |y_{(n)}-y_{(1)}| + \sum_{(y_i,y_j) \in \ord(f)\sqcup \rel(f)}|y_i-y_j| = \tper(f).\]
\end{proof}

\begin{propositionthm}
\label{prop:circle_affine}
    Let $F : \mathbb S^1 \to \rr$ have finitely many isolated critical points. Next, let $f:\bb I \to \rr$ be defined by
    \[ f(x) = F(e^{2\pi ix}) \ \forall x \in [0,1]. \]
    Then, $\dgm (f) = \dgm (F)$ and $\tvar(f) = \tvar(F)$. 
\end{propositionthm}

\begin{proof}
    First, we undergo a horizontal scaling by $2\pi$, for which we repeat the proof for Proposition \ref{prop:hor_scaling}. Then, because $[0,2\pi] / 2\pi \bb Z$ is isomorphic to $\bb S^1$, the persistence diagram and total variation will be the same between identical functions on those domains, as no critical points are created or destroyed, and no critical values are changed.
\end{proof}

\begin{corollarypropthm}
\label{cor:circle_affine}
    If $H : \bb S^1 \to \rr$ is a rotation, horizontal or vertical reflection, vertical shift, or vertical scaling of $F$, then the total variation, total persistence, and persistence diagrams of $F$ and $H$ will vary in accordance with the analogous transformation in Theorem \ref{thm:app_sgnl_trans_thm}.
\end{corollarypropthm}

\begin{proof}
    We define $h$ analogously to $f$, then compare $f$ and $h$ via Theorem \ref{thm:app_sgnl_trans_thm}. 
\end{proof}

\begin{propositionthm}
    \label{prop:discrete_affine}
Let $g$ be a discrete function, and $G$ be the (discrete) function $g$ under an affine transformation. Then, the persistence diagrams of $g$ and $G$ will be related as per Theorem \ref{thm:app_sgnl_trans_thm}. 
\end{propositionthm}

\begin{proof}
    In Section \ref{subsec:Discrete Functions}, we show that the persistence diagram of any function is the same as that of its linear (or otherwise monotonic-respecting) interpolation. Let $f$ and $F$ be such interpolations of $g$ and $G$, respectively. Then, Theorem \ref{thm:app_sgnl_trans_thm} relates the persistence diagrams of $f$ and $F$, and thus also of $g$ and $G$.
\end{proof}
\section{Direct Computation}
\label{sec:app_Direct Computation}
\begin{figure}
    \centering
    \includegraphics[width=\linewidth]{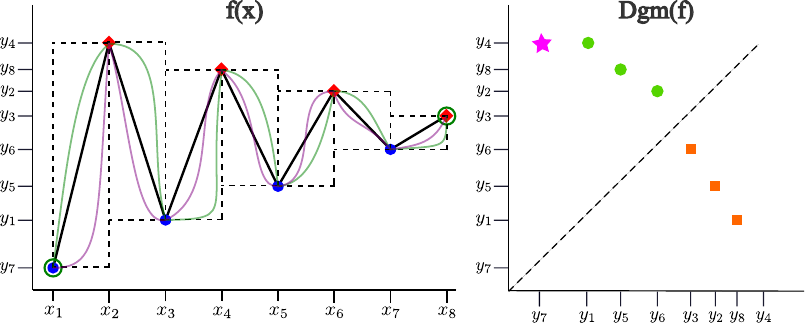}
    \caption{The functions in purple and green share the same sequence of local extrema, and so their persistence diagrams are identical.}
    \label{fig:Monotonic_Perturbation}
\end{figure}

As discussed in the introduction, qualitative changes to the sublevel and superlevel sets only occur at critical levels, specifically the levels of local extrema. Therefore, two functions with the same sequence of extrema will have the same persistence diagram, and therefore the persistence diagram is invariant to changes between critical values which respect monotonicity (so as to not create or change any extrema), as seen in Figure \ref{fig:Monotonic_Perturbation}.

A consequence of this is that the persistence diagram is completely determined by the sequence of local extrema, and so a logical endeavor is to try to find a closed formula for the persistence diagram of a function with a certain sequence of critical values. 

\subsection{PD Preimage Illustrations}
\label{subsec:app_PD Preimage Illustrations}
 \vspace{10pt}

We show example preimages, with the connecting persistence transpositions shown. A dotted line and $\tau_i^*$ denotes a swap which is not a persistence transposition. The black and blue colors denote functions which are reversals of one another. 
 
The preimage of $\dgm = \{(1,8), (2,8),(3,6),(4,5),(5,4),(6,3),(7,2),(8,1)\}$, with connections shown:

\begin{tikzpicture}
\tikzstyle{funcnode}=[font=\small]
\node[funcnode] (n1t)  at (0,4.2)   {18263547};
\node[funcnode] (n1b)  at (0,3.8)   {\textcolor{blue}{74536281}};
\node[funcnode] (n2t)  at (3,4.2)   {18264537};
\node[funcnode] (n2b)  at (3,3.8)   {\textcolor{blue}{73546281}};
\node[funcnode] (n3t)  at (6,4.2)   {18254637};
\node[funcnode] (n3b)  at (6,3.8)   {\textcolor{blue}{73645281}};
\node[funcnode] (n4t)  at (9,4.2)   {18452637};
\node[funcnode] (n4b)  at (9,3.8)   {\textcolor{blue}{73625481}};
\node[funcnode] (n5t)  at (12,4.2)  {15482637};
\node[funcnode] (n5b)  at (12,3.8)  {\textcolor{blue}{73628451}};
\node[funcnode] (n6t)  at (12,2.2)  {\textcolor{blue}{72638451}};
\node[funcnode] (n6b)  at (12,1.8)  {15483627};
\node[funcnode] (n7t)  at (9,2.2)   {\textcolor{blue}{72635481}};
\node[funcnode] (n7b)  at (9,1.8)   {18453627};
\node[funcnode] (n8t)  at (6,2.2)   {\textcolor{blue}{72645381}};
\node[funcnode] (n8b)  at (6,1.8)   {18354627};
\node[funcnode] (n9t)  at (3,2.2)   {\textcolor{blue}{72546381}};
\node[funcnode] (n9b)  at (3,1.8)   {18364527};
\node[funcnode] (n10t) at (0,2.2)   {\textcolor{blue}{74526381}};
\node[funcnode] (n10b) at (0,1.8)   {18362547};
\node[funcnode] (n11t) at (0,0.2)   {16382547};
\node[funcnode] (n11b) at (0,-0.2)  {\textcolor{blue}{74528361}};
\node[funcnode] (n12t) at (3,0.2)   {16384527};
\node[funcnode] (n12b) at (3,-0.2)  {\textcolor{blue}{72548361}};
\node[funcnode] (n13t) at (6,0.2)   {16354827};
\node[funcnode] (n13b) at (6,-0.2)  {\textcolor{blue}{72845361}};
\node[funcnode] (n14t) at (9,0.2)   {16453827};
\node[funcnode] (n14b) at (9,-0.2)  {\textcolor{blue}{72835461}};
\node[funcnode] (n15t) at (12,0.2)  {15463827};
\node[funcnode] (n15b) at (12,-0.2) {\textcolor{blue}{72836451}};
\draw[<->] (0.7,4.0) -- node[midway, above] {$\tau_5$} node[midway, below] {\textcolor{blue}{$\tau_2$}} (2.3,4.0);
\draw[<->] (3.7,4.0) -- node[midway, above] {$\tau_4$} node[midway, below] {\textcolor{blue}{$\tau_3$}} (5.3,4.0);
\draw[<->] (6.7,4.0) -- node[midway, above] {$\tau_3$} node[midway, below] {\textcolor{blue}{$\tau_4$}} (8.3,4.0);
\draw[<->] (9.7,4.0) -- node[midway, above] {$\tau_2$} node[midway, below] {\textcolor{blue}{$\tau_5$}} (11.3,4.0);
\draw[<->] (n5b) -- node[midway, right] {$\tau_5$} node[midway, left] {\textcolor{blue}{$\tau_2$}} (n6t);
\draw[<->] (11.3,2.0) -- node[midway, above] {\textcolor{blue}{$\tau_5$}} node[midway, below] {$\tau_2$} (9.7,2.0);
\draw[<->] (8.3,2.0) -- node[midway, above] {\textcolor{blue}{$\tau_4$}} node[midway, below] {$\tau_3$} (6.7,2.0);
\draw[<->] (5.3,2.0) -- node[midway, above] {\textcolor{blue}{$\tau_3$}} node[midway, below] {$\tau_4$} (3.7,2.0);
\draw[<->] (2.3,2.0) -- node[midway, above] {\textcolor{blue}{$\tau_2$}} node[midway, below] {$\tau_5$} (0.7,2.0);
\draw[<->] (n10b) -- node[midway, right] {$\tau_2$} node[midway, left] {\textcolor{blue}{$\tau_5$}} (n11t);
\draw[<->] (0.7,0.0) -- node[midway, above] {$\tau_5$} node[midway, below] {\textcolor{blue}{$\tau_2$}} (2.3,0.0);
\draw[<->] (3.7,0.0) -- node[midway, above] {$\tau_4$} node[midway, below] {\textcolor{blue}{$\tau_3$}} (5.3,0.0);
\draw[<->] (6.7,0.0) -- node[midway, above] {$\tau_3$} node[midway, below] {\textcolor{blue}{$\tau_4$}} (8.3,0.0);
\draw[<->] (9.7,0.0) -- node[midway, above] {$\tau_2$} node[midway, below] {\textcolor{blue}{$\tau_5$}} (11.3,0.0);
\draw[dotted, thick, <->] (n4b) -- node[midway, right] {$\tau_5^*$} node[midway, left] {\textcolor{blue}{$\tau_3^*$}} (n7t);
\draw[dotted, thick, <->] (n9b) -- node[midway, right] {$\tau_2^*$} node[midway, left] {\textcolor{blue}{$\tau_4^*$}} (n12t);
\draw[dotted, thick, <->] (n1b) -- node[midway, left] {$\tau_4^*$} node[midway, right] {\textcolor{blue}{$\tau_3^*$}} (n10t);
\draw[dotted, thick, <->] (n6b) -- node[midway, left] {$\tau_2^*$} node[midway, right] {\textcolor{blue}{$\tau_3^*$}} (n15t);
\end{tikzpicture}

\vspace{15pt}
The preimage of $\dgm = \{(1,8), (2,7),(3,6),(4,5),(5,4),(6,3),(7,2),(8,1)\}$, with connections shown:

\begin{tikzpicture}
\tikzstyle{funcnode}=[font=\small]
\node[funcnode] (n1t)  at (0,4.2)   {17263548};
\node[funcnode] (n1b)  at (0,3.8)   {\textcolor{blue}{84536271}};
\node[funcnode] (n2t)  at (3,4.2)   {17264538};
\node[funcnode] (n2b)  at (3,3.8)   {\textcolor{blue}{83546271}};
\node[funcnode] (n3t)  at (6,4.2)   {17254638};
\node[funcnode] (n3b)  at (6,3.8)   {\textcolor{blue}{83645271}};
\node[funcnode] (n4t)  at (9,4.2)   {17452638};
\node[funcnode] (n4b)  at (9,3.8)   {\textcolor{blue}{83625471}};
\node[funcnode] (n5t)  at (12,4.2)  {15472638};
\node[funcnode] (n5b)  at (12,3.8)  {\textcolor{blue}{83627451}};
\node[funcnode] (n6t)  at (12,2.2)  {\textcolor{blue}{82637451}};
\node[funcnode] (n6b)  at (12,1.8)  {15473628};
\node[funcnode] (n7t)  at (9,2.2)   {\textcolor{blue}{82635471}};
\node[funcnode] (n7b)  at (9,1.8)   {17453628};
\node[funcnode] (n8t)  at (6,2.2)   {\textcolor{blue}{82645371}};
\node[funcnode] (n8b)  at (6,1.8)   {17354628};
\node[funcnode] (n9t)  at (3,2.2)   {\textcolor{blue}{82546371}};
\node[funcnode] (n9b)  at (3,1.8)   {17364528};
\node[funcnode] (n10t) at (0,2.2)   {\textcolor{blue}{84526371}};
\node[funcnode] (n10b) at (0,1.8)   {17362548};
\node[funcnode] (n11t) at (0,0.2)   {16372548};
\node[funcnode] (n11b) at (0,-0.2)  {\textcolor{blue}{84527361}};
\node[funcnode] (n12t) at (3,0.2)   {16374528};
\node[funcnode] (n12b) at (3,-0.2)  {\textcolor{blue}{82547361}};
\node[funcnode] (n13t) at (6,0.2)   {16354728};
\node[funcnode] (n13b) at (6,-0.2)  {\textcolor{blue}{82745361}};
\node[funcnode] (n14t) at (9,0.2)   {16453728};
\node[funcnode] (n14b) at (9,-0.2)  {\textcolor{blue}{82734561}};
\node[funcnode] (n15t) at (12,0.2)  {15463728};
\node[funcnode] (n15b) at (12,-0.2) {\textcolor{blue}{82736451}};
\draw[<->] (0.7,4.0) -- node[midway, above] {$\tau_5$} node[midway, below] {\textcolor{blue}{$\tau_2$}} (2.3,4.0);
\draw[<->] (3.7,4.0) -- node[midway, above] {$\tau_4$} node[midway, below] {\textcolor{blue}{$\tau_3$}} (5.3,4.0);
\draw[<->] (6.7,4.0) -- node[midway, above] {$\tau_3$} node[midway, below] {\textcolor{blue}{$\tau_4$}} (8.3,4.0);
\draw[<->] (9.7,4.0) -- node[midway, above] {$\tau_2$} node[midway, below] {\textcolor{blue}{$\tau_5$}} (11.3,4.0);
\draw[<->] (n5b) -- node[midway, right] {$\tau_5$} node[midway, left] {\textcolor{blue}{$\tau_2$}} (n6t);
\draw[<->] (11.3,2.0) -- node[midway, above] {\textcolor{blue}{$\tau_5$}} node[midway, below] {$\tau_2$} (9.7,2.0);
\draw[<->] (8.3,2.0) -- node[midway, above] {\textcolor{blue}{$\tau_4$}} node[midway, below] {$\tau_3$} (6.7,2.0);
\draw[<->] (5.3,2.0) -- node[midway, above] {\textcolor{blue}{$\tau_3$}} node[midway, below] {$\tau_4$} (3.7,2.0);
\draw[<->] (2.3,2.0) -- node[midway, above] {\textcolor{blue}{$\tau_2$}} node[midway, below] {$\tau_5$} (0.7,2.0);
\draw[<->] (n10b) -- node[midway, right] {$\tau_2$} node[midway, left] {\textcolor{blue}{$\tau_5$}} (n11t);
\draw[<->] (0.7,0.0) -- node[midway, above] {$\tau_5$} node[midway, below] {\textcolor{blue}{$\tau_2$}} (2.3,0.0);
\draw[<->] (3.7,0.0) -- node[midway, above] {$\tau_4$} node[midway, below] {\textcolor{blue}{$\tau_3$}} (5.3,0.0);
\draw[<->] (6.7,0.0) -- node[midway, above] {$\tau_3$} node[midway, below] {\textcolor{blue}{$\tau_4$}} (8.3,0.0);
\draw[<->] (9.7,0.0) -- node[midway, above] {$\tau_2$} node[midway, below] {\textcolor{blue}{$\tau_5$}} (11.3,0.0);
\draw[dotted, thick, <->] (n4b) -- node[midway, right] {$\tau_5^*$} node[midway, left] {\textcolor{blue}{$\tau_2^*$}} (n7t);
\draw[dotted, thick, <->] (n9b) -- node[midway, right] {$\tau_2^*$} node[midway, left] {\textcolor{blue}{$\tau_5^*$}} (n12t);
\draw[dotted, thick, <->] (n1b) -- node[midway, left] {$\tau_4^*$} node[midway, right] {\textcolor{blue}{$\tau_3^*$}} (n10t);
\draw[dotted, thick, <->] (n6b) -- node[midway, left] {$\tau_2^*$} node[midway, right] {\textcolor{blue}{$\tau_3^*$}} (n15t);
\end{tikzpicture}

\subsection{Computed PD Flowcharts }
\label{subsec:app_Computed PD Flowcharts}

Let $f$ be a generic function on the interval with five local extrema, at $\vec x = [1,2,3,4,5]^T$, with local extrema at levels $y_{(1)}<y_{(2)}<y_{(3)}<y_{(4)}<y_{(5)}$.  Then, the persistence diagram $\dgm(f)$ can be computed as follows: 

\noindent\makebox[\textwidth][c]{%
\begin{forest}
  for tree={
    parent anchor=south,
    child anchor=north,
    align=center,
    edge={-},
    draw, rectangle, rounded corners,
    font=\small,
    s sep=0.3mm,
    inner xsep=1pt,
    l sep+=14mm,
  }
  [{$(y_{(1)},y_{(5)}) \in \ess(f)$}, text width=3.1cm, minimum height=0.8cm
    [{$f(x_1) < f(x_2)$ \\ (local min.)}, text width=2.3cm, minimum height=0.8cm
      [$y_{(3)}$ local max., text width=2.25cm, minimum height=0.8cm, before computing xy={l-=10mm, s+=3.35mm}
        [{$(y_{(2)},y_{(3)}) \in \ord(f)$\\
        $(y_{(4)},y_{(5)}) \in \ord(f)$ \\
        $(y_{(3)}, f(x_3)) \in \rel(f)$}, text width=3.3cm, minimum height=0.8cm, before computing xy={l-=12mm}]
      ]
      [$y_{(4)}$ local max., text width=2.25cm, minimum height=0.8cm, l+=26mm
        [${(y_{(4)},f(x_3)) \in \rel(f)}$, text width=3.35cm, minimum height=0.8cm, s sep=8.8mm
          [{$f(x_1) = y_{(3)}$}, text width=2cm, minimum height=1cm, font=\small, before computing xy={l-=9.5mm}
            [{$(y_{(2)},f(x_4)) \in \ord(f)$ \\
            $(y_{(3)},f(x_2)) \in \ord(f)$}, text width=3.35cm, minimum height=1.1cm, font=\small, before computing xy={l+=10mm}]
          ]
          [{$f(x_3) = y_{(3)}$}, text width=2cm, minimum height=1cm, font=\small, l+=17mm, before computing xy={l-=13.5mm}
            [{$(y_{(2)},y_{(5)}) \in \ord(f)$ \\
            $(y_{(3)},y_{(4)}) \in \ord(f)$}, text width=3.1cm, minimum height=1.1cm, font=\small, before computing xy={l+=10mm}]
          ]
          [{$f(x_5) = y_{(3)}$}, text width=2cm, minimum height=1cm, font=\small, l+=34mm, before computing xy={l-=17.5mm}
            [{$(y_{(2)},f(x_2)) \in \ord(f)$ \\
            $(y_{(3)},f(x_4)) \in \ord(f)$}, text width=3.35cm, minimum height=1.1cm, font=\small, before computing xy={l+=10mm}]
          ]
        ]
      ]
    ]
    [{$f(x_1) > f(x_2)$ \\ (local max.)}, text width=2.3cm, minimum height=0.8cm
      [$y_{(2)}$ local min., text width=2.25cm, minimum height=0.8cm, l+=26mm
       [${(y_{(4)},f(x_3)) \in \rel(f)}$, text width=3.35cm, minimum height=0.8cm, s sep=8.8mm
          [{$f(x_1) = y_{(3)}$}, text width=2cm, minimum height=1cm, font=\small, before computing xy={l-=9.5mm}
            [{$(y_{(2)},f(x_4)) \in \ord(f)$ \\
            $(y_{(3)},f(x_2)) \in \ord(f)$}, text width=3.35cm, minimum height=1.1cm, font=\small, before computing xy={l+=10mm}]
          ]
          [{$f(x_3) = y_{(3)}$}, text width=2cm, minimum height=1cm, font=\small, l+=17mm, before computing xy={l-=13.5mm}
            [{$(y_{(2)},y_{(5)}) \in \ord(f)$ \\
            $(y_{(3)},y_{(4)}) \in \ord(f)$}, text width=3.1cm, minimum height=1.1cm, font=\small, before computing xy={l+=10mm}]
          ]
          [{$f(x_5) = y_{(3)}$}, text width=2cm, minimum height=1cm, font=\small, l+=34mm, before computing xy={l-=17.5mm}
            [{$(y_{(2)},f(x_2)) \in \ord(f)$ \\
            $(y_{(3)},f(x_4)) \in \ord(f)$}, text width=3.35cm, minimum height=1.1cm, font=\small, before computing xy={l+=10mm}]
          ]
        ]
      ]
      [$y_{(3)}$ local min., text width=2.25cm, minimum height=0.8cm, before computing xy={l-=10mm, s-=3.35mm}
         [{$(y_{(2)},y_{(1)}) \in \ord(f)$\\
        $(y_{(4)},y_{(3)}) \in \ord(f)$ \\
        $(y_{(3)}, f(x_3)) \in \rel(f)$}, text width=3.3cm, minimum height=0.8cm, before computing xy={l-=12mm}]
      ]
    ]
  ]
\end{forest}
}

\pagebreak
Conversely, let $f$ be a generic function on the interval with six local extrema, at $\vec x = [1,2,3,4,5,6]^T$, at levels $y_1<y_2<y_3<y_4<y_5<y_6$.  With $y_1$ and $y_6$ denoting the global minimum and maximum, respectively, we label the other (nonglobal) minima as $y_< < y_>$, and the other maxima $y^< < y^>$. 
Referring to Table \ref{app:P6_Flowchart_Legend} for specific $\ord$ and $\rel$ values, the persistence diagram $\dgm(f)$ can be computed with the following flowchart, split among the following two figures:
\\

\noindent\makebox[\textwidth][c]{%
\begin{tikzpicture}[
    node distance=10mm and 10.8mm,
    box/.style={
        draw,
        rectangle,
        minimum width=14mm,
        minimum height=10mm,
        align=center
    },
    every path/.style={thick}
]

\node[box] (A) {$(y_{(1)},y_{(6)}) \in \ess(f)$};

\node[box, below=10mm of A] (B) {$f(x_1) < f(x_2)$ \\ Local Minimum};

\draw (A)--(B);

\node[box, below left=15mm and 15.6mm of B] (B1) {$f(x_4) > f(x_2)$};
\node[box, below right=15mm and 15.6mm of B] (B2) {$f(x_4) < f(x_2)$};

\draw (B)--(B1);
\draw (B)--(B2);

\node[box, below left=15mm and 1.2mm of B1] (B11) {$f(x_5) = y_{(3)}$};
\node[box, below right=15mm and 1.2mm of B1] (B12) {$f(x_5) \neq y_{(3)}$};

\node[box, below left=15mm and 1.2mm of B2] (B21) {$f(x_1) = y_{(3)}$};
\node[box, below right=15mm and 1.2mm of B2] (B22) {$f(x_1) \neq y_{(3)}$};

\draw (B1)--(B11);
\draw (B1)--(B12);
\draw (B2)--(B21);
\draw (B2)--(B22);

\node[circle,draw,below=40mm of B1] (one) {1};
\node[circle,draw,below=40mm of B2] (two) {2};

\node[box, below=80mm of B] (D) {$f(\_) = y_{(4)}$};

\draw (B11) to (one);
\draw (B12) to (two);
\draw (B21) to  (one);
\draw (B22) to  (two);
\draw (one) to (D);
\draw (two) to (D);

\node[box, below left=20mm and 21.6mm of D] (D1) {$f(x_5) \  \rule[0ex]{1.2em}{0.4pt}\ f(x_3)$};
\node[box, below right=20mm and 21.6mm of D] (D2) {$f(x_5) \  \rule[0ex]{1.2em}{0.4pt} \ f(x_3)$};

\node[circle,draw,below=30mm of D] (four) {4};
\node[circle,draw,below=60mm of D] (three) {3};

\draw (D)--(D1) node[midway, left] {{\Large $x_2$ \ \ }};
\draw (D)--(four) node[midway, left] {{\Large $x_4$ \ }};
\draw (D)--(D2) node[midway, right] {{\ \ \Large $x_6$ }};
\draw (D1) -- (three) node[midway, above] {{\Large $\ \ >$}};
\draw (D1) -- (four) node[midway, above] {{\Large $\ <$}};
\draw (D2) -- (three) node[midway, above] {{\Large $ < \ \ $}};
\draw (D2) -- (four) node[midway, above] {{\Large $ > \ $}};

\end{tikzpicture}%
}

\pagebreak
\vspace*{-2cm}
\enlargethispage{2cm}

Second half: \\
\noindent\makebox[\textwidth][c]{%
\begin{tikzpicture}[
    node distance=10mm and 10.8mm,
    box/.style={
        draw,
        rectangle,
        minimum width=14mm,
        minimum height=10mm,
        align=center
    },
    every path/.style={thick}
]

\node[box] (A) {$(y_{(1)},y_{(6)}) \in \ess(f)$};

\node[box, below=10mm of A] (B) {$f(x_1) > f(x_2)$ \\ Local Maximum};

\draw (A)--(B);

\node[box, below left=15mm and 15.6mm of B] (B1) {$f(x_3) > f(x_5)$};
\node[box, below right=15mm and 15.6mm of B] (B2) {$f(x_3) < f(x_5)$};

\draw (B)--(B1);
\draw (B)--(B2);

\node[box, below left=15mm and 1.2mm of B1] (B11) {$f(x_2) = y_{(3)}$};
\node[box, below right=15mm and 1.2mm of B1] (B12) {$f(x_2) \neq y_{(3)}$};

\node[box, below left=15mm and 1.2mm of B2] (B21) {$f(x_6) = y_{(3)}$};
\node[box, below right=15mm and 1.2mm of B2] (B22) {$f(x_6) \neq y_{(3)}$};

\draw (B1)--(B11);
\draw (B1)--(B12);
\draw (B2)--(B21);
\draw (B2)--(B22);

\node[circle,draw,below=40mm of B1] (one) {5};
\node[circle,draw,below=40mm of B2] (two) {6};

\node[box, below=80mm of B] (D) {$f(\rule[0ex]{1.2em}{0.4pt}) = y_{(4)}$};

\draw (B11) to (one);
\draw (B12) to (two);
\draw (B21) to  (one);
\draw (B22) to  (two);
\draw (one) to (D);
\draw (two) to (D);

\node[box, below left=20mm and 21.6mm of D] (D1) {$f(x_2) \  \rule[0ex]{1.2em}{0.4pt} \ f(x_4)$};
\node[box, below right=20mm and 21.6mm of D] (D2) {$f(x_2) \  \rule[0ex]{1.2em}{0.4pt} \ f(x_4)$};

\node[circle,draw,below=30mm of D] (four) {8};
\node[circle,draw,below=60mm of D] (three) {7};

\draw (D)--(D1) node[midway, left] {{\Large $x_5$ \ \ }};
\draw (D)--(four) node[midway, left] {{\Large $x_3$ \ }};
\draw (D)--(D2) node[midway, right] {{\ \ \Large $x_1$ }};
\draw (D1) -- (three) node[midway, above] {{\Large $\ \ >$}};
\draw (D1) -- (four) node[midway, above] {{\Large $\ <$}};
\draw (D2) -- (three) node[midway, above] {{\Large $ < \ \ $}};
\draw (D2) -- (four) node[midway, above] {{\Large $ > \ $}};

\end{tikzpicture}%
}

\begin{table}[h]
\centering
\begin{tabular}{c p{5.5cm} c p{5.5cm}}
\hline
\textbf{Number} & \textbf{$\dgm(f)$ Contribution} & \textbf{Number} & \textbf{$\dgm(f)$ Contribution} \\
\hline
\circnum{1} & $(y_<,f(x_2)),(y_>,f(x_4)) \in \ord(f)$    & \circnum{5} & $(y_<,f(x_5)),(y_>,f(x_3)) \in \ord(f)$ \\[6pt]
\circnum{2} & $(y_<,f(x_4)),(y_>,f(x_2)) \in \ord(f)$  & \circnum{6} & $(y_<,f(x_3)),(y_>,f(x_5)) \in \ord(f)$  \\[6pt]
\circnum{3} & $(y^<,f(x_3)),(y^>,f(x_5)) \in \rel(f)$   & \circnum{7} & $(y^<,f(x_4)),(y^>,f(x_2)) \in \rel(f)$  \\[6pt]
\circnum{4} & $(y^<,f(x_5)),(y^>,f(x_3)) \in \rel(f)$     & \circnum{8} & $(y^<,f(x_2)),(y^>,f(x_4)) \in \rel(f)$   \\[6pt]
\hline
\end{tabular}
\caption{P6 Flowchart Legend. Note that $y_<$ and $y_>$ are the interior local minima, and $y^<$ and $y^>$ are the interior local maxima. }
\label{app:P6_Flowchart_Legend}
\end{table}
\end{document}